\documentclass[oneside, 12pt]{amsart}

\usepackage{amscd, amssymb, amsmath, mathrsfs, amsthm}
\usepackage[english]{babel}
\usepackage[all]{xy}
\usepackage{color}
\usepackage[colorlinks, linktocpage, citecolor = blue, linkcolor = blue]{hyperref}

\newcommand\mylabel[1]{\label{#1}}

\theoremstyle{plain}
\newtheorem{theorem}{Theorem}[section]
\newtheorem*{maintheorem}{Theorem}
\newtheorem{lemma}[theorem]{Lemma}
\newtheorem{proposition}[theorem]{Proposition}
\newtheorem{corollary}[theorem]{Corollary}

\theoremstyle{definition}
\newtheorem{definition}[theorem]{Definition}

\newtheorem*{acknowledgement}{Acknowledgement}

\theoremstyle{remark}

\DeclareFontFamily{U}{wncy}{}
\DeclareFontShape{U}{wncy}{m}{n}{<->wncyr10}{}
\DeclareSymbolFont{mcy}{U}{wncy}{m}{n}
\DeclareMathSymbol{\Sh}{\mathord}{mcy}{"58}

\newcommand{\ZZ}	{\mathbb{Z}}
\newcommand{\QQ}	{\mathbb{Q}}

\newcommand{\PP}	{\mathbb{P}}
\renewcommand{\AA}	{\mathbb{A}}
\newcommand{\GG}	{\mathbb{G}}

\newcommand{\ideala}    {\mathfrak{a}}

\newcommand  {\shA}     {\mathscr{A}}

\newcommand  {\shE}     {\mathscr{E}}
\newcommand  {\shF}     {\mathscr{F}}

\newcommand  {\shI}     {\mathscr{I}}

\newcommand  {\shM}     {\mathscr{M}}

\newcommand  {\shN}     {\mathscr{N}}
\newcommand  {\shL}     {\mathscr{L}}

\newcommand  {\alg}     {{\operatorname{alg}}}

\newcommand  {\APic}    {\operatorname{APic}}
\newcommand  {\Ass}     {\operatorname{Ass}}

\newcommand  {\Br}      {\operatorname{Br}}
\newcommand  {\Bs}      {\operatorname{Bs}}

\newcommand  {\can}     {{\rm \text{can}}}

\newcommand  {\Cl}      {\operatorname{Cl}}

\newcommand  {\Cokernel}{\operatorname{Coker}}

\newcommand  {\Div}     {\operatorname{Div}}

\newcommand  {\edim}    {\operatorname{edim}}

\newcommand  {\Ext}     {\operatorname{Ext}}

\newcommand  {\Frac}    {\operatorname{Frac}}

\newcommand  {\Gal}     {\operatorname{Gal}}

\newcommand  {\Hom}     {\operatorname{Hom}}

\renewcommand  {\k}     {\kappa}
\newcommand  {\Kernel } {\operatorname{Ker}}

\newcommand  {\length}	{\operatorname{length}}

\newcommand  {\loc}     {\text{\rm loc}}
\newcommand  {\dirlim}  {\varinjlim}

\newcommand  {\lra}     {\longrightarrow}

\newcommand  {\Mat}     {\operatorname{Mat}}

\newcommand  {\maxid}   {\mathfrak{m}}

\newcommand  {\NS}      {\operatorname{NS}}

\renewcommand{\O}       {\mathscr{O}}

\newcommand  {\perf}     {\text{{\rm perf}}}

\newcommand  {\pd}    	{\operatorname{pd}}
\newcommand  {\pdeg}    {\operatorname{pdeg}}

\newcommand  {\Pic}     {\operatorname{Pic}}

\newcommand  {\pr}      {\operatorname{pr}}
\newcommand  {\Proj}    {\operatorname{Proj}}

\newcommand  {\quadand} {\quad\text{and}\quad}

\newcommand  {\ra}      {\rightarrow}

\newcommand  {\rank}    {\operatorname{rank}}
\newcommand  {\red}     {{\operatorname{red}}}
\newcommand  {\Reg}     {\operatorname{Reg}}
\newcommand  {\res}     {\operatorname{res}}

\newcommand  {\sep}     {{\operatorname{sep}}}

\newcommand  {\Sing}    {\operatorname{Sing}}

\newcommand  {\Spec}    {\operatorname{Spec}}

\newcommand  {\Supp}    {\operatorname{Supp}}
\newcommand  {\Sym}     {\operatorname{Sym}}

\newcommand  {\trdeg}   {\operatorname{trdeg}}

\newcommand{\tF}{{\tilde{F}}}

\newcommand{\tD}{{\tilde{D}}}
\newcommand{\hD}{{\hat{D}}}

\def\mydate{\number\day\space\ifcase\month \or January\or February\or March\or 
April\or May\or June\or July\or
August\or September\or October\or November\or December\fi \space\number\year}

\DeclareFontFamily{U}{wncy}{}
\DeclareFontShape{U}{wncy}{m}{n}{<->wncyr10}{}
\DeclareSymbolFont{mcy}{U}{wncy}{m}{n}
\DeclareMathSymbol{\Sh}{\mathord}{mcy}{"58}

\begin{document}

\title[Del Pezzo surfaces in positive characteristic]
      {Del Pezzo surfaces and Mori fiber spaces in positive characteristic}

\author[Andrea Fanelli]{Andrea Fanelli}
\address{Institut de Math\'ematiques de Bordeaux, CNRS UMR 5251, Universit\'e de Bordeaux, 33405 Talence CEDEX, France }
\curraddr{}
\email{andrea.fanelli@uvsq.fr}

\author[Stefan Schr\"oer]{Stefan Schr\"oer}
\address{Mathematisches Institut, Heinrich-Heine-Universit\"at, 40204 D\"usseldorf, Germany}
\curraddr{}
\email{schroeer@math.uni-duesseldorf.de}

\subjclass[2010]{14E30, 14G17, 14J26, 14J30,  14J17,  14H20, 14J45}

\dedicatory{Final version, 8 September 2019}

\begin{abstract}
We settle a question that originates from  results and remarks by  Koll\'ar on extremal ray in the minimal model program:
In positive characteristics, there are no Mori fibrations on    threefolds with only terminal singularities    
whose   generic fibers are      geometrically non-normal surfaces.
To show this we establish some general structure results for
del Pezzo surfaces over imperfect ground fields.
This relies on  Reid's classification of non-normal del Pezzo surfaces
over algebraically closed fields, combined with a detailed  analysis of geometrical non-reducedness
over imperfect fields of $p$-degree one. 
\end{abstract}

\maketitle
\tableofcontents

\section*{Introduction}
\mylabel{Introduction}

A smooth proper scheme  $V$ of dimension two
whose dualizing sheaf $\omega_V$ is antiample is called a \emph{del Pezzo surface}.
This notion immediately generalizes from the smooth case to the Gorenstein case.
Smooth del Pezzo surfaces were first studied
by del Pezzo \cite{delPezzo 1887} in the nineteenth century. Over algebraically closed ground fields, they are either   $\PP^1\times\PP^1$ 
or the blowing-up of   $\PP^2$ in at most eight points in general position
(see \cite{Demazure 1976} or \cite{Dolgachev 2012} for this classification in modern language).
In some sense, del Pezzo surfaces are the two-dimensional analogs of the projective line.
The higher-dimensional generalizations are called \emph{Fano varieties}.

Del Pezzo surface and Fano varieties play an important role in the
\emph{minimal model program}, which is a tremendously successful approach
to achieve a classification of higher-dimensional algebraic schemes over
algebraically closed ground fields $k$.
If $f:X\ra B$ is a \emph{Mori fibration}, that is, the  contraction of an extremal ray of fiber type,
then the generic  fiber $V=X_\eta$ is a Fano variety over the function field
$$
F=\O_{B,\eta}=\kappa(\eta)=k(B).
$$
Furthermore, $\dim(V)=\dim(X)-\dim(B)$ 
and the  Picard number is $\rho(V)=1$. In particular,  Mori fibrations
on \emph{threefolds}  $X$ yield   del Pezzo surfaces $V $ over   function fields $F$ of   algebraic curves $B$.
Of course, one may also study more general fibrations $f:X\ra B$ having del Pezzo surfaces or Fano varieties
as generic fiber. One may call them \emph{del Pezzo fibrations} or \emph{Fano fibrations}.
In any case, to understand the geometry of $X$, it is imperative to understand the geometry
of $V=X_\eta$ over the non-closed field $F$.

If the total space $X$ is smooth, that is, the sheaf of K\"ahler differentials $\Omega^1_{X/k}$ is locally free of
rank $\dim(X)$,   the generic fiber $V=X_\eta$ is regular, in the sense that all local rings $\O_{V,a}$
are regular. In characteristic zero, this ensures that $V$ is smooth as well,
and one may understand it in terms of  the base-change $V\otimes_F F^\alg$, together with the action of the Galois group
$\Gal(F^\alg/F)$. This was exploited, for example, in \cite{Mori 1982}, \cite{{Codogni et al 2016}} and \cite{{Codogni et al 2018}}.

In positive characteristics $p>0$, however, this is no longer true, and it may easily happen
that the \emph{geometric generic fiber}  $V\otimes_FF^\alg$ is non-normal, or even non-reduced.
The former already plays an important role in the Enriques classification of surfaces: In characteristic $p=2$
and $p=3$, there are \emph{quasielliptic fibrations}, where the generic fiber is a regular genus-one curve
and the geometric generic fiber is the rational cuspidal curve \cite{Bombieri; Mumford 1976}. The latter easily happens  in higher-dimensions:
The hypersurface $X\subset\PP^n\times\PP^n$ given by the bihomogeneous equation
$$
X:\quad S_0T_0^p+S_1T_1^p+\ldots+ S_nT_n^p=0
$$
is smooth, whereas the geometric generic fiber for the projection $\pr_1:X\ra\PP^n$ is a $p$-fold hyperplane.

Del Pezzo  surfaces in positive characteristics $p>0$ and their log-generalizations have been studied, among others,
by Reid \cite{Reid 1994}, Schr\"oer \cite{Schroeer 2007}, Maddock \cite{Maddock 2016}, Cascini, Tanaka \cite{Cascini; Tanaka 2017}, 
Cascini, Tanaka and Witaszek \cite{Cascini; Tanaka; Witaszek 2017}, Bernasconi \cite{Bernasconi 2017}, and Das \cite{Das 2017}.
The minimal model program was introduced by Mori \cite{Mori 1982}, and originally focused on the situation 
over the complex numbers, although Shepherd-Barron analyzed Fano threefolds in positive characteristics \cite{Shepherd-Barron 1997}.
Recently, MMP  made  tremendous advances in positive characteristic,
for example by the work of Hacon and Xu \cite{Hacon; Xu 2015}  and Tanaka \cite{Tanaka 2015}.

However, many foundational issues remained open.
About 25 years ago,  Koll\'ar's analysis of extremal rays on threefolds
raised the question whether geometric non-normality may appear on Mori fibrations on threefolds
(\cite{Kollar 1991}, Remarks in 1.2), as it does in the Enriques classification of surfaces. 
Our main result  is that this--perhaps surprisingly--does not happen:

\begin{maintheorem} 
{\rm (See   \ref{main fibration}.)}
Suppose $k$ is an algebraically closed ground field of characteristic $p>0$.
Let $X$ be a threefold with only terminal singularities, and $f:X\ra B$ be a Mori fibration 
of relative dimension two. Then the generic fiber $V=X_\eta$ is geometrically normal.
\end{maintheorem}

This generalizes results of Saito \cite{Saito 2003}, who treated the case where $X$ is a Fano threefold
with $\rho(X)=2$, and Patakfalvi and Waldron \cite{Patakfalvi; Waldron 2017}, who ruled
out the cases $p\geq 5$.
Besides the above non-existence result  concerning Mori fibrations,  also show that     del Pezzo fibration with unusual properties
actually do exist:

\begin{maintheorem}
{\rm (See   \ref{main examples}.)}
Let $k$ be an algebraically closed ground field of characteristic $p=2$.
\begin{enumerate}
\item There is a Mori fibration $Y\ra\PP^1$ whose generic fiber is a non-smooth del Pezzo surface
that is geometrically normal.
\item There is a del Pezzo fibration $X\ra \PP^1$ whose generic fiber is geometrically non-normal with
Picard number two.
\end{enumerate}
Here $X$ arises from  $Y$ by some blowing-up whose center $Z\subset Y$ is a horizontal curve.
\end{maintheorem}

Situations like in (i) are probably well-known, and we list it here because it leads to case (ii).

The preceding results are special cases of our analysis of del Pezzo surfaces $V$ over   arbitrary imperfect
fields $F$, by relating the geometry of the surface to the arithmetic of the ground
field. In fact, we develop several  techniques that work for general algebraic schemes over imperfect ground fields,
which should be useful in many other situations.

The key idea is to use   the \emph{$p$-degree} $\pdeg(F)\geq 0$ of the ground field systematically.
This notion was introduced by Teichm\"uller \cite{Teichmueller 1936} under the name \emph{degree of imperfection}, and was further studied by
Becker and MacLane \cite{Becker; MacLane 1940}. It can be seen as
the dimension of the vector space of absolute K\"ahler differentials $\Omega^1_{F/\ZZ}$.
If $F$ is the function field of some integral algebraic scheme $B$ over a perfect field, we just have
$\pdeg(F)=\dim(B)$. The non-existence result concerning Mori fibrations on threefolds comes from the following statement:

\begin{maintheorem}
{\rm (See   \ref{main pdeg}.)}
Let $V$ be a regular del Pezzo surface with Picard number one over a ground field $F$ of   $p$-degree one.
Then $V$ is geometrically normal.
\end{maintheorem}

As a rule of thumb, the higher the $p$-degree, the more unusual the geometry of algebraic $F$-schemes may become.
Some scholars  may regard the  ensuing possibilities as ``pathological'' or ``psychedelic'', but we believe that
these effects are rather natural and deserve systematic further  study.

The   crucial idea for   our results is to study the \emph{locus of non-smoothness} $\Sing(V/F)$ for
regular but geometrically non-normal del Pezzo surfaces. It  carries a natural
scheme structure via Fitting ideals, and we look at its   divisorial part $N\subset\Sing(V/F)$ and its reduction $D=N_\red$. 
This is an effective  Cartier divisor $D\subset V$,
and we analyze how it could fit in with \emph{Reid's Classification} \cite{Reid 1994} of non-normal del Pezzo surface $Y=V_K$ 
obtained by base-changing along sufficiently large finite purely inseparable field extensions $F\subset K$.
This non-normal del Pezzo surface will be analyzed in terms of its \emph{conductor square}
$$
\begin{CD}
R	@>>>	X\\
@VVV		@VV\nu V\\
C	@>>>	Y,
\end{CD}
$$
where $\nu:X\ra Y$ is the normalization, $R\subset X$ is the ramification divisor and $C\subset Y$ is the conductor curve. The latter
is not a Cartier divisor, but has the same support of the Cartier divisor $D_K\subset Y$.

Another important input comes from Serre's \emph{characterization of one-dimension Gorenstein rings}
in terms of length conditions \cite{Serre 1988}, which was extended by Reid in the language of schemes \cite{Reid 1994}.
These Gorenstein conditions can be used to obtain information about the local rings of $\O_{Y,y}$
and the structure of $N$ and its reduction $D=N_\red$. 
Finally, we combine the geometry  of our schemes with the arithmetic of   ground fields by
introducing the \emph{geometric generic embedding dimension}
$$
\edim(\O_{N,\eta}/F)=\edim(\O_{N\otimes_FF^\perf})
$$
for arbitrary integral algebraic schemes $N$. Another  key observation is the following,
which strengthens some general bound of the second author \cite{Schroeer 2010}:

\begin{maintheorem}
{\rm (See   \ref{edim bound})}
Let $F$ be a field with $\pdeg(F)\leq 1$. Then for each proper integral scheme $N$,
we have $\edim(\O_{N,\eta}/F)\leq 1$.
\end{maintheorem}

Another important ingredient is \emph{Maddock's bound} \cite{Maddock 2016}: If $V$
is a normal del Pezzo surface in characteristic $p>0$ with  irregularity $h^1(\O_V)>0$, then
this irregularity is actually bounded from below by 
$h^1(\O_V) \geq \frac{p^2-1}{6} K_V^2$.

\medskip
The paper is organized as follows: 
In the first two sections we   establish various foundational facts about 
generic geometric embedding dimension $\edim(\O_{N,\eta}/F)$ and
the locus of non-smoothness $\Sing(V/F)$.
Section \ref{Local computations} contains some   computations with complete local rings for conductor squares, mainly
in codimension one and two, which gives some   information on the behavior of the locus of non-smoothness.
In Section \ref{Geometrically non-normal} we collect some general facts on proper schemes that
are geometrically non-normal, and give some results on their Picard schemes.
This is applied in Section \ref{Regular del Pezzo} to del Pezzo surfaces that are geometrically non-normal, where we also tabulate
the possibilities according to Reid's classification.
Section \ref{Smooth ramification} treats the case that the ramification divisor $R\subset X$
is smooth: It turns out that only one case of Reid's classification is possible, which has Picard number $\rho(V)=2$.
The cases that the ramification divisor is non-smooth is much more challenging.
Here we start with a preliminary investigation, excluding further possibilities in 
Section \ref{Non-reduced ramification p=2} for characteristic two, and Section \ref{Non-reduced ramification p=3} for characteristic three.
In the remaining cases, the reduced locus of non-smoothness $D=\Sing(V/F)$ are 
curves of low genera with rather peculiar properties.
We investigate the structure of such curves in Sections \ref{Peculiar curves} and \ref{Peculiar higher genus}, where non-uniqueness of
coefficient fields plays a decisive role. 
This is applied in Sections \ref{Non-existence} and \ref{Non-existence irregularity p=2}, where we rule out the remaining cases
in characteristic two. Section \ref{Non-existence p=3} contains a similar analysis for characteristic three.
The paper closes with an Appendix, where we give a general treatment of conductor squares
and Gorenstein conditions, with some slight generalizations of results of Serre, Reid
and many others.

\begin{acknowledgement}
We like  to  thank Paolo Cascini, Ariyan Javanpeykar, J\'anos Koll\'ar,  Zsolt Patakfalvi, Hiromu Tanaka and Joe Waldron for 
useful comments. We also wish to thank the referee for very careful reading,   many valuable suggestions and 
for pointing out some mistakes.
The first-named author was funded by the Deutsche Forschungsgemeinschaft
with the grant PE 2165/1-2 \emph{Gromov-Witten Theorie, Geometrie und Darstellungen}.
This research was conducted in the framework of the   research training group
\emph{GRK 2240: Algebro-geometric Methods in Algebra, Arithmetic and Topology}, which is funded
by the Deutsche Forschungsgemeinschaft. 
\end{acknowledgement}

\section{Geometric generic embedding dimension and  \texorpdfstring{$p$}{}-degree}
\mylabel{Geometric generic}

We start by collecting  some notions and results pertaining to 
local rings at generic points for algebraic schemes.
Let $A$ be a local Artin ring, with maximal ideal $\maxid_A\subset A$
and residue field $\kappa=A/\maxid_A$.
The two basic numerical invariants are   \emph{embedding dimension} and   
\emph{Hilbert--Samuel multiplicity}
$$
\edim(A) = \dim_{\kappa} \maxid_A/\maxid_A^2\geq 0\quadand
e(A) = \length(A)\geq 1.
$$
For the general theory of Hilbert--Samuel multiplicities, we refer to 
\cite{AC 8-9}, Chapter VIII, \S7.
The following immediate observations show  that 
both integers give a     measure for   non-regularity:

\begin{proposition}
\mylabel{artin rings}
The following are equivalent:
\begin{enumerate}
\item The local Artin ring $A$ is regular.
\item The projection $A\ra\kappa$ is bijective.
\item The embedding dimension is  $\edim(A)=0$.
\item The Hilbert--Samuel multiplicity is $e(A)=1$.
\end{enumerate}
\end{proposition}

Now suppose $F$ is a ground field of characteristic $p>0$,
and let $A$ be local ring that  is essentially of finite type as $F$-algebra,
with Krull dimension $\dim(A)=0$. These are precisely the local Artin rings where the
field extension $F\subset \kappa(A)$
is   finitely generated.
They can also be regarded as     stalks $\O_{Y,\eta}$ at the  generic point  $\eta$ of   irreducible  schemes $Y$
that are separated and of finite type. If the latter holds, we also
say that $Y$ is an \emph{algebraic scheme}. 

For each algebraic field extension $F\subset K$, the ring $A_K=A\otimes_FK$ is essentially of finite type over $K$,
hence noetherian, and integral over $A$ and therefore $\dim(A_K)=0$.
It follows that $A_K$ is an Artin ring.
If $F\subset K$ is purely inseparable, the rings stays local,
and we can consider the integers $\edim(A_K)$ and $e(A_K)$ as above.
Of particular interest is  the situation where $K=F^\perf$ is the \emph{perfect closure}.
Let us call the integer
$$
\edim(A/F) = \edim(A\otimes_FF^\perf)\geq 0
$$ 
the \emph{geometric embedding dimension}.
Similarly,   define  
$$
e(A/F) = e(A\otimes_FF^\perf)\geq 1
$$
as the  \emph{geometric Hilbert--Samuel multiplicity}. Furthermore, 
if $Y$ is an irreducible algebraic scheme  with generic point $\eta\in Y$, we call the integers
$$
\edim(\O_{Y,\eta}/F) = \edim(\O_{Y,\eta}\otimes_FF^\perf)\quadand e(\O_{Y,\eta}/F) = e(\O_{Y,\eta}\otimes_FF^\perf).
$$
the \emph{geometric generic embedding dimension} and the
\emph{geometric generic Hilbert--Samuel multiplicity} of the scheme $Y$.
We thus get:

\begin{proposition}
\mylabel{geometrically reduced}
For an irreducible algebraic scheme $Y$ without embedded components, the following are equivalent:
\begin{enumerate}
\item The scheme $Y$ is geometrically reduced.
\item The geometric generic embedding dimension is  $\edim(\O_{Y,\eta}/F) =0$.
\item The geometric generic Hilbert--Samuel multiplicity is $e(\O_{Y,\eta}/F)=1$.
\end{enumerate}
\end{proposition}

\proof
Conditions (ii) and (iii) are equivalent by Proposition \ref{artin rings}.
The assumption that $Y$ has no embedded components means that the scheme satisfies
Serre's Condition $(S_1)$. This also holds for each base-change $Y_K=Y\otimes_FK$,
according to \cite{EGA IVb}, Proposition 6.7.1. It follows that $Y$ is geometrically reduced if and only
if the local Artin ring $A=\O_{Y,\eta}$ is geometrically reduced.

Each   of the Conditions (i)--(iii) imply that the scheme $Y$ and the local Artin ring $A$ are reduced,
so it suffices to treat this situation.
Then the field extension  $F\subset A$ is geometrically reduced if and only if it 
is separable. The latter holds if and only if $A\otimes_FF^\perf$ is reduced,
according to MacLane's Criterion \cite{A 4-7}, Chapter V, \S2, No.\ 4, Theorem 2.
\qed

\medskip
The following notion goes back to Teichm\"uller \cite{Teichmueller 1936}, under the name \emph{degree of imperfection},
and was further studied
by Becker and MacLane \cite{Becker; MacLane 1940}.

\begin{definition}
\mylabel{p-degree}
The dimension of the vector space of absolute K\"ahler differentials $\Omega^1_{F/\ZZ}=\Omega^1_{F/F^p}$ is called the
  \emph{$p$-degree} $\pdeg(F)$ of the field $F$.
\end{definition}

It could be seen as the cardinality of a $p$-basis for $F\subset F^{1/p}$,
or $F^p\subset F$. If these algebraic extensions are finite, the degree is of the form
$[F^{1/p}:F]=p^n$ with exponent $n=\pdeg(F)$.

The second author established in \cite{Schroeer 2010}, Theorem 2.3 
the relation 
$$
\edim(\O_{Y,\eta}/F)<\pdeg(F)
$$
for every proper normal scheme $Y$ with $h^0(\O_Y)=1$ that is not geometrically reduced.
As a consequence, if the ground field has   $\pdeg(F)\leq 1$, then
every proper normal scheme $Y$ with $h^0(\O_Y)=1$
is geometrically reduced. The following  extension to \emph{arbitrary} proper integral schemes,
which could be non-normal or could  have $h^0(\O_Y)\geq 2$,
will be crucial for our applications:

\begin{theorem}
\mylabel{edim bound}
Suppose the ground field $F$ has $p$-degree $\pdeg(F)\leq 1$.
Then for each proper integral scheme $Y$, we have $\edim(\O_{Y,\eta}/F)\leq 1$.
\end{theorem}

\proof
The conclusion is trivial if $Y$ is geometrically reduced.
So we assume that $Y$ is not geometrically reduced.
Let $X\ra Y$ be the normalization. Then $X$ 
is a proper scheme that is integral but not geometrically reduced.
The ring of global sections $K=H^0(X,\O_X)$
is integral with $[K:F]<\infty$, whence $F\subset K$ is a finite field extension.
According to \cite{Schroeer 2010}, Theorem 2.3 we have $\edim(\O_{X,\eta}/K)<\pdeg(K)$.
But $\pdeg(K)=\pdeg(F)=1$ according to \cite{Becker; MacLane 1940}, Theorem 3.
So by Proposition \ref{geometrically reduced}, the algebraic scheme  $X$ is geometrically reduced
over the field $K$. 

Now choose a perfect closure $E=K^\perf$.
The composite extension $F\subset K\subset E$ is a  perfect closure for  $K$.
Write  $A=\O_{X,\eta}=\O_{Y,\eta}$ for the common function field of the integral schemes $X$ and $Y$. Then 
$$
A\otimes_F E = A\otimes_K (K\otimes_F E). 
$$
Since $\pdeg(F)\leq 1$, the finite field extension $F\subset K$ is obtained by
adjoining a single element $\alpha\in K$, such that  $K=F(\alpha)$, according to   \cite{Becker; MacLane 1940}, Theorem 1. 
Let $f\in F[T]$ be the minimal polynomial
of this generator $\alpha\in K$. Then $K\otimes_FE=E[T]/(f)$, thus 
$$
A\otimes_FE=A\otimes_K E[T]/(f)= (A\otimes_KE)[T]/(f).
$$
Since $A\otimes_KE$ is a field, it 
follows that the residue class of the indeterminate $T$ generates the maximal ideal of the local Artin ring $A\otimes_FE$,
and consequently  $\edim(\O_{Y,\eta}/F)\leq 1$.
\qed

\section{Singular loci on algebraic schemes}
\mylabel{Singular loci}

In this section we establish some useful facts on closed subschemes contained in singular loci,
and introduce some general notations along the way.
Let $F$ be a ground field and   $V$ be an algebraic $F$-scheme, which means a separated scheme of finite type.
As customary, we write $\Sing(V)$ for the set of points $a\in V$ where the local ring $\O_{V,a}$ is
not regular, and call it the \emph{singular locus}.
Note that $\Sing(V)\subset V$ is a closed set, according to \cite{EGA IVb}, Corollary 6.12.5.
One may regard it as a closed subscheme, endowed with   reduced scheme structure.
For every field extension $F\subset K$, we have $\Sing(V)\otimes_FK\subset\Sing(V\otimes_FK)$
by loc.\ cit.\  Proposition 6.5.1. Note that this inclusion is an equality provided that $F\subset K$ is separable,
but in general one has a strict inclusion. The following relative version is a remedy for this defect:

\begin{definition}
\mylabel{locus non-smoothness}
The \emph{locus of non-smoothness} $\Sing(V/F)\subset V$ is the set of points $a\in V$ where the local ring $\O_{V,a}$ is not geometrically regular
as $F$-algebra.
\end{definition}

According to \cite{EGA IVb}, Definition 6.7.6, this
means that for some finite field extension $F\subset K$, the resulting semilocal noetherian ring $\O_{V,a}\otimes_FK$ becomes non-regular.
By loc.\ cit.\ Proposition 6.7.7 this already appears  for some purely inseparable
field extension $F\subset K$. The relative and absolute notions are related as follows:

\begin{proposition}
\mylabel{locus as image}
Let $K=F^\perf=F^{1/p^\infty}$ be the perfect closure. Then $\Sing(V/F)$ is the image of $\Sing(V\otimes_FK)$
under the universal homeomorphism $V\otimes_FK\ra V$.
\end{proposition}

\proof
By \cite{EGA IVb}, Corollary 6.7.8, the locus of non-smoothness
$\Sing(V/F)$ coincides with the image of $\Sing(V\otimes_FK/K)$ under the     projection $V\otimes_FK\ra V$. 
But over a perfect fields, the singular locus coincides with the locus of non-smoothness.
\qed

\medskip
In particular, we see that $\Sing(V/F)\subset V$ is a closed set, and it commutes with ground field extensions.
In contrast to the singular locus, it comes with a canonical scheme structure
that is usually non-reduced. This relies on Fitting ideals for coherent sheaves $\shF$.
If $V=\Spec(R)$ is affine, we may choose a finite presentation
\begin{equation}
\label{fitting presentation}
R^{\oplus s}\stackrel{A}{\lra} R^{\oplus r}\lra M \lra 0.
\end{equation}
for the $R$-module $M=\Gamma(V,\shF)$, 
with some  matrix $A\in\Mat_{r\times s}(R)$. For each integer $0\leq n\leq r$, the ideal $\ideala_n\subset R$ generated by
the $(r-n)$-minors of the matrix $A$ is called the \emph{$n$-th Fitting ideal} of the module $M$.
For $n>r$, we set $\ideala_n=R$.
Indeed, this depends only on the module and not on the chosen presentation, confer the discussion
in \cite{Eisenbud 1995}, Section 20.2.
By gluing, we thus get Fitting ideals as coherent ideal sheaves $\shI_n\subset\O_V$ for           $\shF$ on 
arbitrary algebraic schemes $V$.
The corresponding closed set is  the locus of   points $a\in V$ where some presentation \eqref{fitting presentation} 
with $n\leq r$  exists on some affine open neighborhood such that 
the matrix $A\otimes\kappa(a)$  has rank $<r-n$,
in other words  $\dim_{\kappa(a)}(M\otimes_R\kappa(a))>n$.
According to  Nakayama's Lemma, the latter means that $\shF$ needs at least $n+1$ generators on each 
affine open neighborhood of the point $a$.

\begin{proposition}
\mylabel{scheme structure}
Suppose that the algebraic scheme $V$ is equidimensional   of dimension $n=\dim(V)$.
Then the closed set defined by the $n$-th Fitting ideal for the sheaf of K\"ahler differentials $\Omega^1_{V/F}$ coincides with
the locus of non-smoothness $\Sing(V/F)$. 
\end{proposition}

\proof
Let $U\subset V$ be the complementary open set.
As discussed above, this is the set of all points $a\in V$ where the stalk $\Omega^1_{V/F,a}$ can be generated by $n$ elements.

Likewise, let $U'\subset V$ be the complement of the locus of non-smoothness.
Recall that the \emph{local dimension}     $\dim_a(V)$ is the limit of the dimensions of 
open neighborhoods of $a\in V$. It takes  constant value $n$  because $V$ is equidimensional, and 
coincides with the \emph{relative dimension} $\dim_a(f)=\dim_a f^{-1} f(a) $ 
of  the structure morphism $f:V\ra\Spec(F)$.
According to \cite{EGA IVd}, Proposition 17.15.15 the open set $U'$ is the set of points where $\Omega^1_{V/F}$ is locally
free of rank $n=\dim_a(f)$. 

We thus have $U'\subset U$. Seeking a contradiction, we assume that the inclusion is strict.
Making a base-change, we may assume that $F$ is algebraically closed.
By Hilbert's Nullstellensatz there is a rational point $a\in U\smallsetminus U'$, and thus we get
$\Omega^1_{V/F}\otimes\kappa(a) = \maxid_a/\maxid_a^2$. The former vector space has dimension $\leq n$,
because $a\in U$. The latter vector space has dimension $\edim(\O_{V,a})> \dim(\O_{V,a})=n$, because
$a\not\in U'$, contradiction.
\qed
 
\medskip
For equidimensional algebraic schemes $V$, we usually regard the locus of non-smoothness $\Sing(V/F)$ as  a closed subscheme,
endowed with the scheme structure coming from the $n$-th Fitting ideal for the coherent sheaf  $\Omega_{V/F}^1$, where $n=\dim(V)$.
This subscheme is stable under ground field extension, and has the following 
strange property:

\begin{proposition}
\mylabel{reduced subschemes geom non-reduced}
Let  $Z\subset\Sing(V/F)$ be a reduced closed subscheme, and $\eta\in Z$ be a generic point.
Suppose that $\O_{Z,\eta}=\O_{V,\eta}/(f_1,\ldots,f_r)$ for some regular sequence
$f_1,\ldots,f_r\in\O_{V,\eta}$. Then the scheme $Z$ is geometrically non-reduced.
\end{proposition}

\proof
Suppose that $Z$ is geometrically reduced. Base-changing to the algebraic closure,
we may assume that $F=F^\alg$. Then $\eta\in\Sing(V/F)=\Sing(V)$, such that  the local ring $\O_{V,\eta}$ becomes singular.
On the other hand, the   local Artin ring $\O_{Z,\eta}$ is regular.
Since the sequence $f_1,\ldots, f_r$ is regular, the local $\O_{V,\eta}$
must be regular, contradiction.
\qed

\medskip
Since regular subschemes in regular schemes are locally given by regular sequences,
we obtain: 

\begin{corollary}
\mylabel{regular subschemes geom non-reduced}
Suppose that $V$ is regular, and let $Z\subset V$ be some regular subscheme contained in $\Sing(V/F)$.
Then $Z$ is geometrically non-reduced.
\end{corollary}

\medskip
The following special case, which already appears in   \cite{EGA IVd}, Proposition 17.15.1,   will play an important role throughout:

\begin{corollary}
\mylabel{closed points geom non-reduced}
Let $a\in V$ be a closed point contained in   $\Sing(V/F)$
 whose local ring $\O_{V,a}$ is regular. 
Then the finite field extension $F\subset\kappa(a)$ is not separable. In particular, the closed point
$a\in V$ is not rational.
\end{corollary}

Since prime divisors in normal schemes are generically defined by a single equation,
we also have the following consequence:

\begin{corollary}
\mylabel{prime divisors geom non-reduced}
Suppose that $V$ is normal, and let $Z\subset V$ be some prime divisor contained in 
$\Sing(V/F)$. Then $Z$ is geometrically non-reduced.
\end{corollary}

Assume now that $V$ is geometrically integral. Let $\shI\subset\O_V$ be the Fitting ideal
corresponding to the locus of non-smoothness.
There are only finitely many points $a_1,\ldots,a_r\in V$ of codimension one with $\shI_{a_i}\neq 0$.
They admit  a common affine open neighborhood $U=\Spec(R_0)$, even if there is no ample sheaf 
(\cite{Gross 2012},  proof for Theorem 1.5).
Let  $R=S^{-1}R_0$ be the ensuing semilocal ring,
with Fitting ideal $\ideala\subset R$. The schematic image of the morphism $\Spec(R/\ideala)\ra V$
is called the \emph{divisorial part} $N\subset \Sing(V/F)$ for the locus of non-smoothness.
If non-empty, the subscheme $N\subset V$ is purely one-codimensional and without embedded component.
Moreover, its formation commutes with ground field extensions.

Throughout the paper, one main idea is to analyze the divisorial part $N\subset\Sing(V/F)$ and
the resulting reduction $D=N_\red$, which is also purely one-co\-dimensional and without embedded component.
Note, however, that  $D=N_\red$ does not necessarily commute with inseparable ground field extensions.
This phenomenon will play a crucial role.

If $V$ is moreover normal, such that the one-dimensional local rings $\O_{V,a_i}$  are discrete valuation rings,
the Fitting ideals take the form $\shI_{a_i}=\maxid_{a_i}^{m_i}$, and 
one may regard the subscheme $N\subset V$ as the  effective \emph{Weil divisor} $N=\sum m_i D_i$, where $D_i=\overline{\{a_i\}}$.
Its reduction is $D=\sum D_i$.
If furthermore $V$ is locally factorial, the Weil divisors $N$ and $D=N_\red$ are Cartier divisors.
For each base-change $Y=V\otimes_FK$, we get an induced \emph{Cartier divisor} $D\otimes_FK$,
supported by the locus of non-smoothness.

\section{Local computations}
\mylabel{Local computations}

Fix a ground field $K$ of characteristic $p>0$. Let $X$ and $Y$ be algebraic schemes,
and $\nu:X\ra Y$ be a finite modification. We then have the \emph{conductor square}
\begin{equation}
\label{conductor square local}
\begin{CD}
R	@>>>	X\\
@VVV		@VV\nu V\\
C	@>>>	Y,
\end{CD}
\end{equation}
where $R\subset X$ denotes the \emph{ramification locus}, $C\subset Y$ is the \emph{conductor scheme}, and $\nu:R\ra C$ is the \emph{gluing map},
as discussed in detail in   Appendix \ref{Finite modifications}. 
Throughout, we also assume that $X$ is normal, that $Y$ satisfies conditions $(S_2)$ and $(G_1)$,
and that $\nu:X\ra Y$ is a universal homeomorphism. Here $(G_1)$ means that the local rings $\O_{Y,y}$ are Gorenstein
for all points $y\in Y$ of codimension one.

The goal of this section is to study the complete local rings $\O_{Y,y}^\wedge$ 
at points $y\in Y$ of codimension $\leq 2$ contained in the conductor scheme $C$,
and the schematic structure of the locus of non-smoothness $N=\Sing(Y/K)$.
To avoid   technical problems with respect to
K\"ahler differentials of complete local rings, we 
also assume that the ground field $K$ has finite $p$-degree,
compare the discussion in \cite{Schroeer 2017}, Section 1.
Let us  start with the case that the ramification locus is regular.

\begin{proposition}
\mylabel{ramification regular}
If the ramification locus $R\subset X$ is regular, then the conductor scheme $C\subset Y$ is normal,
and    for each point $x\in X$ of codimension one whose image $y\in Y$ is 
contained in  $C$, the field extension $\O_{C,y}\subset\O_{R,x}$
is purely inseparable of degree two. In particular, the characteristic must be $p=2$.
\end{proposition}

\proof
By assumption, we have $\O_{R,x}=\kappa(x)$.
It follows that the   local Artin ring  $\O_{C,y}\subset\O_{R,x}$ is integral, 
whence also $\O_{C,y}=\kappa(y)$. According to Proposition \ref{necessary condition gorenstein}, 
the Gorenstein assumption ensures that the field extension $\kappa(y)\subset\kappa(x)$
has degree two. It is purely inseparable because the gluing map $R\ra C$ is a universal homeomorphism.
\qed

\begin{theorem}
\mylabel{ramification regular equation}
Assumptions as in Proposition \ref{ramification regular}.  Let $x\in X$ be a point of codimension two whose image
$y\in Y$ is  contained in the conductor scheme $C$.
Suppose that the local ring $\O_{X,x}$ is regular, and that the residue field extension 
$\kappa=\kappa(y)\subset\kappa(x)$ is trivial.  Then we have an isomorphism  of complete local rings 
$$
\O_{Y,y}^\wedge\simeq \kappa[[a,b,c]]/(a^2-b^2c).
$$
If furthermore the field extension $K\subset\kappa$ is separable, 
the locus of non-smoothness $N=\Sing(Y/K)$ corresponds
to the subscheme defined by the equation $b^2=0$, and $N$ is formally isomorphic to the spectrum of $\kappa[[a,b,c]]/(a^2,b^2)$.
The preimage $\nu^{-1}(N)\subset X$ coincides with $2R\subset X$ in an open neighborhood of $x\in X$.
\end{theorem}

\proof
Since the local ring $\O_{Y,y}^\wedge $ is complete, there exists a \emph{coefficient field} $\kappa\subset\O_{Y,y}^\wedge$ 
(\cite{AC 8-9}, Chapter IX, \S3, No.\ 3, Theorem 1).  
Recall that this is a subfield that bijects onto the residue field.
By assumption, the schemes $R$ and $X$ are regular   at $x\in X$, so the corresponding 
ideal is generated by a single element $u\in\O_{X,x}$,
and we may extend it to   a regular system of parameters $u,v\in\O_{X,x}$.
Then we get an identification $\O_{R,x}^\wedge=\kappa[[u,v]]/(u)=\kappa[[v]]$, with field of fractions $\kappa((v))$.
The subring $\O_{C,y}^\wedge$ is normal by Proposition \ref{condition S2}.
It contains the coefficient field $\kappa$ by construction, and also the square $v^2$  because the field extension
$\Frac(\O_{C,y}^\wedge)\subset\Frac(\O_{R,x}^\wedge)$
is purely inseparable of degree $p=2$. This gives
$\kappa[[v^2]]\subset \O_{C,y}^\wedge\subset \O_{R,x}^\wedge=\kappa[[v]]$.
The composite extension and the   extension on the right both have  degree two, whence the left extension has degree one.
Thus $\O_{C,y}^\wedge=\kappa[[v^2]]$, because the former is normal.
A computation with formal power series immediately shows that the diagram
$$
\begin{CD}
\kappa[[v]]		@<<<	\kappa[[u,v]]\\
@AAA			@AAA\\
\kappa[[v^2]]		@<<<	\kappa[[u,uv,v^2]]
\end{CD}
$$
is cartesian, where the horizontal maps send $u$ and $uv$ to zero.
The three generators $a=uv$, $b=u $ and $c=v^2$ satisfy the relation $a^2=b^2c$.
With the ring  $B=\kappa[[a,b,c]]/(a^2-b^2c)$, we obtain a surjective homomorphism
$B\ra \kappa[[u,uv,v^2]]$. This map is actually bijective by Krull's Principal Ideal Theorem, 
because the   $a^2- b^2c$ is irreducible, viewed as a polynomial in $a$. This establishes
$\O_{Y,y}^\wedge=B$.

Now suppose that the field extension $K\subset\kappa$ is separable. After shrinking $Y$ we find a closed 
embedding into an affine space $\AA^r_K$ for some $r\geq 0$. This gives a surjection of $K$-algebras
$\O_{\AA^r_K,y}^\wedge\ra \O_{Y,y}^\wedge$.
According to \cite{EGA IVa}, Theorem 19.6.4  there is a   coefficient field in $\O_{\AA^r_K,y}^\wedge$ containing $K$.
This induces a  new   coefficient field $\kappa\subset\O_{Y,y}^\wedge$ containing the ground field $K$. The ring extensions
$K\subset\kappa\subset B$ yields an exact sequence
\begin{equation}
\label{omega sequence}
\Omega^1_{\kappa/K}\otimes_\kappa B\lra \Omega^1_{B/K} \lra \Omega^1_{B/\kappa}\lra 0 
\end{equation}
of K\"ahler differentials. The map on the left is injective. To see this, consider the exact sequence
$\Omega^1_{\kappa/K}\otimes_\kappa \Frac(B)\ra \Omega^1_{\Frac(B)/K} \ra \Omega^1_{\Frac(B)/\kappa}\ra 0$.
Here the map on the left is injective, in fact a direct summand, according to \cite{EGA IVa}, Theorem 20.5.7,
because the  extension $\kappa\subset\Frac(B)$ is separable.
Using the injectivity of  $\Omega^1_{\kappa/K}\otimes_\kappa B\ra \Omega^1_{\kappa/K}\otimes_\kappa \Frac(B)$ 
and the functoriality of the sequences of K\"ahler differentials, we infer that the map on the left in \eqref{omega sequence} is injective.

Set $n=\dim(Y)$. Obviously,  $\kappa=\kappa(y)$ is the field of fractions of the integral closed
subscheme $ \overline{\{y\}}\subset Y$ of dimension $n-2$. 
In turn, the extension $K\subset\kappa$ has   transcendence degree $\trdeg( \kappa/K)=n-2$.
According to \cite{A 4-7}, Chapter V, \S 16, No.\ 7, Theorem 5
the vector space $\Omega^1_{\kappa/K}$ is free of rank $n-2$,
thus the $B$-module $\Omega^1_{\kappa/K}\otimes_\kappa B$ is free of the same rank.
Since $K$ has finite $p$-degree, so does $\kappa$, which ensures that $\Omega^1_{B/\kappa}$ is generated
by the differentials $da,db,dc$ modulo the relation $b^2dc=0$ (compare the discussion in \cite{Schroeer 2017},  Section 1). 
Thus $\Omega^1_{B/K}$ has $n+1$ generators and one relation, and we may use them to compute Fitting ideals.
By Proposition \ref{scheme structure},  the locus of non-smoothness $N=\Sing(Y/K)$ is formally defined by the coefficient  $b^2\in B$,
and the description of $\O_{N,y}^\wedge$ follows.
Since $b^2=u^2$ in $\O_{X,x}^\wedge=k[[u,v]]$,
the statement on $\nu^{-1}(N)$ is also clear.
\qed

\medskip
Next, we examine the case that the ramification locus is non-reduced.

\begin{proposition}
\mylabel{ramification multiple}
If the ramification locus takes the form  $R=2R_\red$, 
then  for each point $x\in X$  of codimension one whose image $y\in Y$ is
contained in the conductor scheme $C$, the extension of local rings $\O_{C,y}\subset\O_{R,x}$
is isomorphic to 
$$
\kappa(y)\subset\kappa(x)[\epsilon] \quad\text{or}\quad \kappa(y)[\epsilon]\subset\kappa(x)[\epsilon],
$$
where $\epsilon$ is an indeterminate subject to $\epsilon^2=0$.
In the first case, we have $\kappa(y)=\kappa(x)$. In the second case, the
residue field extension $\kappa(y)\subset\kappa(x)$ is purely inseparable of degree two and we are in characteristic $p=2$.
\end{proposition}

\proof
Write $\kappa=\kappa(x)$,  and choose a coefficient field  so that the local ring $\O_{R,x}$ becomes isomorphic to $\bar{A}=\kappa[\epsilon]$.
According to Proposition \ref{necessary condition gorenstein},
we have to understand the $K$-subalgebras $\bar{B}\subset \bar{A}$ so that $\bar{A}$ is finite over $\bar{B}$
with $\length_{\bar{B}}(\bar{A})=2\length_{\bar{B}}(\bar{B})$.
Setting $\kappa'=\kappa(y)$, we see that $[\kappa:\kappa']$ is either one or two.

In the former case, $\kappa'=\kappa$ and the length of $\bar{A}$ as a module over itself coincides with its length
as a module over $\bar{B}$. It  follows that $\bar{B}=\kappa'$. Replacing the coefficient field in $\bar{A}$ 
by the image of $\bar{B}$, we   reach
the first alternative.

In the second case, the length of $\bar{A}$ as a module over itself is twice its length
as module over $\bar{B}$, and we infer that both Artin rings have length two.
Hence   $\bar{B}=\kappa'[\epsilon']$. Inside the overring $\bar{A}=\kappa[\epsilon]$,
the element $\epsilon'$ is non-zero and lies in the maximal ideal, hence generates the maximal ideal.
Replacing $\epsilon$ by $\epsilon'$, we arrive at the second alternative.
\qed

\begin{theorem}
\mylabel{ramification multiple equation}
Assumptions as in Proposition \ref{ramification multiple}. Let $x\in X$ be a point of codimension two whose image
$y\in Y$ is contained in the conductor scheme $C$. 
Suppose that the local rings $\O_{X,x}$   and $\O_{C,y}$ are regular. Then  
$$
\O_{Y,y}^\wedge\simeq \kappa[[a,b,c]]/(a^2-b^3),
$$
where $\kappa=\kappa(y)$ is the residue field. If furthermore    the field extension $K\subset\kappa$
is separable, the following holds: 
\begin{enumerate}
\item
In characteristic $p=3$, 
the locus of non-smoothness $N=\Sing(Y/K)$  is   defined by the equation
$a=0$,  such that $\O_{N,y}^\wedge\simeq \kappa[[a,b,c]]/(a,b^3)$,
and the preimage $\nu^{-1}(N)\subset X$ coincides with $\frac{3}{2}R=3R_\red$ in a neighborhood of $x\in X$.

\item
In characteristic $p=2$, the subscheme $N=\Sing(Y/K)$ is defined by  $b^2=0$, such that 
$\O_{N,y}^\wedge\simeq\kappa[[a,b,c]]/(a^2,b^2)$. The preimage $\nu^{-1}(N)\subset X$
coincides with $2R=4R_\red$ in a neighborhood of $x\in X$.
\end{enumerate}
\end{theorem}

\proof
Choose a coefficient field $\kappa\subset\O_{X,x}^\wedge$. Since $\O_{C,y}$ is regular,
we are in the first alternative of Proposition \ref{ramification multiple}, hence $R_\red\ra C$ is birational.
The inclusion $\O_{C,y}\subset\O_{R_\red,x}$ is an equality, because $\O_{C,y}$ is normal,
hence the local ring $\O_{R_\red,x}$ is regular.
Therefore the ideal for the reduced ramification locus is generated by 
a member  $u$ of some regular system of parameters $u,v\in\O_{X,x}^\wedge$.
Whence $\O_{X,x}^\wedge=\kappa[[u,v]]$ and $\O_{R,x}^\wedge=k[[u,v]]/(u^2)$.
The rings $\O_{C,y}$ and $\O_{R,x}$ are regular and Cohen--Macaulay, respectively,
whence the finite extension    $\O_{C,y}\subset\O_{R,x}$ of degree two is flat (\cite{EGA IVb}, Proposition 6.1.5).
Thus we may write $\O_{R,x}=\O_{C,y}[U]/(U^2+\mu U+\xi)$ for some generator $U$. By adding an element
from $\O_{C,y}=\O_{R_\red,x}$, we may assume that $U$ is nilpotent.
Then $\mu=\xi=0$, because    polynomial rings over discrete valuation rings are  factorial.
Choose a uniformizer $V\in\O_{C,y}$.
Replacing $u,v\in\O_{X,x}^\wedge $ by representatives of the classes of $U$ and $V$,
we may assume that $\O_{C,y}^\wedge\subset\O_{R,x}^\wedge$ is given by $k[[v]]\subset\k[[u,v]]/(u^2)$.
A computation with formal power series shows that the diagram
$$
\begin{CD}
\kappa[[u,v]]/(u^2)	@<<<	\kappa[[u,v]]\\
@AAA			@AAA\\
\kappa[[v]]		@<<<	\kappa[[u^2,u^3,v]]
\end{CD}
$$
is cartesian. The generators $a=u^3$, $b=u^2$ and $c=v$ satisfy the relations $a^2=b^3$,
and as in the proof for Theorem \ref{ramification regular equation} we infer $\O_{Y,y}^\wedge=\kappa[[a,b,c]]/(a^2-b^3)$.
Likewise, we get the statements on the locus of non-smoothness $N=\Sing(Y/K)$ and its preimage $\nu^{-1}(N)\subset X$.
\qed

\medskip
For the applications we have in mind, we do not have to bother for primes $p\geq 5$, due to the following observation on 
the tangent sheaf:

\begin{proposition}
\mylabel{ramification multiple characteristic}
Assumptions as in Theorem \ref{ramification multiple equation}, with $K\subset \kappa$ separable.
The stalk of the  tangent  sheaf $\Theta_{Y/K,y}$ is   free if and only we  are in characteristic $p\leq 3$.
\end{proposition}

\proof
Set $n=\dim(Y)$, choose a transcendence basis $\xi_1,\ldots,\xi_{n-2}\in \kappa$ for the field extension $K\subset \kappa$,
and write $B=\O_{Y,y}^\wedge=\kappa[[a,b,c]]/(a^2-b^3)$. As in the last paragraph of the proof for Theorem \ref{ramification regular equation},
the $B$-module $\Omega^1_{B/K}$ is generated by $da,db,dc$ together with $d\xi_1,\ldots,d\xi_{n-2}$, modulo the relation
$2ada-3b^2db=0$. Thus we have an exact sequence
$$
B\lra B^{\oplus n+1}\lra\Omega^1_{B/K}\lra 0,
$$
where the map on the right sends the standard basis vectors to the differentials
$da,db,dc$, and the map on the left is given by the $(n+1)\times 1$-matrix $(2a,-3b^2,0,\ldots,0)^t$.
Dualizing gives an exact sequence
$$
0\lra \Theta_{B/K}\lra B^{\oplus n+1}\lra B\lra B/(2a, -3b^2)\lra 0.
$$
In characteristic $p=2$ and $p=3$, the ideal on the right becomes principal, whence the residue class module $M= B/(2a, -3b^2)$
has finite projective dimension. The Auslander--Buchsbaum Formula gives $\pd(M)\leq\dim(B)=2$, 
so the syzygy $\Theta_{B/K}$ must be free.

Conversely, suppose that $\Theta_{B/K}$ is free. Then the $B$-module $B/(2a,-3b^2)$ has finite projective dimension.
Seeking a contradiction, we assume that $p\geq 5$. 
It follows that   $B/I$ has finite projective dimension as well, where $I=(a,b^2,c)$.
The latter ideal has finite colength $l=2$.
Now consider its Frobenius power $I'=(a^p,b^{2p},c^p)\subset B$.
According to a result of Miller \cite{Miller 2003}, Corollary 5.2.3
such a Frobenius power $I'\subset B$ has colength $l'=lp^2=2p^2$.
On the other hand, we may identify $B$ with the subring of $\kappa[[t^2,t^3,c]]\subset \kappa[[t,c]]$
via the identification $a=t^3$ and $b=t^2$.
Then $a^p=t^{3p}$ and $b^{2p}=t^{4p}$. From this we infer that the $3p^2$ elements
$$
t^ic^j\in B/I',\quad i=0,2,3,\ldots, 3p-1,3p+1\quadand  0\leq j\leq p-1 
$$
form a $\kappa$-basis. Thus $l'=3p^2$, contradiction.
\qed

\medskip
It remains to treat the case where the conductor scheme is non-reduced:

\begin{theorem}
\mylabel{conductor multiple equation}
Assumptions as in Proposition \ref{ramification multiple}. Let $x\in X$ be a point of codimension two whose image
$y\in Y$ is contained in the conductor scheme $C$. 
Suppose that the local rings $\O_{X,x}$ and $\O_{R_\red,x}$ and $\O_{C_\red,y}$ are  regular, 
that $\O_{C,y}$ is non-reduced, and that the ring extension $\O_{C,y}\subset\O_{R,x}$ is flat, with
trivial residue field extension  $\kappa=\kappa(y)\subset\kappa(x)$. Then  
$$
\O_{Y,y}^\wedge\simeq \kappa[[a,b,c]]/(a^2-cb^4).
$$
If furthermore    the field extension $K\subset\kappa$
is separable, then the  locus of non-smoothness $N=\Sing(Y/K)$  is   defined by the equation
$b^4=0$,   becomes formally isomorphic to $\kappa[[a,b,c]]/(a^2,b^4)$, and the preimage $\nu^{-1}(N)\subset X$
coincides with $2R=4R_\red$ in a neighborhood of $x\in X$.
\end{theorem}

\proof
Choose a coefficient field $\kappa\subset \O_{Y,y}^\wedge$.
Since $\O_{C,y}$ is non-reduced, so is the overring $\O_{R,x}$, and
we are in the second alternative of Proposition \ref{ramification multiple}, and in particular $p=2$.
Thus we may choose a regular system of parameters $u,v\in\O_{X,x}^\wedge $ so that
$\O_{R,x}^\wedge=\kappa[[u,v]]/(u^2)$. Since $\O_{C_\red,y}$ and $\O_{R_\red,x}$ are regular and Cohen--Macaulay, respectively,
the degree two extension $\O_{C_\red,y}\subset\O_{R_\red,x}$ is flat (\cite{EGA IVb}, Proposition 6.1.5). Choose a uniformizer
$\pi\in\O_{C_\red,y}^\wedge$, such that $\O_{C_\red,y}^\wedge=\kappa[[\pi]]$,
and   write $\O_{R_\red,x}^\wedge= \kappa[[\pi,V]]/(V^2-\varphi)$, for some formal power series $\varphi\in \kappa[[\pi]]$.
The latter is not a square, because the ring is reduced, and becomes  a square in $\kappa$, because the
residue field extension $\kappa=\kappa(y)\subset\kappa(x)$ is trivial.
Replacing $V$, we may assume that $\varphi$ lies in the maximal ideal.
Since $\O_{R_\red,x}^\wedge$ is regular, the element $\varphi\in \kappa[[\pi]]$ is a uniformizer, and
we may assume $ \pi=\varphi$. Clearly, the image $V\in\O_{R_\red,x}^\wedge$ is a uniformizer. 

The nilradical $\mathfrak{n}\subset\O_{C,y}$ 
is a torsion-free module of rank one over $\O_{C_\red,y}$, whence free. Let $U\in \mathfrak{n}$ be
a basis. By assumption the ring extension $\O_{C,y}\subset\O_{R,x}$ is flat, hence the ring $\O_{R,x}/(U)$ is torsion-free.
It is generically reduced by the local description in Proposition \ref{ramification multiple}, hence reduced.
Replacing $u,v$ by representatives
of $U,V$ we may assume that $\O_{R,x}^\wedge=\kappa[[u,v]]/(u^2)$ and
$\O_{C,y}^\wedge=\kappa[[u,v^2]]/(u^2)$.
A computation with formal power series shows that the diagram
$$
\begin{CD}
\kappa[[u,v]]/(u^2)	@<<<	\kappa[[u,v]]\\
@AAA			@AAA\\
\kappa[[u,v^2]]/(u^2)@<<<	\kappa[[u,u^2v,v^2]]
\end{CD}
$$
is cartesian. The generators $a=u^2v$ and $b=u$ and $c=v^2$ satisfy the relation
$a^2= b^4c$. The statement on the locus of non-smoothness $N=\Sing(Y/K)$ and its preimage on $X$
follows as in Theorem \ref{ramification multiple equation}.
\qed

\section{Geometrically non-normal   schemes}
\mylabel{Geometrically non-normal}

Let $F$ be a ground field of characteristic $p>0$,
and $V$ be an algebraic scheme, that is, the structure morphism  $V\ra\Spec(F)$ is separated and of  finite type. We use the  letters $F$ and $V$
because in applications, the former frequently arises as a \emph{function field} and the latter  becomes the \emph{generic fiber} of 
some fibration. Recall that $V$ is called \emph{geometrically normal} if $V\otimes_FK$ is normal
for all field extensions $F\subset K$, with similar locution for other scheme-theoretic properties of $V$.

Given a field extension $F\subset K$, we set $Y=(V\otimes_FK)_\red $ and write $\nu:X\ra Y$
for the normalization. 
Note that for the schemes $X$ and $Y$,  we    regard $K$ rather than $F$
as the ground field.
In our applications, we are mainly interested
in the situation where $V$ is geometrically integral  but not geometrically normal.
In this section, however, we make some fairly general observations.
Let us start with the following well-known fact:

\begin{lemma}
\mylabel{geom-normal normalization}
There is a finite purely inseparable field extension $F\subset K$ so
that  $X$ is geometrically normal and $Y$ is geometrically reduced.
\end{lemma}

\proof
Choose a perfect closure $F\subset F^\perf$.
Then $Y_\infty=(V\otimes_FF^\perf)_\red$ is geometrically reduced, and its normalization $X_\infty$ is geometrically normal.
Let $F\subset F_\lambda\subset F^\perf$, $\lambda\in L$ be the filtered ordered set of finite subextensions.
According to \cite{EGA IVc}, Theorem 8.8.2, for some index $\lambda$ there is a closed subscheme
$Y_\lambda\subset V\otimes_FF_\lambda$ and a finite morphism $X_\lambda\ra Y_\lambda$
inducing $X_\infty\ra Y_\infty\subset V\otimes_FF^\perf$.
It follows that $X_\lambda$ is geometrically normal  and $Y_\lambda$ is geometrically reduced.
This yields the desired finite purely inseparable field extension $K=F_\lambda$.
\qed

\medskip
Now suppose that   $F\subset K$  is a field extension so that 
the   $X$ is geometrically normal and $Y$ is geometrically reduced.
As explained in Appendix \ref{Finite modifications}, 
the finite birational morphism $\nu:X\ra Y$ comes with  the  conductor square 
$$
\begin{CD}
R	@>>>	X\\
@VVV		@VV\nu V\\
C	@>>>	Y,
\end{CD}
$$
where $C\subset Y$ is the  conductor scheme  and $R\subset X$ is the  ramification locus.
Both correspond to  the coherent sheaf $\mathfrak{C}$ defined as the annihilator
of $(\nu_*\O_X)/\O_Y$, which is an ideal sheaf in both $\O_Y$ and $\nu_*(\O_X)$.
For each further field extension $K\subset K'$, the scheme $Y'=Y\otimes_KK'$ is
reduced, and $X'=X\otimes_KK'$ is normal, whence the induced map $X'\ra Y'$ is the normalization.
Since kernels for homomorphisms between quasicoherent sheaves,
and in particular annihilator ideals for coherent sheaves, commute with ground field extensions,
we see that the base-change of the above diagram along $K\subset K'$
is the conductor square for the normalization $X'\ra Y'$.
Applying Lemma \ref{geom-normal normalization} to the conductor scheme and the ramification locus,
we obtain:

\begin{proposition}
\mylabel{geom-reduced conductor}
There is a finite purely inseparable field extension $F\subset K$ such that
$X$ is geometrically normal, and the schemes $Y$, $C_\red$, $R_\red$ are geometrically reduced.
\end{proposition}

Now suppose additionally that the algebraic  scheme $V$ is proper, and furthermore satisfies
the condition $h^0(\O_V)=1$.
Then we can form the  \emph{Picard scheme} $\Pic_{V/F}$, its connected component of the origin $\Pic^0_{V/F}$
and the resulting \emph{N\'eron--Severi   scheme} $\NS_{V/F}$.
These are commutative \emph{group schemes} sitting in a short exact sequence 
$$
0\lra \Pic^0_{V/F}\lra \Pic_{V/F}\lra \NS_{V/F}\lra 0.
$$
The cokernel  is a smooth zero-dimensional group scheme, whence is determined by  the group
$\NS_{V/F}(F^\sep)$, where $F^\sep$ is some separable closure, together with
the action of the Galois group $\Gal(F^\sep/F)$.
The group of rational points   is denoted  by $\NS(V/F)=\NS_{V/F}(F)$. 
Inside, we have the usually smaller subgroup
$\NS(V)$ of rational points coming from   invertible sheaves   on $V$.

We say that $V$ \emph{has completely constant N\'eron--Severi group scheme} if  
the inclusion $\NS(V)\subset\NS_{V/F}(F^\sep)$ is an equality.
In other words, the group scheme is constant, and each point actually corresponds to
an invertible sheaf. In turn, the canonical map $\Pic(V)\ra\NS_{V/F}(F^\sep)$
is surjective. The kernel is $\Pic^0(V)=\Pic(V)\cap \Pic^0_{V/F}(F)$, and we get   identifications
$$
\NS(V)=\NS(V/F)=\NS_{V/F}(F^\sep)=\Pic(V)/\Pic^0(V),
$$
as customary over algebraically closed ground fields. The following observation will be useful:

\begin{lemma}
\mylabel{completely constant}
There is a finite separable field extension $F\subset F'$ so that   the base-change $V'=V\otimes_FF'$
has completely constant   N\'eron--Severi group scheme $\NS_{V'/F'}$.
\end{lemma}

\proof
By finiteness of the base number, the underlying set of the group scheme $\NS_{V/F}$ is countable.  
For each   point $l\in \NS_{V/F}$ the residue field extension  $F\subset \kappa(l)$ 
is finite and separable. Choose
an enumeration $l_n\in\NS_{V/F}$, $n\geq 0$ of them.
The resulting \'etale algebras $A_n=\kappa(l_0)\otimes_F\ldots\otimes_F\kappa(l_n)$
form an increasing sequence $A_0\subset A_1\subset\ldots$, and we write $A=\dirlim_{n\geq 0}A_n$ for their direct limit.
Choose a residue field $F'=A/\maxid$. Clearly, the field extension $F\subset F'$ is separable algebraic.
The images $F_n\subset F'$ 
of the \'etale algebras $A_n$ form an increasing sequence $F_0\subset F_1\subset\ldots$
of finite separable field extensions.
If  this sequence is stationary, the field $F'=F_n$, $n\gg 0$ solves our problem.
Seeking a contradiction, we assume that the sequence is not stationary. Choose an embedding
of $F'$ into some separable closure $F^\sep$.
Then the images $H_n\subset\NS_{V/F}(F^\sep)$ of the groups $\NS_{V/F}(F_n)$ form
a   sequence $H_0\subset H_1\subset\ldots$ that is not stationary, but lies inside
the finitely generated abelian group $\NS_{V/F}(F^\sep)$, contradiction.
\qed

\medskip
Note that the corresponding statement for the Picard scheme does not hold,
for example if $F$ is finite and  $\dim\Pic^0_{V/F}\geq 1$,
or if $F$ is imperfect and the Picard scheme contains a copy of $\GG_a$ or $\GG_m$.

The  ramification locus $R\subset X$ has in general $h^0(\O_R)\neq 1$,
so forming the Picard scheme   is problematic.
The same  difficulty occurs for the conductor scheme $C\subset Y$.
To avoid cumbersome statements, we say that a proper $K$-scheme $Z$ has 
\emph{completely constant N\'eron--Severi group scheme}
if the connected components $Z_1,\ldots,Z_r\subset Z$ satisfy $h^0(\O_{Z_i})=1$
and have completely constant $\NS_{Z_i/K}$.
Note that the condition $h^0(\O_{Z_i})=1$ automatically holds if $Z_i$ is geometrically connected and
geometrically reduced.

\begin{proposition}
\mylabel{adapted}
There is a finite separable field extension $F\subset F'$ and
a finite purely inseparable field extension $F\subset K'$ so that,
for the resulting finite field extension $K=F'\otimes_FK'$, the following holds:
\begin{enumerate}
\item The normalization  $X$ of $Y=(V\otimes_FK)_\red$ is geometrically normal.
\item The inclusion $\Sing(Y)\subset \Sing(Y/K)$ is an equality of closed sets.
\item The schemes $Y$, $C_\red$ and $R_\red$ are geometrically reduced.
\item The schemes $X$, $Y$,    $C_\red$ and  $R_\red$ have completely constant N\'eron--Severi group scheme.
\end{enumerate}
\end{proposition}

\proof
According to Proposition \ref{geom-reduced conductor}, there is a finite purely inseparable field extension $F\subset K'$ so that 
the conditions (i) and (iii) hold for any  separable field extension $F\subset F'$.
In light of Proposition \ref{locus as image}, we may also achieve (ii).

For the last property, choose a separable closure $F\subset F^\sep$, and consider the
filtered order set of finite subextensions $F\subset F_\lambda\subset F^\sep$, $\lambda\in L$.
Write $Z_1,\ldots,Z_r\subset X\otimes_FF^\sep$
for the connected components.
According to \cite{EGA IVc}, Theorem 8.8.2 there  are closed subschemes 
$Z_{i,\lambda}\subset X\otimes_FF_\lambda$ for some index $\lambda\in L$ inducing the $Z_i\subset X\otimes_FF^\sep$.
Clearly, the $Z_{i,\lambda}$ are geometrically connected and geometrically reduced, thus $h^0(\O_{Z_{i,\lambda}})=1$.
According to Lemma \ref{completely constant}, we may enlarge the index $\lambda\in L$ to achieve that
the $Z_{i,\lambda}$ have completely constant N\'eron--Severi group scheme.
The same arguments apply for $Y$,  $C_\red$ and $R_\red$.
Summing up, we may choose the index $\lambda\in L$ so large that, for $F'=F_\lambda$ 
also condition (iv) holds.
\qed

\medskip
The following notion turns out to be useful:

\begin{definition}
\mylabel{definition adapted}
A proper scheme $V$, with structure morphism $V\ra\Spec(F)$, and a purely inseparable
field extension $F\subset K$ are called \emph{adapted} if
  conditions (i)--(iv) of Proposition \ref{adapted} hold.
\end{definition}

According to Proposition \ref{adapted}, for any proper $F$-scheme $V$ there is a finite
separable extension $F\subset F'$ and a finite purely inseparable extension $F\subset K'$
such that   $V\otimes_FF'\ra\Spec(F')$ and $F'\subset F'\otimes_FK'$ are adapted.
We usually reduce to this situation when analyzing the existence or non-existence
of certain proper $F$-schemes $V$. The following observation will be   useful:

\begin{proposition}
\mylabel{neron-severi sequence}
Suppose that the proper scheme $V$ and the finite purely inseparable extension $F\subset K$
are adapted, and that the ramification locus $R\subset X$ satisfies  $h^0(\O_R)=1$ and $h^1(\O_R)=0$.
Then  
$$
h^0(\O_X)=h^0(\O_Y)=h^0(\O_C)=1 \quadand  h^1(\O_Y)=h^1(\O_X)+h^1(\O_C).
$$
If moreover $h^2(\O_Y)=0$, the canonical  sequence
$$
0\lra\NS(Y)\lra \NS(X)\oplus\NS(C)\lra \NS(R)
$$
of N\'eron--Severi groups is exact.
\end{proposition}

\proof
The gluing map $R\ra C$ is schematically dominant, which implies that
the homomorphism $H^0(C,\O_C)\ra H^0(R,\O_R)$ is injective, so
$h^0(\O_C)=1$.
The short exact sequence \eqref{exact coherent} yields an exact sequence
$$
0\lra H^0(Y,\O_Y)\lra H^0(X,\O_X)\oplus H^0(C,\O_C) \lra H^0(R,\O_R)\lra 0.
$$
The map on the right is indeed  surjective, because $h^0(\O_R)=1$.
Since the normalization $\nu:X\ra Y$ is schematically dominant, the
map $H^0(Y,\O_Y)\ra H^0(X,\O_X)$ is injective.
This map is actually bijective, by the  preceding exact sequence and $h^0(\O_C)=h^0(\O_R)$.
We thus have an identification $H^0(X,\O_X)=H^0(Y,\O_Y)$. 
The scheme $V$ is geometrically connected, by our overall assumption $h^0(\O_V)=1$,
so the reduction   $Y=(V\otimes_FK)_\red$ is geometrically connected as well.
It is also geometrically reduced, because  $V$ and $F\subset K$ are adapted.
Consequently   $h^0(\O_X)=h^0(\O_Y)=1$.  
The short exact sequence \eqref{exact coherent} also yields an exact sequence
$$
0\lra H^1(Y,\O_Y)\lra H^1(X,\O_X)\oplus H^1(C,\O_C)\lra 0,
$$
whence $h^1(\O_Y)=h^1(\O_X)+h^1(\O_C)$. 
The above cohomology groups are the Lie algebras for
the respective Picard groups, and we also infer that
the restriction map 
\begin{equation}
\label{group schemes}
\Pic_{Y/K}^0\lra\Pic_{X/K}^0\oplus\Pic_{C/K}^0
\end{equation}
of group schemes has finite \'etale kernel. 
Now suppose that also  $h^2(\O_Y)=0$. Then  the group scheme $\Pic_{Y/K}^0$ is smooth  (\cite{Mumford 1966}, Lecture 27), 
of dimension $h^1(\O_Y)=h^1(\O_X)+h^1(\O_C)$, so the restriction map is surjective.
The short exact sequence \eqref{exact multiplicative} gives an exact sequence
$$
0\lra \Pic(Y)\lra \Pic(X)\oplus\Pic(C)\lra\Pic(R)
$$
of Picard groups. The map on the left is indeed injective, because the mapping
$K^\times=H^0(X,\O_X^\times)\ra H^0(R,\O_R^\times)=K^\times$ is surjective. 
This also holds for all extension fields, and so \eqref{group schemes} is actually an isomorphism of groups schemes.
The assumption 
$H^1(R,\O_R)=0$ ensures that $\Pic^0_{R/K}=0$, in particular $\NS(R)=\Pic(R)$.
In turn, we obtain    a commutative diagram
$$
\begin{CD}
0	@>>> 	\Pic^0(Y)	@>>>	\Pic^0(X)\oplus\Pic^0(C)	@>>> 0		@>>> 0\\
@.		@VVi' V			@VVi V				@VVi''V\\
0	@>>> 	\Pic(Y)		@>>> 	\Pic(X)\oplus\Pic(C)	@>>>	\Pic(R)
\end{CD}
$$
with exact rows. Consequently, the desired exact sequence of N\'eron--Severi groups
$0\ra\Cokernel(i')\ra\Cokernel(i)\ra\Cokernel(i'')$   follows from the Snake Lemma.
\qed

\medskip
In the above situation, we can regard $\NS(Y)$ as 
the kernel of some homomorphism between finitely generated abelian groups.
If moreover the N\'eron--Severi groups of $X,C,R$ are torsion-free, we thus have
$$
\NS(Y)=\Kernel(\Psi)\subset\ZZ^{\oplus n}
$$
for some \emph{integral matrix} $\Psi\in\Mat_{m\times n}(\ZZ)$,
where the size of the matrix is given by  $m=\rho(R)$ and $n=\rho(X)+\rho(C)$.
Here $\rho(X)=\rank \NS(X)$ etc.\ denotes   \emph{Picard numbers}.

In the next sections, we shall exploit this in the following way: every Cartier divisor $D\subset V$
induces a Cartier divisor $D_K\subset 								Y$, which represents an integral vector in
the kernel of the matrix $\Psi$. This applies in particular for the canonical divisor $D=K_V$,
if the scheme $V$ is Gorenstein, or the reduction $D=N_\red$ of the divisorial part $N\subset\Sing(V/F)$ for the locus
of non-smoothness, if the scheme $V$ is  locally factorial.
As we shall see, sometimes geometric reasons  preclude
the existence of such integral solutions for the system of linear equations 
$$\Psi
\begin{pmatrix}
x_1\\
\vdots\\
x_n	\\
\end{pmatrix} = 
\begin{pmatrix}
0\\
\vdots\\
0	\\
\end{pmatrix}.
$$
We close this section with a well-known useful observation:

\begin{proposition}
\mylabel{universal homeomorphism}
Suppose that the local rings $\O_{V,a}$, $a\in V$ are geometrically unibranch.
Then the normalization map $\nu:X\ra Y=(V\otimes_FK)_\red$  is a universal homeomorphism.
\end{proposition}

\proof
The conditions on the local rings means that the spectrum of the strict henselization $\O_{V,a}^s$
is irreducible. If $F\subset K$ is a finite purely inseparable field extension,
then the ring $\O_{V,a}^s\otimes_FK$ remains local,   indeed strictly local,
and coincides with the strictly local ring at the point $b\in V\otimes_FK$ corresponding to $a\in V$,
according to \cite{EGA IVd}, Proposition  18.6.8.
In particular, the normalization $A'$ of the strictly local 
integral domain $A=(\O_{V\otimes_FK,b}^s)_\red$
yields a universal homeomorphism $\Spec(A')\ra\Spec(A)$.
It follows that the normalization $X\ra Y=(V\otimes_FK)_\red$ is a universal homeomorphism.
\qed

\medskip
Note that the conclusion holds, in particular, when $V$ is normal.

\section{Regular del Pezzo surfaces}
\mylabel{Regular del Pezzo}

Let $F$ be a ground field. Throughout the paper, we use the following  general notion of del Pezzo surfaces:

\begin{definition}
A \emph{del Pezzo  surface} is a proper  two-dimensional
scheme $V$ with $h^0(\O_V)=1$ that is Gorenstein, and whose dualizing sheaf $\omega_V$
is antiample. 
\end{definition}

The most immediate examples are  surfaces of degree two or three
in $\PP^3$, or complete intersections of two quadrics in $\PP^4$.
The antiample sheaf $\omega_V$ has no sections, at least if $V$ is integral,
hence $h^2(\O_V)=0$. The selfintersection number $(\omega_V\cdot\omega_V)>0$ is called
the \emph{degree} of the del Pezzo surface, and the integer $h^1(\O_V)\geq 0$ is commonly referred to as  the \emph{irregularity}.
Examples of del Pezzo surfaces with irregularity $h^1(\O_V)>0$ were constructed
by Reid \cite{Reid 1994}, the second author \cite{Schroeer 2007} and Maddock \cite{Maddock 2016}.

In what follows, we suppose that the ground field $F$ has characteristic $p>0$,
and assume that the del Pezzo surface $V$ is normal, locally factorial and geometrically integral,  but   geometrically non-normal. 
We are mainly interested in the case that $V$ is even regular, but the weaker
assumptions of local factoriality lies at the core  for most our arguments.
The locus of non-smoothness $\Sing(V/F)$ 
contains a non-empty divisorial part $N\subset\Sing(V/F)$, which
is an effective  Cartier divisor. Let $D=N_\red$ be its reduction. 
The recurrent idea of this paper is to study the behavior of the effective Cartier divisor $D\subset V$
under base-change, in particular if the ground field has $\pdeg(F)=1$.
 
After replacing the ground field $F$ by some finite separable extension,
we may  assume that there is some finite purely inseparable field extension $F\subset K$
so that $V\ra\Spec(F)$ and $F\subset K$ are adapted, according to Proposition \ref{adapted}.
As in the preceding section, we write   $\nu:X\ra Y$ for the normalization of
$Y=V\otimes_FK$.  Note that $X$ is Cohen--Macaulay,
but not necessarily Gorenstein. Let
\begin{equation}
\label{conductor_square}
\begin{CD}
R	@>>>	X\\
@VVV		@VV\nu V\\
C	@>>>	Y
\end{CD}
\end{equation}
be the conductor square, where $C\subset Y$ is the \emph{conductor curve} and $R\subset X$ is the \emph{ramification divisor},
as discussed in Appendix \ref{Finite modifications}.
Both are Weil divisor
that are not necessarily  Cartier. Moreover, the vertical maps are 
universal homeomorphisms.
Note that  we regard $K$ as a new ground field for  $Y$ and $X$.

Recall that the \emph{Hirzebruch surface with numerical invariant $e\geq 0$}
is the smooth surface $S=\Proj(\Sym\shE)$ for the locally free sheaf $\shE=\O_{\PP^1}\oplus\O_{\PP^1}(e)$.
It comes with a ruling $r:S\ra\PP^1$ and a section $E=\Proj(\Sym\shL)$ with $\shL=\O_{\PP^1}(e)$, which has 
selfintersection number $E^2=-e$. 
If $e>0$ the curve $E\subset S$ is negative-definite and contracts to a rational singularity.
The resulting \emph{contracted Hirzebruch surface} can be regarded as the weighted projective space 
$\PP(1,1,e)=\Proj K[T_0,T_1,T_2]$, with grading  $\deg(T_0)=\deg(T_1)=1$ and $\deg(T_2)=e$.  
One may view it also as a toric variety.

By abuse of notation, we say that a proper curve  $Z$  a \emph{split conic} if 
it is isomorphic to a divisor of degree two inside $\PP^2$ given by one of the following
three homogeneous equations:
$$
T_0^2+T_1T_2=0\quad\text{or}\quad T_0T_1=0\quad\text{or}\quad (T_0+T_1+T_2)^2=0.
$$
In the first case $Z$ is isomorphic to the \emph{projective line} $\PP^1$.
In the second case, we say that $Z$ is a \emph{pair of lines} $\PP^1\cup\PP^1$.
In the last case, $Z$ is isomorphic to the \emph{split ribbon}
$\PP^1\oplus\O_{\PP^1}(-1)$  on the projective line $\PP^1$
with ideal sheaf $\O_{\PP^1}(-1)$. Any such ribbon is split, because
$\Ext^1(\O_{\PP^1}(-1),\O_{\PP^1})=0$.
We refer to \cite{Bayer; Eisenbud 1995} for the general theory of ribbons.

\begin{proposition}
\mylabel{twisted form}
The normal surface $X$ is   either the projective plane $\PP^2$, a contracted Hirzebruch surface $\PP(1,1,e)$
with numerical invariant $e\geq 2$, or a Hirzebruch surface $S$ with   $e\geq 0$.
The ramification curve  $R$ is a split conic, in particular  $h^0(\O_R)=1$ and $h^1(\O_R)=0$.
\end{proposition}

\proof
By assumption,  $V$ is  geometrically integral and  $X$ is  geometrically normal, respectively.
In light of  Reid's Classification (\cite{Reid 1994}, Theorem 1.1),
the assertion holds if we pass to some finite field extension $K\subset K'$.
It remains to check that the assertion already holds over $K$.
For this we use the assumption that $V$ and $F\subset K$ are adapted.

We first treat the ramification curve $R$ and set $R'=R\otimes_KK'$.
If $R'$ is a projective line, we choose an invertible sheaf $\shL$ on $R$ of degree one,
which defines an isomorphism $R\ra\PP^1$.
If $R'=R'_1\cup R'_2$ is a pair of lines, there are invertible sheaves $\shL_1$ and $\shL_2$
on $R$ with $\deg(\shL_i|R'_j)=\delta_{ij}$. The resulting morphisms $f_i:R\ra\PP^1$
reveal that $R$ is a pair of lines.
Now suppose that $R'$ is a split ribbon.
The scheme $R_\red$ is geometrically reduced, and its ideal sheaf $\shI\subset\O_R$  
is   invertible of degree one on $R_\red$. It follows $R_\red=\PP^1$.
Such ribbons   split, so $R=\PP^1\oplus\O_{\PP^1}(-1)$.

In all cases, we may regard $R\subset \PP^2$ as a divisor of degree two.
The resulting short exact sequence $0\ra\O_{\PP^2}(-2)\ra\O_{\PP^2}\ra \O_R\ra 0$
induces a long exact sequence
$$
 H^0(\O_{\PP^2})\ra H^0(\O_R)\ra H^1(\O_{\PP^2}(-2))\ra H^1(\O_{\PP^2})\ra H^1(\O_R)\ra H^2(\O_{\PP^2}(-2)).
$$
The cohomology groups $H^i(\PP^2,\O_{\PP^2}(-2))$ vanish for all $i\geq 1$, and
it follows that $h^0(\O_R)=1$ and $h^1(\O_R)=0$.

We now turn to the normal surface $X$.
Suppose that $X'=X\otimes F$ is isomorphic to a Hirzebruch surface. Fix a ruling
$r':X'\ra\PP^1_{K'}$ and let $\shL'$ be the preimage of $\O_{\PP^1_{K'}}(1)$.
This invertible sheaf descend to to an invertible sheaf $\shL$ on $X$.
The latter is semiample, and defines a ruling $r:X\ra\PP^1$.
Now choose an invertible sheaf $\shN'$ that has degree one on the fibers of $r':X'\ra\PP^1_{K'}$.
This also descends, and it follows that $X$ is a Hirzebruch surface.
The case that $X'$ is the projective plane is treated in a similar way.

Now suppose that $X'$ is a contracted Hirzebruch surface, and let  $S'\ra X'$ be the 
minimal resolution of singularities. Let $a'\in X'$ be the singularity and $a\in X$ be its image.
Using \cite{Artin 1966}, Theorem 4, 
we infer that $S'\ra X'$ is the blowing-up of
the reduced center $a'\in X'$, and that the exceptional divisor $E'\subset S'$
is a projective line.
According to Reid's Classification, the ramification divisor $R'\subset X'$
is linearly equivalent to $2H'$, where $H'\subset X'$ is the image
of any fiber $F'\subset S'$ from the ruling. 
Note that $eH'$ generates the Picard group, but  that $H'$ is not Cartier.

Suppose first that there is a Weil divisor $H\subset X$ whose base-change is linearly
equivalent to $H'\subset X'$.
Since $h^0(\O_X(H))=2$, there are two such Weil divisors $H_1\neq H_2$.
Their intersection $H_1\cap H_2$ has length one  and contains $a\in X$.
So its blowing-up $S\ra X$ yields a twisted form of the Hirzebruch surface,
and the exceptional divisor $E\subset S$ is a twisted form of the projective line.
The latter intersects the strict transform of $H$ in a rational point,
so $E=\PP^1$. Arguing as above, one easily  see that $S$ is a Hirzebruch surface
and hence $X=\PP(1,1,e)$.

Seeking a contradiction, we now assume that there is  no Weil divisor on $X$ inducing $H'\subset X'$.
Recall that the ramification divisor $R$ is a split conic.
It must be a projective line if $a\in R$. This yields a contradiction,
because on $X'$ there are no lines linearly equivalent to $2H'$ passing through
the singularity.
Thus $R\subset X$ lies in the smooth locus. Hence $R$ and the linearly equivalent $2H$ are  Cartier, and we must have $e=2$.
It also follows that the split conic $R$ is a projective line, with selfintersection  $R^2=2$, and one easily computes $h^0(\O_X(R))=4$.
Since $\Sing(Y)=\Sing(Y/K)$, the local ring $\O_{X,a}=\O_{Y,a}$ is singular.
It must be a  twisted form of the rational double point of type $A_1$ given by the equation
$z^2-xy=0$. We see that the  local Artin scheme $\Sing(X/K)$ has length two,
and we conclude with Lemma \ref{regular twisted form} that $a\in X$ is a rational point.
To proceed, choose two further  rational points $b,c\in R$  and consider the reduced closed
subscheme $Z=\{a,b,c\}$. The invertible sheaf $\shL=\O_X(R)$ has $h^0(\shL)=4$.
For dimension reasons, the  restriction map $H^0(X,\shL)\ra H^0(Z,\shL|Z)$ is not injective,
so there is an effective Cartier divisor $A\subset X$ linearly equivalent to $R$ and containing $Z$.
Examining its base-change to $X'$, we easily infer that $A$ is reducible.
Consequently, there is a Weil divisor on $X$ inducing $H'\subset X'$, contradiction.
\qed

\medskip
In what follows, we shall use the following notation: 
If $X=\PP^2$ is the projective plane,   write $H\subset X$ for a hyperplane. Then $\Pic(X)$ is freely generated by the class of $H$,
the intersection pairing is given by $H^2=1$, and the canonical class is $K_X=-3H$.

If $X=\PP(1,1,e)$ is the contraction of a  Hirzebruch surface $S$ with numerical invariant $e\geq 2$, 
we write $E\subset S$ for the section with   $E^2=-e$
and $F\subset S$ for a fiber of the ruling $r:S\ra\PP^1$. 
(It should be clear from the context weather the symbol $F$ means a fiber for the ruling or
a ground field.)
The Picard group $\Pic(S)$ is freely generated by the classes of $E$ and $F$,
the intersection form has Gram matrix $(\begin{smallmatrix}-e&1\\1&0\end{smallmatrix})$, and the canonical class is $K_S=-2E-(e+2)F$.
Let $f:S\ra X$ be the contraction of the negative-definite curve $E\subset S$ 
and write $H=f(F)$ for the image of the fiber. The Weil divisor $H\subset X$ is not Cartier,
but $eH$ is Cartier, with $f^{*}(eH)=eF+E$, and freely generates $\Pic(X)$. The intersection pairing
is given by $(eH)^2=e$, and the canonical class is $K_X=-(e+2)H$. Note that $X$ is Gorenstein if and only if $e=2$.

If $X=S$ is a Hirzebruch surface with numerical invariant $e\geq 0$, we likewise write
$E,F\subset S$ for the section with $E^2=-e$ and the fiber for a  ruling. Here the canonical class is given by  $K_X=-2E-(e+2)F$.

We now tabulate the five possibilities that  follow from  Reid's Classification of non-normal
del Pezzo surfaces over algebraically closed ground fields (\cite{Reid 1994}, Theorem 1.1):

\newcounter{Reid}\setcounter{Reid}{0}
\renewcommand{\theReid}{\rm \roman{Reid}}
\newcommand{\lReid}[1]{\refstepcounter{Reid}({\rm \roman{Reid}})\label{#1}}
\newcommand{\refReid}[1] {{\rm (\ref{#1})}}
\begin{theorem}
\mylabel{reid's classification}
If the del Pezzo surface $V$ and the field extension $F\subset K$ are adapted, then the possibilities
for the normalization $X\ra Y=V_K$ are as follows:
$$
\begin{array}[t]{l| *{5}{|c} }
\text{\rm Case}			& \lReid{P,2H}	& \lReid{P,H}	& \lReid{cS,2H}	& \lReid{S,E+F}		& \lReid{S,E}\\
\hline				&		&		&  			& &\\[-2ex]
X				& \PP^2		& \PP^2		& \PP(1,1,e)		& S			& S \\
e				& 1 		& 1		& \geq 2		& \geq 0		& \geq 0 \\
\O_X(R)				& 2H		& H		& 2H			& E+F			& E\\
\omega_X			& -3H		& -3H		& -(e+2)H		& -2E-(e+2)F		& -2E-(e+2)F\\
\nu^*(\omega_Y)			& -H		& -2H		& -eH			& -E-(e+1)F		& -E-(e+2)F\\
K_Y^2				& 1		& 4		& e			& e+2			& e+4\\
\end{array}
$$
\end{theorem}

\medskip 
We refer to this table as \emph{Reid's Classification}.
Note that Case \refReid{P,2H} and Case \refReid{cS,2H} may be treated on the same footing,
because $\PP^2=\PP(1,1,1)$. It is remarkable that in all cases we have $\rho(V)\leq 2$.
In particular, the dualizing sheaf ceases to be antiample if we blow-up $V$ in two closed points.
For later use, we record the following fact:

\begin{proposition}
\mylabel{pullback injective}
With the exception of case (v) in Reid's Classification, the pullback map $\Pic(V)\ra\Pic(C)$, $\shL\mapsto\shL_C$ is injective.
\end{proposition}

\proof
In any case, the projection $Y=V\otimes_FK\ra V$ induces  an inclusion $\Pic(V)\subset\Pic(Y)$,
and the short exact sequence \eqref{exact multiplicative} yields an  exact sequence
$$
0\lra \Pic(Y)\lra\Pic(X)\oplus\Pic(C)\lra\Pic(R).
$$
The map on the left is indeed injective, because $h^0(\O_R)=1$.
For the cases in question, the restriction map $\Pic(X)\ra\Pic(R)$ is injective, and it follows
that the restriction map $\Pic(Y)\ra\Pic(C)$ is injective as well.
\qed

\medskip
Since $h^0(\O_R)=1$ and $h^1(\O_R)=0$ and $h^2(\O_Y)=0$, we   can apply Proposition \ref{neron-severi sequence}  and conclude that we also
have an exact sequence
\begin{equation}
\label{neron-severi exact sequence}
0\lra \NS(Y)\lra\NS(X)\oplus\NS(C)\stackrel{\Psi}{\lra} \NS(R).
\end{equation}
Note that all terms in this sequence are finitely generated free abelian groups,
so we may regard $\NS(Y)$ as the kernel of some matrix $\Psi\in\Mat_{s\times (r+s)}(\ZZ)$,
where $r=\rho(X)=\rho(Y)$ and $s=\rho(C)=\rho(R)$.  
Since $Y$ is Gorenstein,
the class of the dualizing sheaf $\omega_Y$ defines an element  in this  kernel. This already gives an
important restriction:

\begin{proposition}
\mylabel{no pair of lines}
The ramification divisor $R\subset X$ is either the projective line $\PP^1$ or
the split ribbon $\PP^1\oplus\O_{\PP^1}(-1)$.
\end{proposition}

\proof
Seeking a contradiction, we assume that  $R=R_1\cup R_2$ is a pair of lines.
By Proposition \ref{universal homeomorphism}, 
the map $R\to C$ is a universal homeomorphism. Setting $C_i=\nu(R_i)$ we get $C=C_1\cup C_2$.  
According to Proposition \ref{necessary condition gorenstein} the induced morphisms $R_i\ra C_i$ have degree two,
and the characteristic must be $p=2$.

We now go through the five cases of Reid's Classification:
In Case  \refReid{P,2H}, the matrix $\Psi$ whose kernel gives $\NS(Y)$ takes the form 
$$
\Psi=
\begin{pmatrix}
1	& -2	& 0\\
1	& 0	& -2\\
\end{pmatrix}
\in\Mat_{2\times 3}(\ZZ).
$$
The entries in the first column are the intersection numbers $(H\cdot R_i)=1$,
whereas the other non-zero entries are the negative degrees of the gluing maps $\nu\colon R_i\ra C_i$.
We have $\nu^*(K_Y)=-H$, so the class of $K_Y\in\Kernel(\Psi)$ is a column vector of the form
$-(1,m,n)^t$ for some integers $m,n\in\ZZ$. The only solution  is $m,n=1/2$, which is not integral, contradiction.

In Case \refReid{cS,2H}, the N\'eron--Severi group $\NS(X)$ is generated by $eH$, which has intersection numbers
$(eH\cdot R_i)=1$. Consequently, the matrix describing $\NS(Y)$ is as in the previous paragraph,
and we get a contradiction again.
In the  cases \refReid{P,H} and \refReid{S,E} the ramification divisor $R$ is not linearly
equivalent to the sum of two effective divisors, so these cases are impossible as well.

It remains to deal with case \refReid{S,E+F}, where $X=S$ is a Hirzebruch surface with numerical invariant $e\geq 0$.
Then $\Pic(X)$ is freely generated by $E,F\subset X$. The ramification divisor is $R=E\cup F$, say
with $R_1=E$ and $R_2=F$. Now the matrix describing $\NS(Y)$ takes the form
$$
\Psi=
\begin{pmatrix}
-e	& 1	& -2	& 0\\
1	& 0	& 0	& -2\\
\end{pmatrix}
\in\Mat_{2\times 4}(\ZZ),
$$
where the entries on the left half are given by the intersection numbers
$$
(E\cdot R_1) = -e,\quad (E\cdot R_2)=1,\quad (F\cdot R_1) =1\quadand (F\cdot R_2)=0.
$$
We have $\nu^*(K_Y)=-E-(e+1)F$, whence $K_Y\in \Kernel(\Psi)$ corresponds to a column vector
of the form $-(1,e+1,m,n)^t$ for some integers $m,n\in \ZZ$. However, the only solution has
$n=1/2$, again a  contradiction.
\qed

\medskip
Recall that $N\subset \Sing(V/F)$ denotes the divisorial part of the locus of non-smoothness.
Then $N\subset V$ is a Cartier divisor, since our normal surface $V$ is assumed to be locally factorial.
In turn, we get a Cartier divisor $N_K\subset V_K=Y$
whose support coincides with the conductor curve $C\subset Y$. 
Another   consequence of the Gorenstein condition is the following:

\begin{proposition}
\mylabel{reduced C}
The conductor curve $C$ is geometrically integral.
\end{proposition}

\proof
Proposition \ref{no pair of lines} ensures that  $C$ is geometrically irreducible.
Since the del Pezzo surface $V $ and the field extension $F\subset K$ are adapted, the scheme $C_\red$ is geometrically reduced.
It thus suffices to verify that $C$ is reduced.
This is obvious if the ramification curve $R$ is a line.
Assume now that $R$ is a split ribbon, 
so we are in   Case \refReid{P,2H} or Case \refReid{cS,2H} in Reid's Classification.
Let us  treat both   at the same time, allowing $e=1$ for $X=\PP^2$.
Seeking a contradiction, we now assume that $C$ is not reduced.
According to Proposition \ref{ramification multiple}, the morphism
$R_\red\ra C_\red $ has degree two.
To proceed, consider again the exact sequence of N\'eron--Severi groups
$$
0\lra \NS(Y)\lra\NS(X)\oplus\NS(C)\stackrel{\Psi}{\lra} \NS(R).
$$
According to \cite{Bosch; Luetkebohmert; Raynaud 1990}, Proposition 9.1.5
for each proper geometrically irreducible curve $Z$ the restriction map  $\NS(Z)\ra\ZZ$ given by $\shL\mapsto \deg(\shL|Z_\red)$
is injective. Choosing identifications $\NS(C)=\ZZ$ and $\NS(R)=\ZZ$,
we see that the linear map in the above exact sequence becomes the matrix $\Psi=(1,-2)$, because $(eH\cdot R_\red)=1$ 
and $R_\red\ra C_\red $ has degree two.
By Reid's Classification, the invertible sheaf $\nu^{-1}(\omega_Y)$ is given by $-eH$,
so the corresponding element in $\NS(X)\oplus\NS(C)$ is of the form $-(1,m)^t$. 
Since this lies in $\Kernel(\Psi)$, we must have   $m=1/2$, contradiction. 
\qed

\medskip
Next, we restrict the possible characteristics:

\begin{proposition}
\mylabel{double line}
If the ramification divisor  $R$ is a  split ribbon,  then the characteristic must be   $p\leq 3$
and the induced morphism $\nu:R_\red\ra C$ is birational.
If $p=3$, the normal surface $X$ is either $\PP^2$ or $\PP(1,1,3)$.
In case $p=2$, the only possibilities are  $\PP^2$, $\PP(1,1,2)$ or $\PP(1,1,4)$.
\end{proposition}

\proof
Clearly, the  tangent sheaf $\Theta_{V/F}=\underline{\Hom}(\Omega^1_{V/F},\O_V)$ on the normal surface $V$ is locally free in codimension one,
whence locally free at almost all points $a\in V$. In turn, $\Theta_{Y/K}$ is locally free
at almost all points $y\in Y$. Since $C_\red$ is geometrically reduced, we can use Bertini Theorems to find
a closed point $y\in C_\red$ in the regular locus so that the field extension $K\subset\kappa(y)$ is separable
(for example \cite{EGA V}, Proposition 4.3 or \cite{Jouanolou 1983}, Chapter I, Theorem 6.3).
Now we can apply Proposition \ref{ramification multiple characteristic} and deduces $p\leq 3$.

In Reid's Classification, the ramification divisor $R\subset X$ is linearly equivalent to a multiple
of the same effective divisor only in case (i) and (iii).
In other words, $X=\PP(1,1,e)$ is a contracted Hirzebruch surface
with numerical invariant $e\geq 1$, and ramification divisor $R=2H$.
Consider the divisorial part $N\subset\Sing(V/F)$ and its preimage $N_K\subset  Y$.
If $p=3$, according to Proposition \ref{ramification multiple equation}, 
the resulting Cartier divisor on $X$ is  $\nu^{-1}(N_K) =3H$. Since $\Pic(X)$ is freely generated by $eH$, we must have $e|3$. 
If $p=2$, Proposition \ref{ramification multiple equation} implies that $\nu^{-1}(N_K) =4H$,
therefore $e|4$. The assertion follows.
\qed

\begin{proposition}
\mylabel{line}
If the ramification divisor $R$ is the projective line, then the
characteristic must be $p=2$,    and the gluing map $\nu: R\ra C$ is isomorphic
to a   finite flat morphism  $\PP^1\ra \PP^1$ that is radical of degree two.
Furthermore, the inclusion $\O_C\subset\nu_*(\O_R)$ is a direct summand of $\O_C$-modules.
The normal surface $X$ is either $\PP^2$, $\PP(1,1,2)$  or $\PP(1,1,4)$.
\end{proposition}

\proof 
From Proposition \ref{ramification regular} we deduce $p=2$, that the ramification curve $C$ is regular,
and that  $\nu:R\ra C$ is finite, flat and radical.
Since this   holds true if we replace $K$ by any finite purely inseparable field extension,
the curve $C$ must be smooth. By L\"uroth's Theorem, it is  a projective line $\PP^1$.

The cokernel for the inclusion $\O_C\subset\nu_*(\O_R)$ is an invertible sheaf $\O_{\PP^1}(n)$
for some integer $n$, with Euler characteristic $\chi(\O_R)-\chi(\O_C)=0$. Thus $n=-1$.
Using $\Ext^1(\O_{\PP^1}(-1),\O_{\PP^1})= H^1(\PP^1,\O_{\PP^1}(1))=0$,
we infer that the inclusion $\O_C\subset\nu_*(\O_R)$ is a direct summand.

According to Reid's Classification, $X=\PP(1,1,e)$ is a contracted Hirzebruch surface with invariant $e\geq 1$,
and $\O_X(R)=\O_X(2H)$. According to Proposition \ref{ramification regular equation}, we have $\nu^{-1}(N_K)=2R$, and thus $e|4$.
\qed

\medskip
Note that the gluing map $\nu:R\ra C$ coincides with the Frobenius map on $\PP^1$ if the ground field $F$ is perfect, such that the function field $k(C)$ has
$p$-degree one. Combining the preceding results, we already obtain a non-existence result,
which generalizes \cite{Patakfalvi; Waldron 2017}, Corollary 1.4:

\begin{theorem}
\mylabel{higher p}
Normal del Pezzo surfaces that  are locally factorial, geometrically integral but
not geometrically normal do not exist in characteristic $p\geq 5$.
\end{theorem}

\medskip
Recall that $N\subset\Sing(V/F)$ denotes  the divisorial part of the locus of non-smoothness.
Its preimage $N_K\subset Y$ has the same support as the conductor curve $C\subset Y$.
According to Proposition \ref{reduced C}, we actually have $(N_K)_\red = C$.
It follows that $N$ is irreducible, and that $D=N_\red$ is integral.

\begin{proposition}
\mylabel{edim and e for p=2}
Suppose  $p=2$. Then the    geometric generic embedding dimension
and  geometric generic Hilbert--Samuel multiplicity  for the curve $N$ are  
$$
\edim(\O_{N,\eta}/F)=2\quadand e(\O_{N,\eta}/F)=4.
$$
The  preimage on normal surface  $X$ is given by   $\nu^{-1}(N_K)=2R$. 
If $N$ is not integral, then $N=2D$, and  $D=N_\red$ satisfies 
$$
\edim(\O_{D,\eta}/F)= 1 \quadand
e(\O_{D,\eta}/F)=  2.
$$
\end{proposition}

\proof 
Suppose that the ramification divisor $R\subset X$ is a line. By Proposition \ref{line}, 
the same holds for the conductor curve $C\subset Y$. Choose a rational point $x\in R$.
Its image $y\in C$ is rational as well. Thus we may apply
Theorem \ref{ramification regular equation}, and the values
for $\edim(\O_{N,\eta}/F)$ and $e(\O_{N,\eta}/F)$ and the multiplicity in $\nu^{-1}(N_K)=2R$
follow.
If the ramification locus $R$ is a double line, we argue analogously with Theorem \ref{ramification multiple equation}.
The case described in Theorem \ref{conductor multiple equation} does not occur, thanks to Proposition \ref{reduced C}. 

Now suppose that the curve $N$ is non-reduced, with $N=nD$. 
The multiplicities satisfy $e(\O_{N,\eta}/F)= n\cdot e(\O_{D,\eta}/F)$, hence $n|4$.
If $n=4$, the integral curve $D$ would be generically smooth, hence
$V$ would be smooth at some point contained in $D$. Contradiction. 
Thus we have $n=2$ and $e(\O_{D,\eta}/F)=  2$.
The value    $\edim(\O_{D,\eta}/F)=1$ again follows from the descriptions in 
Theorem \ref{ramification regular equation} and Theorem \ref{ramification multiple equation}.
\qed

\medskip
In the same way, one verifies the case of characteristic three:

\begin{proposition}
\mylabel{edim and e for p=3}
Suppose  $p=3$. Then the    geometric generic embedding dimension
and  geometric generic Hilbert--Samuel multiplicity  for the curve $N$ are  
$$
\edim(\O_{N,\eta}/F)=1\quadand e(\O_{N,\eta}/F)=3.
$$
The  preimage on  the normal surface $X$ is given by   $\nu^{-1}(N_K)=3R_\red$,
and the curve $N$ is integral.
\end{proposition}

\medskip
According to Proposition \ref{reduced C}, the conductor curve $C\subset Y$ coincides
with the reduced divisorial part of $\Sing(Y)$.
It will be crucial to understand the   effective Cartier divisors on the non-normal surface $Y$ supported by $C$.
For every invertible sheaf $\shL$ on $Y$, the short exact sequence  
$0\ra\shL\ra \shL_X\oplus\shL_C\ra\shL_R\ra 0$ of coherent sheaves on $Y$ 
yields a long exact sequence
\begin{equation}
\label{section sequence}
0\lra H^0(Y,\shL)\lra H^0(X,\shL_X)\oplus H^0(C,\shL_C)\lra H^0(R,\shL_R).
\end{equation}
In other words, the global sections of $\shL$ can be regarded  as pairs $s=(s_X,s_C)$ with $s_X|R=s_C|R$.

\begin{proposition}
\mylabel{sections vanishing}
The mapping
$$
H^0(X,\shL_X(-R))\lra H^0(Y,\shL),\quad s_X\longmapsto (s_X,0)
$$
is a bijection between the global section of the reflexive sheaf  $\shL_X(-R)$ and
the global sections of the invertible sheaf  $\shL$ that vanish along $C=N_\red$.
\end{proposition}

\proof
The  global sections $s=(s_X,s_C)$ that vanish along the conductor curve   have $s_C=0$, and hence $s_X|R=0$.
The assertion follows using the exact sequence $0\ra \shL(-R)\ra\shL_X\ra \shL_R\ra 0$.
\qed

\section{Smooth ramification loci}
\mylabel{Smooth ramification}

Let $F$ be a ground field of characteristic $p>0$,
and $V$ be a normal  del Pezzo surface that is locally factorial, geometrically integral but 
geometrically non-normal. After making a finite separable ground field extension,
we may assume that there is a finite purely inseparable field extension $F\subset K$
so that $V\ra\Spec(F)$ and $F\subset K$ are adapted. Let $\nu:X\ra Y$
be the normalization, $C\subset Y$ be the conductor curve,
and $R\subset X$ the ramification divisor. We then have a commutative diagram
$$
\begin{CD}
R	@>>> 	X\\
@VVV		@VV\nu V\\
C	@>>>	Y.
\end{CD}
$$
\emph{In this section, we treat the case that the ramification locus $R$ is smooth.}
By Proposition \ref{line},  we have $p=2$ and $R=\PP^1$,
and the gluing map $R\ra C$ is isomorphic to  a finite flat morphism $\PP^1\ra\PP^1$ that is radical of degree two.
The possibilities for the embedding $R\subset X$ are given by Reid's Classification in 
Theorem \ref{reid's classification}. Furthermore, we have a proper birational morphism $S\ra X$ from
a Hirzebruch surface $S$ with numerical invariant $e\geq 0$.
Let us start with a vanishing result:

\begin{proposition}
\mylabel{vanishing}
We have $H^1(V,\omega_V^{\otimes n})=0$ for all integers $n$.
\end{proposition}

\proof
By Serre duality, it suffices to check this for $n\geq 1$.
Set $\shL=\omega_V^{\otimes n}$. 
The short exact sequence $0\ra\shL_Y\ra\shL_X\oplus\shL_C\ra\shL_R\ra0$
yields an exact sequence
$$
H^0(R,\shL_R)\lra H^1(Y,\shL_Y)\lra H^1(X,\shL_X)\oplus H^1(C,\shL_C)\lra H^1(R,\shL_R).
$$
The term on the left vanishes, because  $\shL_R=\O_{\PP^1}(-m)$ for some integer $m\geq 1$. 
According to Proposition \ref{line}, the inclusion $\O_C\subset\nu_*(\O_R)$ is a direct summand, and it follows that $H^1(C,\shL_C)\ra H^1(R,\shL_R)$
is injective.
In turn, the canonical mapping  $H^1(Y,\shL_Y)\ra H^1(X,\shL_X)$ is injective as well.
From Reid's Classification in Theorem \ref{reid's classification}, one sees   that $\shL_X\simeq\O_X(-A)$
for some effective Cartier  divisor $A\subset X$ that is geometrically connected and geometrically reduced,
which ensures $h^0(\O_A)=1$. The normal surface $X$ is either the projective plane, a contracted Hirzebruch surface
or a Hirzebruch surface, which all have $H^1(X,\O_X)=0$.
The exact sequence $0\ra\shL_X\ra \O_X\ra\O_A\ra 0$ yields an exact sequence
$$
H^0(X,\O_X)\lra H^0(A,\O_A)\lra H^1(X,\shL_X)\lra H^1(X,\O_X),
$$
hence the term  $H^1(X,\shL_X)$ vanishes.
\qed

\medskip
Let $N\subset\Sing(V/F)$ be the divisorial part of the  locus of non-smoothness.
This irreducible curve  has geometric generic embedding dimension 
$\edim(\O_{N,\eta}/F)=2$, and the preimage on the normal surface $X$ is given by 
$\nu^{-1}(N_K)=2R$, by Proposition \ref{edim and e for p=2}.
The idea now is to study the integral curve $D=N_\red$:

\begin{proposition}
\mylabel{D smooth R}
The  possibilities for the integral curve $D=N_\red$ and the numerical invariant $e\geq 0$ of the
Hirzebruch surface $S$ are given by the following table,
according to the five cases (i)--(v) in Reid's Classification:
$$
\begin{array}[t]{l|*{12}{|c} }
\text{\rm Case}			& \multicolumn{2}{c|}{\refReid{P,2H}}	& \refReid{P,H}	& \multicolumn{3}{c|}{\refReid{cS,2H}} 			& \multicolumn{3}{c|}{\refReid{S,E+F}}	 		& \multicolumn{3}{c}{\refReid{S,E}}\\
\hline				& \multicolumn{2}{c|}{}			&		& \multicolumn{3}{c|}{}					& \multicolumn{3}{c|}{}			 		& \multicolumn{3}{c}{}\\[-2ex]
X				& \multicolumn{2}{c|}{\PP^2}		& \PP^2		& \multicolumn{3}{c|}{\PP(1,1,e)}			& \multicolumn{3}{c|}{S}		 		& \multicolumn{3}{c}{S}\\
\hline				& \multicolumn{2}{c|}{}			&		& \multicolumn{2}{c|}{}	 &				& \multicolumn{2}{c|}{}			& 		& 	&	& \\[-2ex]
e				& \multicolumn{2}{c|}{1}		& 1		& \multicolumn{2}{c|}{2} & 4				& \multicolumn{2}{c|}{0}		& 1		& 0	&1	& 2\\
\O_X(R)				& \multicolumn{2}{c|}{2H}		& H		& \multicolumn{2}{c|}{2H}& 2H				& \multicolumn{2}{c|}{E+F}		& E+F		& E	& E	& E\\
\hline&&&&&&&&&&&\\[-2ex]
h^0(\O_{D})			& 1 		& 1			& 1		& 1		& 1	& 1				& 1		& 1			& 1		& 1	& 4	& 2\\
h^1(\O_{D})			& 2		& 7			& 1		& 1		& 3	& 1				& 1		& 3			& 1		& 0	& 0	& 0\\
\edim(\O_{D,\eta}/F)		& 1		& 2			& 2		& 1		& 2	& 2				& 1		& 2			& 2		& 1	& 2	& 1\\
D				& \frac{1}{2}N	& N			& N		& \frac{1}{2}N	& N	& N				& \frac{1}{2}N	& N			& N		& \frac{1}{2}N  & N	& \frac{1}{2}N\\
\nu^{-1}(D_K)			& R		& 2R			& 2R		& R		& 2R	& 2R				& R		& 2R			& 2R		& R 	 & 2R 	& R\\
\end{array}
$$
\end{proposition}

\proof
According to Proposition \ref{edim and e for p=2}, we either have $D=\frac{1}{2}N$ or $D=N$, and the respective values for 
the geometric generic embedding dimension $d=\edim(\O_{D,\eta}/F)$ are  $d=1$ or $d=2$. 
We first treat the cases \refReid{P,2H}--\refReid{cS,2H} from Reid's Classification. 
According to Proposition \ref{line}, we have $e=2^\nu$ with $0\leq \nu\leq 2$.
Moreover,  $\Pic(V)$ is cyclic, and generated by the dualizing sheaf, so we have $\O_V(-D)=\omega_V^{\otimes m}$ for some integer $m\geq 1$.
The    short exact sequence $0\ra\omega_V^{\otimes m}\ra \O_V\ra \O_D\ra 0$ yields an exact sequence
$$
H^0(V,\O_V)\lra H^0(D,\O_D)\lra H^1(V,\omega_V^{\otimes m}).
$$
The term on the right vanishes, by Proposition \ref{vanishing},   thus $h^0(\O_D)=1$.
Riemann--Roch gives $2h^1(\O_D)=2-2\chi(\O_D)=2+\deg(\omega_D)$. By the Adjunction Formula, this integer equals  
$$
2+(K_V+D)\cdot D= 2+(\nu^{-1}(K_Y) + \nu^{-1}(D_K))\cdot \nu^{-1}(D_K),
$$
and the tabulated values for $h^1(\O_D)$   easily follow in each of the three cases \refReid{P,2H}--\refReid{cS,2H}.

In case \refReid{P,H},   $D=\frac{1}{2}N$ does not occur. Here $X=\PP^2$, and the ramification locus $R\subset X$ is a line.
In light of the exact sequence \eqref{neron-severi exact sequence},
we may regard $\Pic(Y)$ as the kernel of the matrix $\Psi=(1,-2)$.
The first matrix entry is the intersection number $(H\cdot R)=1$,
and the second entry is the negative   degree for the gluing map $\nu\colon R\ra C$.
Since $\nu^{-1}(N_K)=2R$, we may regard $N_K$ as the    vector $(\begin{smallmatrix}2\\1\end{smallmatrix})$,
which is primitive. Thus $N_K\in\Pic(Y)$ and hence $N\in\Pic(V)$ are generators, and
therefore $N=D$.  In case \refReid{cS,2H}, $D=\frac{1}{2}N$ with $X=\PP(1,1,4)$ does not occur for similiar reasons.

Next, we analyze case \refReid{S,E+F}.
Here $X=S$ is a Hirzebruch surface with numerical invariant $e\geq 0$,
and the ramification curve $R\subset X$ is linearly equivalent to $E+F$.
Seeking a contradiction, we assume $e\geq 2$. Then $E\cdot (E+F)=-e+1<0$.
The exact sequence $0\ra\O_X(F)\ra\O_X(E+F)\ra\O_E(E+F)\ra 0$
gives an exact sequence
$$
0\lra H^0(X,\O_X(F))\lra H^0(X,\O_X(E+F))\lra H^0(E,\O_E(E+F)).
$$
The term on the right vanishes, and the Projection Formula applied to the ruling $r:X\ra\PP^1$
gives $H^0(X,\O_X(F))=H^0(\PP^1,\O_{\PP^1}(1))$, whence $R$ must be a pair of lines, contradiction.
Thus   $e=0$ or $e=1$.
 
If $e=0$ then $X=S=\PP^1\times\PP^1$, and  the Picard group $\Pic(Y)$ can be regarded 
as the kernel of the matrix $\Psi=(1,1,-2)$. The first two entries are the intersection numbers
$(E\cdot R)=E\cdot (E+F)=1$ and $(F\cdot R)=F\cdot (E+F)=1$,
and the last entry is the negative degree of the gluing map $R\ra C$. Since $\nu^*(\omega_Y)=\O_X(-E-F)$ and $\nu^{-1}(N_K)=2R$,
the dualizing sheaf $\omega_Y$ is given by the column vector $-(1,1,1)^t$, 
whereas $N_K$ corresponds to the column vector $(2,2,2)^t$. Thus $\O_V(N)$ is a multiple of $\omega_V$,
and we compute the possible cohomological invariants for $D$ as above.
 
Now suppose the numerical invariant is  $e=1$,
that is, $X=S$ is the projective plane blown-up   in a rational point.
Then $\Pic(Y)$ is the kernel
of the matrix $\Psi=(0,1,-2)$. The first two entries are the intersection numbers $(E\cdot R)=E\cdot(E+F)=0$ and $(F\cdot R)=F\cdot (E+F)=1$,
and the last entry is the negative degree of the gluing map $R\ra C$.
Now $\nu^*(\omega_Y)=\O_X(-E-2F)$ and $\nu^{-1}(N_K)=2R$,
so  the  dualizing sheaf $\omega_Y$ corresponds to the column vector $-(1,2,1)^t$,
whereas the locus of non-smoothness $N_K$ gives $(2,2,1)^t$. The latter is primitive,
so  $D=N$. 
To proceed, consider the invertible sheaf  $\shL=\omega_V(N)$. Then $\shL_Y$ corresponds to the column vector $(1,0,0)^t$, thus $\shL_X=\O_X(E)$.
The Adjunction Formula gives
\begin{align*}
\deg(\omega_D) 	= (K_V+D)\cdot D=   E\cdot (2E+2F) = 0. 
\end{align*}
Therefore $h^0(\O_D)=h^1(\O_D)$.  To compute this  number, consider the short exact
sequence $0\ra \O_Y(-N_K)\ra \O_Y\ra\O_{N_K}\ra 0$. It gives a long exact sequence
$$
H^1(Y,\O_Y)\lra H^1(N_K,\O_{N_K})\lra H^2(Y,\O_Y(-N_K))\lra H^2(Y,\O_Y),
$$
in which the outer terms vanish, by Proposition \ref{neron-severi sequence}. Moreover,  $h^2(\O_Y(-N_K))=h^0(\shL_Y)$ by Serre Duality.
Since $\shL_X=\O_X(E)$ we have $\shL_R=\O_R$ and   $\shL_C=\O_C$. As  $E\subset X$ is negative-definite,
so $h^0(\O_X(E))=1$. Using the exact sequence 
\eqref{section sequence}, we infer $h^0(\shL_Y)=1$, therefore $h^1(\O_{N_K})=1$ and   $h^1(\O_D)=1$.
 
We come to the final case \refReid{S,E}.
Here $X=S$ is a Hirzebruch surface with numerical invariant $e\geq 0$, and the 
ramification locus is $R=E$. 
Now we may regard $\Pic(Y)$ as the kernel of the matrix $(-e,1,-2)$.
The first two entries are the intersection numbers $(E\cdot E)=-e$ and $(F\cdot E)=1$,
and the last entry is the negative degree of the gluing map $R\ra C$.
The dualizing sheaf $\omega_Y$ is given by the column vector $-(1,e+2,1)^t$,
whereas $N_K\subset Y$ corresponds to $(2,0,-e)^t$. Note that we have $N^2=(2E)^2=-4e$.

First suppose that the integer $e\geq 0$ is odd.
Then the column vector $(2,0,-e)^t$ is primitive, 
and it follows that   $D=N$. We compute
$$
\deg(\omega_D)= (K_Y+N_K)\cdot N_K = (E-(e+2)F)\cdot 2E = -2(2e+2).
$$
Then  $h^1(\O_D)=0$ and $h^0(\O_D)=2e+2$, by Riemann--Roch.
Hence the purely inseparable  field extension $F\subset H^0(D,\O_D)$ has degree $2e+2$.
Proposition \ref{edim and e for p=2} gives $4=e(\O_{D,\eta}/F)\geq h^0(\O_D)=2e+2$,
hence $e=1$ and $h^0(\O_D)=4$.

Second we assume that $e\geq 0$ is even: the column vector $(2,0,-e)^t$ attached to the Cartier divisor $N_Y\subset Y$ is therefore divisible.
Let $\shL$ be the invertible sheaf on $V$ so that $\shL_Y$ corresponds to the primitive vector $(1,0,-e/2)^t$.
Then we have $\shL_X=\O_X(E)$ and   $\shL_R=\O_{\PP^1}(-e)$ and $\shL_C=\O_{\PP^1}(-e/2)$.
Using the exact sequence \eqref{section sequence} we infer that $h^0(\shL)=1$, and any non-zero
section $s\in H^0(V,\shL)$ defines the unique effective Cartier divisor $D=\frac{1}{2}N$.
It follows $D^2=-e$, and the  Adjunction Formula gives
$$
\deg(\omega_D) = (K_Y + D_K)\cdot D_K = -(e+2)F\cdot E = -2(e/2+1)<0,
$$
thus $h^1(\O_D)=0$ and $h^0(\O_D)=e/2+1$.
On the other hand, Proposition \ref{edim and e for p=2} gives $2=e(\O_{D,\eta}/F)\geq h^0(\O_D)=e/2+1$.
The only solutions are $e=0$ or $e=2$, such that $h^0(\O_D)=1$ and $h^0(\O_D)=2$, respectively.
\qed
 
\medskip
So far, we have exploited that the local rings $\O_{V,a}$ are factorial.
If they are regular, one can say even more:

\begin{proposition}
\mylabel{D regular}
Suppose $D=\frac{1}{2}N$. If the  del Pezzo surface $V$ is regular,   the curve $D$ is regular as well.
\end{proposition}

\proof
Consider the   finite morphism   $\varphi:X\ra V$
obtained by composing  the normalization $X\ra Y=V_K$ with the projection $V_K\ra V$.
It is   flat, because  $V$ is regular  and $X$ is Cohen--Macaulay  
(for example  \cite{EGA IVb}, Proposition 6.1.5). We have $\nu^{-1}(N)=2R$, thus $\nu^{-1}(D)=R$.
Since here the scheme $R$ is regular, the scheme $D$ must be regular as well, by  \cite{EGA IVb}, Proposition 2.1.13.
\qed

\medskip
Now we  relate the geometric information provided by Proposition \ref{D smooth R} to 
the arithmetic of the ground field $F$. 

\begin{theorem}
\mylabel{main smooth ramification}
Suppose the del Pezzo surface $V$
is regular and the ground field $F$ has $\pdeg(F)=1$. Then $V$ has Picard number $\rho(V)=2$ and belongs to  case \refReid{S,E} from Reid's Classification,
with $e=2$.
In particular, $X=S$ must be a Hirzebruch surface with invariant $e=2$, and the ramification divisor is $R=E$.
\end{theorem}

\proof
We go through the table in Proposition \ref{D smooth R}.
Since $\pdeg(F)\leq 1$ the integral curve $D$ has $\edim(\O_{D,\eta}/F)\leq 1$,
according to Theorem \ref{edim bound}. This rules out all cases with $D=N$.
So $D=\frac{1}{2}N$, and the integral curve $D$ is regular by Proposition \ref{D regular}.
If furthermore $h^0(\O_D)=1$ then $D$ is geometrically reduced,
by \cite{Schroeer 2010}, Theorem 2.3. But this contradicts Corollary \ref{prime divisors geom non-reduced}
(compare also Proposition \ref{edim and e for p=2}).
The only possibility left is case \refReid{S,E}, with $e=2$.
From Proposition \ref{neron-severi sequence} we infer that $\rho(Y)=2$.
Since $Y\ra V$ is a universal homeomorphism, we also have  $\rho(V)=2$.
\qed

\medskip
Now suppose we are in case \refReid{S,E} of Reid's Classification, which is the only
possibility if the ground field has $\pdeg(F)=1$.
The Hirzebruch surface $X=S$ comes with a ruling $r:X\ra\PP^1$,
together with  the contraction $X\ra\PP(1,1,2)$  with exceptional divisor $E=R$.
The latter factors over the non-normal surface $Y$, but the former does not.
However, if we write $\psi:\PP^1\ra\PP^1$ for the composition of 
the gluing map $\nu:R\ra C$ and the inverse of $R=E\subset S=X\ra \PP^1$,
we get a commutative diagram
$$
\begin{CD}
X	@>\nu>>		Y\\
@VVV			@VVV\\
\PP^1	@>>\psi>	\PP^1.
\end{CD}
$$
Over each geometric point $\bar{a}:\Spec(\Omega)\ra\PP^1$, the fibers, $X_{\bar{a}}=\PP^1_{\Omega[\epsilon]}$ and $Y_{\bar{a}}$ sit in 
the conductor square
$$
\begin{CD}
\Spec(\Omega[\epsilon])	@>>>	\PP^1_{\Omega[\epsilon]}\\
@VVV				@VVV\\
\Spec(\Omega)		@>>>	Y_{\bar{a}},
\end{CD}
$$
where $\epsilon$ is an indeterminate modulo $\epsilon^2=0$.
In other words, all geometric fibers $Y_{\bar{a}}$ are   split ribbons $\PP^1\oplus\O_{\PP^1}(-1)$.
Let us  call  such morphisms \emph{quasirulings}.
In other words, these are   conic bundles where all fibers are geometrically non-reduced.

Since the projection $Y=V\otimes_FK\ra V$ is a universal homeomorphism,
the two morphism $Y\ra\PP(1,1,2)$ and $Y\ra\PP^1$ are induced from   morphisms 
$$
c:V\lra W\quadand f:V\lra\PP^1.
$$
The latter indeed goes to the projective line rather than a twisted form,
because $\Pic(V)=\Pic(Y)$. The image $W$ of the former is a twisted form of the contracted Hirzebruch surface,
with $W\otimes_FK=\PP(1,1,2)$. In particular, $W$ is a del Pezzo surface. 

\begin{proposition}
\mylabel{exceptional divisor}
The del Pezzo surface $W$ is regular, and the exceptional divisor 
for the contraction $c:V\ra W$ is $D=\frac{1}{2}N$, which is  a projective line $\PP^1_{F'}$ over some
purely inseparable quadratic field extension $F\subset F'$.
\end{proposition}

\proof
According to Proposition \ref{D smooth R}, the exceptional divisor must be the integral curve
$D=\frac{1}{2}N$, which has $h^0(\O_D)=2$ and $h^1(\O_D)=0$. 
In turn, $F'=H^0(D,\O_D)$ is a purely inseparable field extension of degree two.
Moreover, the curve $D$ is regular, by Proposition \ref{D regular}.
Since $D^2$ coincides with $R^2=-e=-2$, the scheme $D$ contains an $F'$-rational point, and it follows that $D=\PP^1_{F'}$.
According to Castelnuovo's Contraction Theorem, \cite{Shafarevich 1966}, Chapter 6, page 102, the curve $D\subset V$ is an exceptional divisor of the first kind, and  $W$ is regular.
\qed

\medskip
We shall see in Section \ref{Del Pezzo fibrations} that such regular del Pezzo surfaces $W$ and $V$ actually do exist
over any imperfect ground field $F$ of characteristic $p=2$.

Note that the quasiruling $f:V\ra\PP^1$ and the contraction $c:V\ra W$ exists even without the assumption
that $V\ra\Spec(F)$ and $F\subset K$ are adapted: If $F\subset F'$ is a Galois extension,
the Galois action on $\Pic(V\otimes_FF')$ must be trivial, because it respects the cone of curves,
and the two extremal rays have different self-intersection.
The only difference is that the range of the  quasiruling may  be a twisted form of the projective line,
rather than the projective line. Such twisted forms correspond to elements of order two
in the Brauer group $\Br(F)$.

\section{Non-reduced ramification locus in characteristic two}
\mylabel{Non-reduced ramification p=2}

In this section, we study the situation  over  ground fields $F$ of characteristic $p=2$
where the ramification divisor $R\subset X$ is non-smooth, hence a split ribbon $R=\PP^1\oplus\O_{\PP^1}(-1)$.
Notation is as in Section \ref{Smooth ramification}: The normal del Pezzo surface
$V$   is locally factorial, geometrically integral but geometrically non-normal. 
We assume that the structure morphism $V\ra\Spec(F)$ and the field extension $F\subset K$ are adapted,
and consider the normalization $\nu:X\ra Y=V_K$ and the ensuing conductor square
$$
\begin{CD}
R	@>>> 	X\\
@VVV		@VV\nu V\\
C	@>>>	Y.
\end{CD}
$$
Recall that  we have a proper birational morphism $S\ra X$ from
a Hirzebruch surface $S$ with numerical invariant $e\geq 0$.
Furthermore,   $N\subset\Sing(V/F)$ denotes the divisorial part
of the locus of non-smoothness, and we write $D=N_\red$ for its reduction.

\begin{proposition}
\mylabel{ramification non-reduced p=2}
The possible numerical invariants
$h^i(\O_{D})\geq 0$ of the integral curve $D$, together with 
its  geometric generic embedding dimension $\edim(\O_{D,\eta}/F)$,
the relation to the divisor $N\subset V$ and the   preimage $\nu^{-1}(D_K)\subset X$ 
are given by the following table:
$$
\begin{array}[t]{l|*{10}{|c}}
\text{\rm Case}			& \multicolumn{3}{c|}{\refReid{P,2H}}		& \multicolumn{7}{c}{\refReid{cS,2H}} 					\\
\hline				
				& \multicolumn{3}{c|}{}				& \multicolumn{7}{c}{}							\\[-2ex]
X				& \multicolumn{3}{c|}{\PP^2}			& \multicolumn{7}{c}{\PP(1,1,e)}					\\
\hline				
				& \multicolumn{3}{c|}{}				& \multicolumn{4}{c|}{}							& \multicolumn{3}{c}{}\\[-2ex]
e				& \multicolumn{3}{c|}{1}			& \multicolumn{4}{c|}{2}						& \multicolumn{3}{c}{4}\\
\O_X(R)				& \multicolumn{3}{c|}{2H}			& \multicolumn{4}{c|}{2H}						& \multicolumn{3}{c}{2H}\\
\hline				
				&   &				& 		&   & 				&   & 					&   &   & 	\\[-2ex]
h^0(\O_{D})			& 1 & 2	 			& 1		& 1 & 2				& 1 & 2					& 1 & 2 & 4	\\
h^1(\O_{D})			& 7 & 8	 			& 2		& 3 & 4				& 1 & 2					& 1 & 2 & 4	\\
\hline				
				& \multicolumn{2}{c|}{}		&		& \multicolumn{2}{c|}{}		& \multicolumn{2}{c|}{}			& \multicolumn{3}{c}{}						\\[-2ex]
\edim(\O_{D,\eta}/F)		& \multicolumn{2}{c|}{2}	& 1		& \multicolumn{2}{c|}{2}	& \multicolumn{2}{c|}{1}		& \multicolumn{3}{c}{2}\\
D				& \multicolumn{2}{c|}{N}	& \frac{1}{2}N	& \multicolumn{2}{c|}{N}	& \multicolumn{2}{c|}{\frac{1}{2}N}	& \multicolumn{3}{c}{N}\\
\nu^{-1}(D_K)			& \multicolumn{2}{c|}{2R}	& R		& \multicolumn{2}{c|}{2R}	& \multicolumn{2}{c|}{R}		& \multicolumn{3}{c}{2R}
\end{array}
$$
\end{proposition}

\proof
According to Proposition \ref{double line}, only    case \refReid{P,2H},  and case \refReid{cS,2H} with $e=2$ or $e=4$ 
from Reid's Classification in Theorem \ref{reid's classification}
are possible. We may treat them 
simultaneously, by regarding the normalization as $X=\PP(1,1,e)$ with $e=2^\nu$ with $0\leq \nu\leq 2$.
In any case,  we have  $\nu^*(K_Y)= -eH$.
According to Theorem  \ref{ramification multiple equation} and Theorem \ref{conductor  multiple equation}, 
the preimage $\nu^{-1}(N_K)=2R$ is linearly equivalent to $4H$.
Without restriction, we may assume that $R=2H$, such that $R_\red=H$.
Using \eqref{neron-severi exact sequence}, we have an exact sequence of N\'eron--Severi groups
$$
0\lra \NS(Y)\lra\NS(X)\oplus\NS(C)\stackrel{\Psi}{\lra} \NS(R_\red).
$$
The induced morphism $\nu:R_\red\ra C_\red$ has degree $d\leq 2$,
and the group $\Pic(X)$ is freely generated by $eH\subset X$, which has $(eH\cdot R_\red)=(eH\cdot H)=1$.
Thus $\NS(Y)$ can be regarded as the kernel of the
matrix $\Psi=(1,-d)\in\Mat_{1\times 2}(\ZZ)$.
The numerical class of $K_Y$ is of the form 
$-(\begin{smallmatrix} a\\b\end{smallmatrix})$ with $a=1$ and $-a+db=0$. 
The existence of an integral solution forces $d=1$, so the morphism $\nu:R_\red\ra C_\red $ is birational.
Furthermore, the numerical class of $K_Y$
is given by $-(\begin{smallmatrix} 1\\1\end{smallmatrix})$.

The numerical class of $N_K$ is given by the vector
$\frac{4}{e}(\begin{smallmatrix} 1\\1\end{smallmatrix})$.
We have $H\subset\nu^{-1}(D_K)$, but equality is impossible, because otherwise $D_K$ and whence $D$ would be
geometrically reduced.
This yields, in each of the three cases $e=2^\nu$, the indicated two possibilities  
$D=\frac{1}{d}N$ with $d=1$ and $d=2$. Consequently,  we have 
\begin{equation}
\label{pullback divisors}
\nu^{-1}(D_K)=\frac{4}{de}eH\quadand \nu^{-1}(K_Y)=-eH.
\end{equation}
Since the coefficient $4/de$ is an integer, the case with $d=2$ and $e=4$ is impossible.
The Adjunction Formula, base-change and the Projection Formula yield
$$
\deg(\omega_D)=(K_V+D)\cdot D = (K_Y+D_K)\cdot D_K = (\nu^{-1}(K_Y)+\nu^{-1}(D_K))\cdot\nu^{-1}(D_K).
$$
Together with \eqref{pullback divisors} we infer
$$
\deg(\omega_D) = \left(-1+\frac{4}{de}\right)\frac{4}{de} (eH)^2 = \frac{4(4-de)}{d^2e}.
$$
This yields
$$
\chi(\O_D)=-\frac{1}{2}\deg(\omega_D) = 
\begin{cases}
-6	& \text{if $e=1$, $d=1$;}\\
-1	& \text{if $e=1$, $d=2$;}\\
-2	& \text{if $e=2$, $d=1$;}\\
0	& \text{else.}
\end{cases}
$$
Using that $h^0(\O_D)$ is a $p$-power   and divides $\chi(\O_D)$, we obtain the possibilities indicated in  the table.
\qed

\medskip
This narrows down further if we take the arithmetic of the ground field $F$ into account:

\begin{proposition}
\mylabel{pdeg ramification non-reduced}
If the ground field has $\pdeg(F)\leq 1$, then the del Pezzo surface $V$ belongs to case (iii)
of Proposition \ref{ramification non-reduced p=2}, and the curve $D\subset V$ has numerical invariants
$h^0(\O_D)=h^1(\O_D)=1$ or $h^0(\O_D)=h^1(\O_D)=2$.
\end{proposition}

\proof
Without restriction, we may assume that the ground field $F$ is separably closed.
Then the length of any  finite irreducible $F$-scheme   is a 2-power.
According to Theorem \ref{edim bound}, we must have $\edim(\O_{D,\eta}/F)\leq 1$,
which rules out all but three sub-cases in Proposition \ref{ramification non-reduced p=2}.
In these remaining cases, $D=\frac{1}{2}N$ is a proper integral curve
inside the del Pezzo surface
with geometric generic embedding dimension $\edim(\O_{D,\eta}/F)=1$.
Furthermore, the curve $D$ contains no $F$-rational points. 

In the first sub-case, the curve $D$ has 
numerical invariants $h^0(\O_D)=1$ and $h^1(\O_D)=2$.
We show that such a curve does not exist:
According to \cite{Schroeer 2010}, Theorem 2.3, the curve $D$ is non-normal.
Let $\tilde{D}\ra D$ be the normalization map. Since $\tilde{D}$ is normal and
$\edim(\O_{\tD,\eta}/F)=1$, the same result ensures that $h^0(\O_{\tilde{D}})=2^\nu$
for some exponent $\nu\geq 1$. In particular, $h^0(\O_{\tilde{D}})$ is an even number.
Now consider the conductor square
$$
\begin{CD}
A	@>>> 	\tilde{D}\\
@VVV		@VVV\\
B	@>>>	D.
\end{CD}
$$
From this we obtain a long exact sequence
$$
0\ra H^0(\O_D)\ra H^0(\O_{\tilde{D}})\oplus H^0(\O_B)\ra H^0(\O_A)\ra
H^1(\O_D)\ra H^1(\O_{\tilde{D}})\ra 0.
$$
The numbers $h^1(\O_{\tilde{D}})$ and $h^0(\O_A)$ are even, because $h^0(\O_{\tilde{D}})$ is even.
Since $D(F)=\varnothing$ and  $F$ is separably closed,
the number   $h^0(\O_B)$ is even as well. Thus all terms in the exact sequence have even vector space dimension over the
ground field $F$, except $H^0(\O_D)$, which is one-dimensional.
On the other hand,  the alternating sum of  vector space dimensions  in the exact sequence vanishes, 
contradiction.
\qed

\section{Non-reduced ramification locus in characteristic three}
\mylabel{Non-reduced ramification p=3}

In this section, we study the situation  over  ground fields $F$ of characteristic $p=3$.
Notation is as in Section \ref{Smooth ramification}: The normal del Pezzo surface
$V$   is locally factorial, geometrically integral but geometrically non-normal. 
We assume that the structure morphism $V\ra\Spec(F)$ and the field extension $F\subset K$ are adapted,
and consider the normalization $\nu:X\ra Y=V_K$ and the ensuing conductor square
$$
\begin{CD}
R	@>>> 	X\\
@VVV		@VV\nu V\\
C	@>>>	Y.
\end{CD}
$$
We saw in Proposition \ref{no pair of lines} and Proposition \ref{line} that the ramification divisor   is 
a split ribbon $R=\PP^1\oplus\O_{\PP^1}(-1)$.
Recall that  we have a proper birational morphism $S\ra X$ from
a Hirzebruch surface $S$ with numerical invariant $e\geq 0$.
Furthermore,   $N\subset\Sing(V/F)$ denotes the divisorial part
of the locus of non-smoothness, and we write $D=N_\red$ for its reduction.

\begin{proposition}
\mylabel{ramification non-reduced p=3}
Assumptions as above. Then   $N$ is reduced, and   the curve $D=N$ has 
$\edim(\O_{D,\eta}/F)=1$  and $\nu^{-1}(D_K)=\frac{3}{2}R$.
Only two cases in Reid's Classification are possible, and the  numerical invariants are given by
the following table:
$$
\begin{array}[t]{l|*{4}{|c}}
\text{\rm Case}			& \multicolumn{2}{c|}{\refReid{P,2H}}		& \multicolumn{2}{c}{\refReid{cS,2H}} 					\\
\hline				
				& \multicolumn{2}{c|}{}				& \multicolumn{2}{c}{}							\\[-2ex]
X				& \multicolumn{2}{c|}{\PP^2}		& \multicolumn{2}{c}{\PP(1,1,3)}\\
\hline				
				& \multicolumn{2}{c|}{}				& \multicolumn{2}{c}{}							\\[-2ex]
e				& \multicolumn{2}{c|}{1}			& \multicolumn{2}{c}{3}						  \\
\O_X(R)				& \multicolumn{2}{c|}{2H}			& \multicolumn{2}{c}{2H}						\\
\hline				
				&   &				& 		&  	\\[-2ex]
h^0(\O_{D})			& 1  & 3	 			&\hspace{.2cm}1\hspace{.2cm} 	& 3\\
h^1(\O_{D})			&4 & 6	 			&1& 3\\	
\end{array}
$$
\end{proposition}

\proof
The curve $N$ is reduced by Proposition \ref{edim and e for p=3}, with $\edim(\O_{N,\eta}/F)=1$.
We proceed by going through the cases in  Reid's classification.
Since the ramification divisor $R$ is the split ribbon,
only Case \refReid{cS,2H} with $e=3$ and Case \refReid{P,2H}  
from Theorem \ref{reid's classification}
are possible. We may treat them 
simultaneously, by regarding the normalization as $X=\PP(1,1,e)$ with $e=3$ and  $e=1$, respectively.

The split ribbon has $h^0(\O_R)=1$ and $h^1(\O_R)=0$,
and  the del Pezzo surface has  $h^2(\O_Y)=0$. Proposition \ref{neron-severi sequence} gives an exact sequence
$$
0\lra \NS(Y)\lra\NS(X)\oplus\NS(C)\stackrel{\Psi}{\lra} \NS(R ).
$$
Furthermore,  the Weil divisor $R_\red\subset X$
is linearly equivalent to the line $H\subset X$, and the induced morphism $\nu:R_\red\ra C$ is 
birational. Furthermore, the restriction map $\Pic(R)\ra\Pic(R_\red)$ is bijective.
The group $\Pic(X)$ is freely generated by $eH\subset X$, which has $(eH\cdot R_\red)=(eH\cdot H)=1$.
Thus $\NS(Y)$ can be regarded as the kernel of the
matrix $\Psi=(1,-1)\in\Mat_{1\times 2}(\ZZ)$.
Then $\nu^*(K_Y)= -eH$ generates the Picard group $\Pic(X)$,
and $\nu^{-1}(N_K)$ is linearly equivalent to $3H=\frac{3}{e}\cdot eH$.
The numerical class of $\omega_Y$ corresponds to the vector 
$-(\begin{smallmatrix} 1\\1\end{smallmatrix})$, whereas the numerical class of $N_K$ is given by the vector
$\frac{3}{e}(\begin{smallmatrix} 1\\1\end{smallmatrix})$. The Adjunction Formula for $D=N$ yields
$$
\deg(\omega_D)=(K_V+D)\cdot D = (K_Y+N_K)\cdot N_K = \left(-eH + \frac{3}{e}eH\right)\cdot \frac{3}{e}eH = 3(3/e-1).
$$
In    Case \refReid{P,2H} we have $e=1$ and thus  $\deg(\omega_D)=6$. Now we use that  $h^0(\O_D)$ is a $p$-power that
divides $\chi(\O_D)=-\frac{1}{2}\deg(\omega_D)=-3$, and get the indicated values for $h^i(\O_D)$.
In Case \refReid{cS,2H}  we have $e=3$,   thus $\deg(\omega_D)=0$ and $h^0(\O_D)=h^1(\O_D)$.
Here  $h^0(\O_D)$ is a $p$-power that divides the selfintersection number $D^2=(eH)^2=3$,
which again gives the possible values for $h^i(\O_D)$.
\qed

\medskip
We proceed to  rule out all the two sub-cases with $h^0(\O_D)>1$.
This relies on a general result of Maddock on normal del Pezzo surfaces with irregularity $h^1(\O_V)>0$:
According to \cite{Maddock 2016}, Corollary 1.4 the bound
\begin{equation}
\label{maddock bound}
h^1(\O_V)\geq \frac{p^2-1}{6} K_V^2 
\end{equation}
then holds. Let us check that this result indeed applies:

\begin{proposition}
\mylabel{char 3 irregularity}
Our del Pezzo surface $V$ has irregularity $h^1(\O_V)>0$ provided the curve $D=N$ has $h^0(\O_D)>1$.
\end{proposition}

\proof
Seeking a contradiction we assume  $h^1(\O_V)=0$ and thus $h^1(\O_C)=0$.
The short exact sequence $0\ra\O_V(-D)\ra\O_V\ra\O_D\ra 0$ yields a long exact sequence
\begin{equation}
\label{sequence p=3}
H^1(V,\O_V)\lra H^1(D,\O_D)\lra H^2(V,\O_V(-D))\lra H^2(V,\O_V),
\end{equation}
where the outer terms vanish. Together with Serre duality we get $h^1(\O_D)=h^0(\shL)$ for the invertible sheaf $\shL=\omega_V(D)$.
In case $X=\PP(1,1,3)$ we have $h^1(\O_D)=3$ and  $\shL_X=\O_X$, whence $\shL$ is numerically trivial and  $h^0(\shL)\leq 1$, contradiction.

Now suppose we are in case $X=\PP^2$. Then $\shL_X=\O_{\PP^2}(2)$. Since the morphism $R_\red\ra C$ is birational,
we have $\deg(\shL_C)=2$ and with Riemann--Roch $h^0(\shL_C)=3$. The short exact sequence $0\ra\shL_Y\ra\shL_X\oplus\shL_C\ra\shL_R\ra 0$
yields an exact sequence
$$
0\lra H^0(Y,\shL_Y)\lra H^0(X,\shL_X) \oplus H^0(C,\shL_C) \lra H^0(R,\shL_R).
$$
Since $\shL_X(-R)=\O_{\PP^2}$ we see that the restriction map $H^0(X,\shL_X) \ra H^0(R,\shL_R)$
is surjective, with one-dimensional kernel. So $h^1(\O_D)=h^0(\shL)=1+h^0(\shL_C)=4$.
Hence we are in the case of Proposition \ref{ramification non-reduced p=3} where  $h^0(\O_D)=1$, contradiction.
\qed

\begin{proposition}
\mylabel{no subcases p=3}
The two  sub-cases in Proposition \ref{ramification non-reduced p=3} with  $h^0(\O_D)>1$ do not exist.
\end{proposition}

\proof
The inclusion $C\subset D_K\subset Y$ yield a commutative diagram
\begin{equation}
\label{sequence diagonal}
\begin{gathered}
\begin{xy}
\xymatrix{
H^1(Y,\O_Y)\ar[dr]_\simeq\ar[r]	& H^1(D_K,\O_{D_K})\ar[d]\ar[r]		& H^2(Y,\O_Y(-D_K))\ar[r]	& 0\\
				& H^1(C,\O_C)
}
\end{xy}
\end{gathered}
\end{equation}
with exact row. The vertical map is surjective, and the diagonal map is bijective.
Serre duality gives $h^2(\O_Y(-D_K)) = h^0(\shL_Y)$ for the invertible sheaf $\shL=\omega_V(D)$.

First suppose that we are in case (iii) of Proposition \ref{ramification non-reduced p=3}. Then $K_V^2=3$. Moreover, $\shL$ is numerically trivial,
thus $h^0(\shL)\leq 1$,  and $h^1(\O_D)=3$. It follows that $h^1(\O_V)=h^1(\O_C)\leq 3$.
On the other hand, Maddock's bound \eqref{maddock bound} gives $h^1(\O_V)\geq 8/6\cdot K_V^2=4$, contradiction.

Finally suppose we are in case (i) of   Proposition \ref{ramification non-reduced p=3}. Then $K_V^2=1$, so Maddock's bound
ensures $h^1(\O_V)\geq 2$. On the other hand, the ring of global sections $F'=H^0(D,\O_D)$ is a purely inseparable field extension $F\subset F'$
of degree $[F':F]=3$, and the cohomology group $H^1(D,\O_D)$ is a two-dimension vector space over $F'$.
In turn, $H^1(D_K,\O_{D_K})$ is a free module of rank two over the local Artin ring $A'=F'\otimes_FK$, which has length two.
The inclusion $C\subset D_K$ of curves induces a surjection $A'=H^0(D_K,\O_{D_K})\ra H^0(C,\O_C)$ whose kernel 
is an ideal $\ideala\subset A'$ of length one. This ideal annihilates the image of the surjection $A'\oplus A'\simeq H^1(D_K,\O_{D_K})\ra H^1(C,\O_C)$,
and we conclude that $\dim_K H^1(C,\O_C)\leq 2$. In turn, our del Pezzo surface $V$ has irregularity $h^1(\O_V)\leq 2$.
Summing up, we have $h^1(\O_C)=h^1(\O_V)=2$.

Using the exact sequence in \eqref{sequence diagonal}, we conclude that $h^0(\shL)\geq 6-2= 4$.
On the other hand, the arguments in the proof for Proposition \ref{char 3 irregularity} show that $h^0(\shL)=1+h^0(\shL_C)$.
According to Proposition \ref{ramification multiple}, the  conductor curve $C$ is integral.
In our situation we have  $h^0(\O_C)=1$ and $h^1(\O_C)=2$, and $\deg(\shL_C)=2$.
Riemann--Roch gives $\chi(\shL_C)=1$. Furthermore, $\deg(\omega_C)=2$, hence $\omega_C\otimes\shL_C^\vee$
is numerically trivial, and $h^1(\shL_C)\leq 1$.
Combining these observations we get  $h^0(\shL_C)= \chi(\shL_C)+h^1(\shL_C)\leq 2$,
and therefore $h^0(\shL)\leq 3$, contradiction.
\qed

\begin{proposition}
\mylabel{vanishing char 3}
In case (iii) of Proposition \ref{ramification non-reduced p=3} with $h^1(\O_D)=1$,  the  del Pezzo surface has $H^1(V,\O_V)=0$.
\end{proposition}

\proof 
Seeking a contradiction, we assume that $h^1(\O_V)>0$. 
According to Maddock's bound \eqref{maddock bound}, the irregularity  $h^1(\O_V)$ is at least two.
Now recall that the restriction map $H^1(Y,\O_Y)\ra H^1(C,\O_C)$ is bijective.
The inclusion of curves $C\subset D_K$ induces a surjection $H^1(D_K,\O_{D_K})\ra H^1(C,\O_C)$ 
on cohomology groups, and we thus have $h^1(\O_D)\geq 2$, contradiction.
\qed

\section{Some peculiar genus-one curves}
\mylabel{Peculiar curves}

In this section we study some rather peculiar algebraic curves of arithmetic genus one.
These curves  showed up on certain geometrically non-normal del Pezzo surfaces
for $p=2$ and $p=3$, and their   structure will lead to non-existence results.
These  genus-one curves, however, are of interest in all characteristics.

Let $F$ be  an imperfect ground field of characteristic $p>0$,
and $F\subset\tilde{F}$ be a purely inseparable extension of degree $ [\tF:F]=p$.
Choose an element $\beta\in F^\times$ so that $\tF=F(\beta^{1/p})$.
Set $\tD=\PP^1_\tF$, write $a=(0:1)$ for the origin,
and let $A\subset\PP^1_{\tF}$ be its first infinitesimal neighborhood.
The structure sheaf $\O_A$ is a skyscraper sheaf supported by $a\in\tD$.
By abuse of notation, we also  write $\O_A$ when we mean the stalk $\O_{A,a}$ or
the ring of global sections $H^0(A,\O_A)$. This said, we may write
$$
\O_A=\tF[\epsilon] = \tF\oplus\tF\epsilon.
$$
As customary, $\epsilon$ denotes an  indeterminate  subject to $\epsilon^2=0$.
Consider the $F$-subalgebra $\O_B$ generated by $\beta^{1/p}+\epsilon\in\O_A$. 
Using $(\beta^{1/p}+\epsilon)^p=\beta+\epsilon^p=\beta$, we see that the resulting map
$\tilde{F}\ra\O_B$ given by  $\beta^{1/p}\mapsto \beta^{1/p}+\epsilon$ is bijective.
The crucial point here is that the images of
$H^0(\O_\tD)$ and $H^0(\O_B)$ define two different copies of the field $\tilde{F}$ inside
the $F$-algebra $H^0(\O_A)$, which become the same after projecting onto the residue field $\kappa(a)=\tF$. 
In other words, the construction relies on the  non-uniqueness of coefficient fields.
The inclusion $\O_B\subset\O_A$ yields a finite morphism $A\ra B$,
and the resulting   cocartesian square
\begin{equation}
\label{conductor square curve}
\begin{CD}
A	@>>>	\tD\\
@VVV		@VV\nu V\\
B	@>>>	D
\end{CD}
\end{equation}
defines a proper integral curve $D$ containing a unique singularity $b=\nu(a)\in D$
with residue field $\kappa(b)=\tF$. 
The normalization map $\nu:\tD\ra D$ is a universal homeomorphism, and
the cocartesian square yields a long exact sequence
\begin{equation}
\label{cohomology sequence}
0\ra H^0(\O_D)\ra H^0(\O_\tD)\oplus H^0(\O_B)\ra H^0(\O_A)\ra H^1(\O_D) \ra 0,
\end{equation}
from which one can deduce the structure of the curve:

\begin{proposition}
\mylabel{peculiar gorenstein}
The curve $D$ is Gorenstein and  has cohomological invariants $h^0(\O_D)=1$ and $h^1(\O_D)=1$.
Moreover, the reduced base-change $(D\otimes_F\tF)_\red$ is isomorphic to the projective line  $\PP^1_\tF$.
For $p=2$, the base-change $D\otimes_F\tF$ is the split ribbon on $\PP^1_\tF$ with 
respect to the ideal sheaf $\O_{\PP^1_\tF}(-2)$.
\end{proposition}

\proof
The length of $\O_A$ as an $\O_B$-module is two, hence $D$ is Gorenstein, by Proposition \ref{sufficient condition gorenstein}.
The conductor square \eqref{conductor square curve} gives
$$
\chi(\O_D) = \chi(\O_\tD) + \chi(\O_B) - \chi(\O_A) = p+ p-2p=0,
$$
and consequently $h^1(\O_D)=h^0(\O_D)$.
Consider the   field extensions $F\subset H^0(\O_D)\subset H^0(\O_\tD)=\tF$,
which has degree $p$. Suppose the first inclusion is not an equality.
Then the second inclusion must be an equality, and this also holds
for the inclusion $H^0(\O_D)\subset H^0(\O_B)$. By the exactness
of \eqref{cohomology sequence}, the two maps 
$$
H^0(\O_D)=H^0(\O_\tD)\lra H^0(\O_A) \longleftarrow H^0(\O_B)
$$
must have the same image. But  the  former map factors over the first summand in $\O_A=\tF\oplus \tF\epsilon$, whereas
the latter does not,  contradiction.
Thus we have  $H^0(\O_D)=F$, and hence $h^1(\O_D)=h^0(\O_D)=1$.

To understand the base-change,  write $\PP^1_\tF=\Proj\tF[U,V]$, and look at the affine open subset $D_+(V)=\Spec\tF[u]$, with $u=U/V$.
The conductor square \eqref{conductor square curve} yields the  cartesian diagram
\begin{equation}
\label{peculiar square}
\begin{CD}
\tF[u]/(u^2) 		@<<<	\tF[u]\\
@AAA				@AAA\\
F[\beta^{1/p}+u]	@<<<	R,
\end{CD}
\end{equation}
where $\Spec(R)$ defines an affine open neighborhood of the singularity $b\in D$.
Base-changing the diagram along $F\subset\tF$, we get a new cartesian diagram
$$
\begin{CD}
\tF\otimes_F\tF[u]/(u^2) 		@<<<	\tF\otimes_F\tF[u]\\
@AAA				@AAA\\
\tF\otimes_FF[\beta^{1/p}+u]	@<<<	\tF\otimes_FR,
\end{CD}
$$
Clearly, the birational morphism $\varphi:\PP^1_\tF=(\tD_\tF)_\red\ra (D_\tF)_\red$
is an isomorphism over the complement of the the point $\tilde{b}\in D_\tF$ corresponding to the singularity $b\in D$.
The schematic fiber $\varphi^{-1}(\tilde{b})\subset\PP^1_\tF$ is given by the spectrum $\tF\otimes_F\tF[u]/(u^2) $ modulo the ideal generated
by
$$
\beta^{1/p}\otimes1-1\otimes \beta^{1/p}\quadand \beta^{1/p}\otimes 1-1\otimes(\beta^{1/p}+\epsilon).
$$
Clearly, they generate the maximal ideal of the local Artin ring $\tF\otimes_F\tF[u]/(u^2) $,
and we conclude that $(D_\tF)_\red=\PP^1_\tF$. 

Finally, suppose $p=2$. Then the nilradical $\shI$ for the scheme $D_\tF$
is torsion-free of rank one, with $\shI^2=0$, whence $\shI=\O_{\PP^1}(n)$ for some integer $n$.
We have $n=-2$ because  $\chi(\shI)=-1$. Such a  ribbon is necessarily split, because
$\Ext^1(\Omega^1_{\PP^1},\O_{\PP^1}(-2)) =0$.
\qed

\medskip
The curve $D$ does not seem to have a nice description in terms of equations. However, 
the situation improves when we pass to   more singular models:

\begin{proposition}
\mylabel{birational model}
There is a finite birational universal homeomorphism $D\ra D'$
to the effective Cartier divisor $D'\subset \PP(1,1,p)$ in   weighted projective space defined by
the homogeneous equation $t^p-\beta(x^p-\beta y^p)^p=0$ of degree $p^2$.
In characteristic two, this birational morphism  is actually an isomorphism.
\end{proposition}

\proof
Let us continue to use the notation from the preceding proof.
From the cartesian square \eqref{peculiar square} we see that the subring $R\subset\tF[u]$ comprises all polynomials of the form
$$
\lambda(\beta^{1/p}+u) + \lambda_2u^2+\ldots+\lambda_nu^n
$$
with coefficients $\lambda\in F$ and $\lambda_i\in\tF$.
Consider the two particular polynomials 
\begin{equation}
\label{normalization substitution}
x=\beta^{1/p}+u		\quadand 	t=\beta^{1/p}u^p\in R,
\end{equation}
which satisfy the relation $\beta x^{p^2} = \beta^{p+1} + t^p$, hence
\begin{equation}
\label{inhomogeneous equation}
t^p-\beta(x^p - \beta)^p = 0.
\end{equation}
Regarding $x,t$ as indeterminates and using the Reduction Criterion for the field of fraction
$F[x,t]\subset F(x,t)$, we see that 
the polynomial $t^p-\beta(x^p - \beta)^p$ is irreducible. It follows that the  canonical map
$$
R'=F[x,t]/(t^p-\beta(x^p - \beta)^p) \lra R\subset\tF[u]
$$
is injective. To see that the resulting inclusion of function fields is bijective, consider as denominator $s=x^p - \beta$.
Then $(t/s)^p=\beta$, such that $(R')_s=\tF[x]=\tF[u]$.
We conclude that the morphism $\Spec(R)\ra\Spec(R')$ is a universal homeomorphism, and becomes
an isomorphism upon removal of the  singular point $b\in\Spec(R)$ and its image $b'\in \Spec(R')$.

Homogenization of \eqref{inhomogeneous equation} yields the  homogeneous equation $T^p - \beta(X^p-\beta Y^p)^p=0$,
where the degrees of the generators are $\deg(X)=\deg(Y)=1$ and $\deg(T)=p$.
Inside the weighted projective space $\PP=\PP(1,1,p)$, it defines a closed subscheme $D'$.
This is a Cartier divisor, because 
the coherent sheaf $\O_\PP(p^2)$ is invertible. It lies in the smooth locus of $\PP(1,1,p)$,
because $(0:0:1) $ is the only singular point of the weighted projective space.
By construction, we have an identification $D'\cap D_+(Y)=\Spec(R')$.
In particular, the scheme $D'$ is regular outside the singular point $b'$.
The situation is similar on $D'\cap D_+(X)$, and we infer $\Sing(D')=\{b'\}$.
Summing up, the morphism $D\ra D'$ yields an isomorphism between the regular loci.

Now consider the case $p=2$. Here we claim that the inclusion $R'\subset R$ is actually an equality.
Indeed, we have 
$$
u^2=x^2-\beta,\quad \beta^{1/2}u^2= t\quadand \beta^{1/2}u^3= xt + \beta u^2,
$$
which implies $\beta^{1/2}u^n\in R'$ for all $n\geq 2$, and thus $R'=R$.
Thus our birational morphism $D\ra D'$ is an isomorphism.
\qed

\medskip
Conversely, given some $\beta\in F^\times$ that is not a $p$-power, 
we can consider
the curve $D'$ defined by the homogeneous equation $T^p-\beta(X^p-\beta Y^p)^p=0$
inside the weighted projective space $\PP=\PP(1,1,p)$. 
We already observed that the homogeneous polynomial is irreducible, such that the scheme $D'$ is integral.

\begin{proposition}
\mylabel{peculiar invariants}
The curve $D'$ has  cohomological invariants
$$
h^0(\O_{D'})=1\quadand h^1(\O_{D'})=(p^3-p^2-2p+2)/2.
$$
In characteristic $p=2$, this means $h^0(\O_{D'})=h^1(\O_{D'})=1$.
\end{proposition}

\proof
Consider the long exact cohomology sequence coming from 
the short exact sequence $0\ra\O_\PP(-p^2)\ra\O_\PP\ra\O_{D'}\ra 0$. It starts with 
$$
H^0(\PP,\O_\PP(-p^2))	\lra H^0(\PP,\O_\PP)	\lra H^0({D'},\O_{D'})	\lra H^1(\PP,\O_\PP(-p^2)).
$$
The outer terms vanish, and we get $h^0(\O_{D'})=1$.
The long exact sequence continues with 
$$
H^1(\PP,\O_\PP)\lra H^1({D'},\O_{D'})\lra H^2(\PP,\O_\PP(-p^2))\lra H^2(\PP,\O_\PP).
$$
Again, the outer terms vanish, according to \cite{Dolgachev 1981}, Theorem in Section 1.4, 
and Serre duality gives $h^2(\O_\PP(-p^2))=h^0(\O_\PP(n))$ with $n=p^2-p-2=(p+1)(p-2)$.
The vector space $H^0(\PP,\O_\PP(n))$ equals the homogeneous component  of degree $n$ in the graded ring $S=F[X,Y,T]$.
A basis for this vector space is given by the monomials $X^iY^jT^k$ with exponents satisfying 
$$
i+j+pk=(p+1)(p-2).
$$
This equation has  solutions only for  $0\leq k\leq p-2$,  and for   fixed such $k$
there are $s_k=(p+1)(p-2)-pk+1$ solutions. It follows that $h^1(\O_{D'})$ equals
$$
\sum_{k=0}^{p-2}s_k = (p-1)(p+1)(p-2)  - p\frac{(p-1)(p-2)}{2} + (p-1) = (p^3-p^2-2p+2)/2.
$$
Setting $p=2$   gives the value $h^1(\O_{D'})=1$.
\qed

\medskip
The Jacobian Criterion immediately reveals that $D'$ contains no smooth point.
The subscheme $D'\subset\PP(1,1,p)$ lies in $D_+(X)\cup D_+(Y)$, thus is contained in the smooth locus of the
weighted projective space.
The intersections
$D'\cap V_+(Y)$ and $D'\cap V_+(X)$ are   the two non-rational points $(1:0:\beta^{1/p})$ and $(0:1:-\beta^{1+1/p})$,
respectively.
Both residue fields are isomorphic to $\tF$, and the subscheme is reduced, hence the corresponding local rings are regular.
To understand the singular locus of $D'$, it suffices to look at the affine chart $D_+(Y)$, which is the 
spectrum of $k[t,x]$, with $t=T/Y^p$ and $x=X/Y$, and our homogeneous equation becomes
$$
t^p - \beta(x^p-\beta)^p=0.
$$
Let $u$ be new indeterminate, and consider the polynomial ring $\tF[u]$ 
over the field extension $\tF=F(\beta^{1/p})$. Then \eqref{normalization substitution} defines an injective homomorphism
of $F$-algebras
\begin{equation}
\label{normalization}
R'=F[x,t]/(t^p - \beta(x^p-\beta)^p) \lra \tF[u].
\end{equation}
Localizing with denominator $s=x^p-\beta$, one sees that 
$$
t/s\longmapsto \beta^{1/p}\quadand x- t/s\longmapsto u,
$$
such that $(R')_s=\tF[u]_s$. It follows that the normalization of $D'\subset \PP(1,1,p)$ 
is given by \eqref{normalization}. In particular, $b'=(\beta^{1/p}:1:0)$ is the only singular point,
and its preimage $a\in \tD$ is given by the origin $a=(0:1)\in\PP^1_\tF$.
We now consider the normalization $\tD\ra D'$ and the resulting cartesian square  
$$
\begin{CD}
A'	@>>>	\tD\\
@VVV		@VVV\\
B'	@>>>	D',
\end{CD}
$$
where $B'\subset D'$ is defined by the conductor ideal, and $A'\subset\tD$ is the preimage.

\begin{proposition}
\mylabel{peculiar ramification}
The ramification locus $A'\subset\tD$ is given by $\O_{A'}=\tF[u]/(u^{p(p-1)})$,
the $F$-subalgebra $\O_{B'}\subset\O_{A'}$ is generated by $\beta^{1/p}+u$ and $\beta^{1/p}u^p$,
and   the morphisms $\nu':\tD\ra D'$ factors uniquely over $\nu:\tilde{D}\ra D$.
\end{proposition}

\proof
The conductor square gives a long exact cohomology sequence
$$
0\ra H^0(\O_{D'})\ra H^0(\O_\tD)\oplus H^0(\O_{B'}) \ra H^0(\O_{A'}) \ra H^1(\O_{D'})\ra 0.
$$
We have $h^0(\O_{\tD})=p$ and $h^0(\O_{A'})=pl$ for some integer $l\geq 1$,
and our task is to verify $l=p^2-p$.
Being a complete intersection, the curve $D'$ is Gorenstein, thus $h^0(\O_{B'})=pl/2$.
Taking alternating vector space dimension, the above sequence yields
$$
1- (p + pl/2) + pl - (p^3-p^2-2p+2)/2=0,
$$
and the first assertion follows. Since $\O_{D'}\ra\O_{B'}$ is surjective and $\O_{B'}\ra \O_{A'}$
is injective, the second assertion follows from the description \eqref{normalization substitution}  of the morphism
\eqref{normalization}.
For the last assertion, we use the universal properties of cocartesian squares of schemes
and the corresponding cartesian squares of rings:
The image of the composite map  $\O_{D'}\subset\O_\tD\ra\O_{A'}\ra\O_A$ is generated by  $\beta^{1/p}+u$,
and thus factors over the subring $\O_B\subset\O_A$.
\qed

\medskip
We now establish a structure result on certain algebraic curves of arithmetic genus one.
Recall that for ground fields with $\pdeg(F)=1$, there is precisely one purely inseparable
field extension $\tF$ with $[\tF:F]=p$, namely $\tF=F^{1/p}$.

\begin{theorem}
\mylabel{peculiar classification}
Let $D$ be a proper integral curve with cohomological invariants
$h^0(\O_D)=h^1(\O_D)=1$. Assume $D$ in not regular, contains no rational points, and the local rings $\O_{D,a}$
are geometrically unibranch.
Assume further that $\Pic^p_{D/F}$ contains a rational point, and that the ground field has
$\pdeg(F)=1$.
Then  the curve $D$ is given by a cocartesian square as in \eqref{conductor square curve}, for some
non-zero scalar $\beta\in\tF=F^{1/p}$.
For $p=2$, this $D$ is isomorphic to the curve inside $\PP(1,1,2)$ given
by the homogeneous equation $t^2-\beta(x^2-\beta y^2)^2=0$.
\end{theorem}

\proof
Let
$\nu:\tilde{D}\ra D$ be the normalization map, with  conductor square
\begin{equation}
\label{conductor square curve 2}
\begin{CD}
A	@>>>	\tilde{D}\\
@VVV		@VV\nu V\\
B	@>>>	D.
\end{CD}
\end{equation}
Then $\tilde{D}$ is another proper integral curve, and we write
$\tilde{F}=H^0(\tilde{D},\O_{\tilde{D}})$ for the field of global sections.
This gives a finite field extension $F\subset\tilde{F}$.
Let
\begin{equation}
\label{cohomology sequence1}
0\ra H^0(\O_D)\ra H^0(\O_{\tilde{D}})\oplus H^0(\O_B)\ra H^0(\O_A)\ra H^1(\O_D) \ra H^1(\O_{\tilde{D}})\ra 0.
\end{equation}
be the exact sequence resulting from the conductor square.
Our proof proceeds in a sequence of little steps:

\medskip
{\bf Step 1:} 
\emph{We have $h^0(\O_\tD)=p$, $h^1(\O_{\tilde{D}})=0$, $h^0(\O_B)=p$ and  $h^0(\O_A)=2p$.}
To see this, set $d=h^0(\O_{\tD})=[\tilde{F}:F]$ for the degree of the field extension $F\subset\tF$,
and write $h^0(\O_A)=dl$ for some integer $l\geq 1$.
Since the curve  $D$ is Gorenstein, we then have $h^0(\O_B)=dl/2$. 
Clearly  $h^1(\O_\tD)\leq 1$. If equality holds,
then $d=1$.
Using that the alternating sum of vector spaces dimensions in \eqref{cohomology sequence1} vanishes, we
get   $l=0$, contradiction. Therefore $h^1(\O_\tD)=0$.
In turn, we have  $(d+dl/2)-dl =0$, hence $l=2$.
If $d=1$ we get $h^0(\O_B)=1$, which implies
that $D$ contains a rational point, contradiction. This shows $d>1$.
In order to check $d=p$, we may base-change to the separable closure of the ground field $F$,
and assume that $F$ is separably closed.
Since $h^0(\O_D)=1$ we have an exact sequence
$$
\Pic(D)\lra \Pic_{D/F}(F)\lra\Br(F).
$$
Since the Brauer group vanishes, 
the rational point on $\Pic^p_{D/F}$ comes from an  invertible sheaf on $D$ of degree $p$,
and its preimage $\shL$ on $\tD$ likewise has $\chi(\shL)-\chi(\O_\tD)=\deg(\shL)=p$.
This degree is a multiple of $d>1$, and we conclude $d=p$.

\medskip
{\bf Step 2:} 
\emph{The field extension $F\subset\tilde{F}$ is purely inseparable.}
Suppose this does not hold. Choose some $\alpha\in \tilde{F}$ whose minimal polynomial $f\in F[T]$
has at least two roots $\omega_1\neq \omega_2$ in  some algebraic closure $\Omega$.
This gives a subfield $F[T]/(f)\subset\tilde{F}$ such that the spectrum  of $\tilde{F}\otimes_F\Omega$ 
is disconnected. It follows that the curve $\tilde{D}$ is geometrically reducible, whence
the same holds for $D$. On the other hand,   $D$ is geometrically connected because $h^0(\O_D)=1$.
In turn, there must be a closed point on $D\otimes_F\Omega$ whose local ring has reducible spectrum.
Its image $a\in D$ is a closed point where the local ring $\O_{D,a}$ is not geometrically unibranch,
contradiction.

\medskip
{\bf Step 3:}
\emph{The schemes $A$ contains only one point.}
Suppose this is not the case. Then   the condition $h^0(\O_A)=2p$ implies that 
$\O_A=\tilde{F}\times\tilde{F}$, hence the subring $\O_B\subset\O_A$
is reduced and we have $A=\{a_1,a_2\}$.
The scheme $B$ contains at least one and at most two points.
The case  $B=\{b\}$ is impossible, because then the local ring $\O_{D,b}$ is not unibranch.
Thus  $B=\{b_1,b_2\}$, such that $\O_B=F_1\times F_2$ for some intermediate  fields 
$F\subset F_i\subset\tilde{F}$.
Since the local rings $\O_{D,b_i}$ are  Gorenstein, we have $[\tilde{F}:F_i]=[F_i:F]=p/2$, which forces $p=2$
and $[F_i:F]=1$. The latter implies that $D$ contains rational points, contradiction.

\medskip
{\bf Step 4:}
\emph{We have $\O_A=\tilde{F}[\epsilon]$ for some indeterminate $\epsilon$ subject to $\epsilon^2=0$.}
The task is to show that $A$ is non-reduced. 
Seeking a contradiction, we assume that $A$ is integral. 
Since $\pdeg(F)=1$, we must have $\O_A=F^{1/p^2}$, and there  is but one intermediate field
between $F\subset\O_{A}$, namely $F^{1/p}$.
The difference map $H^0(\O_{\tD})\oplus H^0(\O_B)\ra H^0(\O_A)$ is not surjective,
by the exact sequence \eqref{cohomology sequence1}, hence its image is the intermediate field,
and the cokernel $H^1(\O_D)$ has dimension $p$, contradiction.

\medskip
{\bf Step 5:}
\emph{The composite map $\O_B\subset\O_A\ra \O_A/\maxid_A$ is bijective.}
Consider the inclusions of fields
$F\subset \O_B/\maxid_B \subset \O_A/\maxid_A =\tilde{F}$.
The first inclusion is strict, because $D$ contains no rational points.
The composite extension has degree $p$, hence $\O_B/\maxid_B=\tilde{F}$.
Since $h^0(\O_B)=p$, we must have   $\maxid_B=0$.

\medskip
{\bf Step 6:} 
\emph{The normal curve $\tD$ is isomorphic  to $\PP^1_\tF$.}
We saw in Step 4 that there is an  invertible sheaf $\shL$ on $\tD$ of degree $\deg(\shL)=p$.
Let us  temporarily regard $\tilde{F}$ as the ground field for $\tD$, such that $\deg(\shL)=p/[\tF:F]=1$
and $\deg(\omega_\tD)=-2\chi(\O_\tD)=-2$. It follows that $h^1(\shL)=h^0(\shL^\vee\otimes\omega_\tD)=0$,
and Riemann--Roch gives $h^0(\shL)=\chi(\shL)=\deg(\shL)+\chi(\O_\tD)=2$.
Each non-zero section $s\in H^0(\tD,\shL)$ vanishes at a unique $\tF$-rational point,
and two linearly independent sections vanish at different points. It follows that
$\shL$ is globally generated, thus defines a morphism $\tD\ra\PP^1_\tF$ of degree one,
which thus must be an isomorphism.

\medskip
The proof concludes as follows: Write $\tF=F(\beta^{1/p})$ for some $\beta\in F^\times$ that is not a $p$-power.
The conductor scheme $B$ is canonically isomorphic to $\Spec(\tF)$, via the
composite map in Step 6, but the inclusion $\O_B\subset\O_A$ is given by 
$$
\beta^{1/p}\mapsto \beta^{1/p}+\gamma\epsilon
$$
for some $\epsilon\in F^\times$. 
Replacing the generator $\beta^{1/p}\in\O_B$ by $\gamma^{-1}\beta^{1/p}$, we may assume $\gamma=1$.
Thus the scheme $D$ obtained from $\tD=\PP^1_\tF$ by the denormalization in \eqref{conductor square curve}.
According to Proposition \ref{birational model}, there exists a finite birational universal homeomorphism
$D\to D'$, where $D'\subset \PP(1,1,p)$ is given by the equation $t^p-\beta(x^p-\beta y^p)^p=0$. Moreover, if $p=2$,
the finite birational universal homeomorphism is actually an isomorphism. This means that $D$ is isomorphic to the
curve in weighted projective space given by $t^2-\beta(x^2-\beta y^2)^2=0$.
\qed

\section{Non-existence without irregularity in characteristic two}
\mylabel{Non-existence}

Let $F$ be a ground field of characteristic $p=2$. Throughout, we use the
notation from Section \ref{Regular del Pezzo}:
Let $V$ be a regular del Pezzo  surface that is geometrically integral
but geometrically non-normal. 
Let $N\subset \Sing(V/F)$ be the divisorial part of the locus of non-smoothness,
and $D=N_\red$ the underlying reduced scheme. 
We now settle the sub-case from Proposition \ref{ramification non-reduced p=2} 
occurring in Proposition \ref{pdeg ramification non-reduced}:

\begin{theorem}
\mylabel{non-existence without irregularity}
Assume that   $h^0(\O_D)=h^1(\O_D)=1$, and 
that over the algebraic closure $\Omega=\bar{F}$,
the normalization of $V_\Omega$ is
the contracted Hirzebruch surface $\PP(1,1,2)$, with non-reduced
ramification locus $R=2H$. Then the ground field $F$ must have
$\pdeg(F)\geq 2$.
\end{theorem}

\proof
By assumption, we are   in case (iii) of Reid's Classification.
After making a finite separable extension of $F$, we may assume
that there is a finite purely inseparable extension $F\subset K$ so that
the normalization $\nu:X\ra Y=V_K$ is given by $X=\PP(1,1,2)$,
with non-reduced ramification divisor $R=2H$. The latter has selfintersection number $R^2=2$.
According to Proposition \ref{ramification non-reduced p=2}, we have $\nu^{-1}(D_K)=R$,
thus $D^2=2$, in particular there is an invertible sheaf of degree two on $D$.
Moreover, the curve $D$ is geometrically irreducible.
Being the divisorial part of non-smoothness on some regular scheme,
the integral curve $D$ is geometrically non-reduced  and  contains no rational points,
by Corollaries \ref{prime divisors geom non-reduced} and \ref{closed points geom non-reduced}.

Seeking a contradiction, we now assume that $\pdeg(F)\leq 1$.
In light of \cite{Schroeer 2010}, Theorem 2.3 the curve $D$ 
is not regular.
Thus Theorem \ref{peculiar classification} applies, hence $D$ is given by the homogeneous equation $t^2-\beta(x^2-\beta y^2)^2=0$ inside the 2-dimensional weighted projective space $\PP(1,1,2)$ for some non-square $\beta\in F^\times$.
Starting from this information, we shall show that the regular del Pezzo surface $V$
is a Cartier divisor inside the 3-dimensional weighted projective space $\PP(1,1,1,2)$,
and that its defining equation must produce a singularity $v\in V$.

To achieve this goal, we first establish that $H^1(V,\O_V)=0$. The conductor square
\eqref{conductor_square} gives an exact sequence
$$
H^0(\O_X)\oplus H^0(\O_C)\ra H^0(\O_R)\ra H^1(\O_Y)\ra H^1(\O_X)\oplus H^1(\O_C).
$$
The map on the left is surjective, because the ramification curve $R=2H$ is   the split ribbon $\PP^1\oplus\O_{\PP^1}(-1)$,
hence   $H^0(R,\O_R)=F$.
The contracted Hirzebruch surface $X=\PP(1,1,2)$ has $H^1(X,\O_X)=0$,
and it suffices to check $H^1(C,\O_C)=0$.
According to Proposition \ref{reduced C}, the curve $C$ is geometrically integral.
The two curves $C,D_K\subset Y$ have the same support, and thus $C\subset (D_K)_\red$.
By Proposition \ref{peculiar gorenstein}, we have $(D_K)_\red=\PP^1$ and thus $C=\PP^1$.

Next, consider the invertible sheaf $\shL=\O_V(D)$: the sheaf $\shL$ is ample because
its pullback to $X=\PP_K(1,1,2)$ is ample. The restriction   has
$\deg(\shL_D)=2$, whereas the dualizing sheaf has $\deg(\omega_D)=-2\chi(\O_D)=0$.
In turn, $h^1(\shL_D)=0$ and $h^0(\shL_D)=2$. Choose a vector space basis $s,s'\in H^0(D,\shL_D)$.
Since there is no rational point, the sections vanish at some non-rational points $a,a'\in D$
whose residue fields must be isomorphic to $\tF=F^{1/2}$.
The local rings $\O_{D,a}$ and $\O_{D,a'}$ then must be regular.
Consequently we have $a\neq a'$ because $Fs\neq Fs'$ are different linear systems.
It follows that $\shL_D$ is globally generated.

The short exact sequence $0\ra\O_V\ra\shL\ra\shL_D\ra 0$ yields an exact sequence
$$
0\lra H^0(V,\O_V)\lra H^0(V,\shL)\lra H^0(D,\shL_D)\lra H^1(V,\O_V).
$$
The term on the right vanishes, and we conclude that $\shL$ is globally generated,
with $h^0(\shL)=3$. 
We get a finite morphism $f:V\ra\PP^2$. It must be surjective,
of degree $\deg(f)=\deg(\shL)=2$.
This morphism is actually flat, because $\PP^2$ is regular and $V$ is Cohen--Macaulay.
Setting $\shA=f_*(\O_V)$, we get a short exact sequence
\begin{equation}
\label{infinitesimal covering}
0\lra \O_{\PP^2}\lra \shA\lra \O_{\PP^2}(d)\lra 0
\end{equation}
for some integer $d$. We claim that $d=-2$. To see this, choose a global section $s\in H^0(V,\shL)$
with zero-locus $D\subset V$. Restricting the above   sequence
to the   line $\PP^1=V_+(s) $ 
gives a short exact sequence $0\ra\O_{\PP^1}\ra\nu_*(\O_D)\ra\O_{\PP^1}(d)\ra 0$,
thus $d+1=\chi(\O_{\PP^1}(d))=\chi(\O_D)-\chi(\O_{\PP^1})=-1$, and therefore   $d=-2$.
The short exact sequence \eqref{infinitesimal covering} splits as
an extension of coherent sheaves, because
$$
\Ext^1(\O_{\PP^2}(d), \O_{\PP^1})= H^1(\PP^2,\O_{\PP^2}(-d))=0.
$$
Consequently $\shA=\O_{\PP^2}\oplus\O_{\PP^2}(-2)$, and the algebra structure
is determined  by some linear map  
$$
(\Psi,\Phi):\O_{\PP^2}(-2)\otimes \O_{\PP^2}(-2)\lra \shA=\O_{\PP^2}\oplus\O_{\PP^2}(-2),
$$
where we regard $\Psi,\Phi\in F[x,y,z]$ as  homogeneous polynomials of degree four and two, respectively.
It follows that $V$ is isomorphic to the Cartier divisor in the  weighted projective space $\PP(1,1,1,2)$
defined by the  equation $P(x,y,z,t)=0$, for the homogeneous polynomial
$$
P(x,y,z,t)= t^2 - \Phi(x,y,z) t  - \Psi(x,y,z) 
$$
of degree four, where $\deg(t)=2$.
Without restriction we may assume that   the  curve $D\subset V$ is given by the equation $z=0$.
The locus of non-smoothness $\Sing(V/F)\subset\PP(1,1,1,2)$ is given by the Jacobian ideal
$\ideala=(P,P_t,P_x,P_y,P_z)$ generated by $P$ and its partial derivatives, because
the unique singularity $(0:0:0:1)$ of the weighted projective space lies outside $V$,
and the reflexive sheaf $\O_{\PP(1,1,1,2)}(-4)$ is invertible.
According to Proposition \ref{ramification non-reduced p=2}, the   divisorial part is given by 
$2D\subset\Sing(V/F)$, such that $\ideala\subset(z^2)$. Using   $P_t=\Phi$,
we see that  $z^2|\Phi$, and thus the homogeneous quadratic polynomial is of the form   $\Phi=\alpha z^2$ for some
scalar $\alpha\in F$. According to Theorem \ref{peculiar classification}, we have
$$
\Psi(x,y,z) = \beta(x^2-\beta y^2)^2 + zR(x,y,z)
$$
for some homogeneous polynomial $R(x,y,z)$ of degree three.
Thus $P_z = zR_z + R$. Since the jacobian ideal $\ideala$ is contained in $(z^2)$ and in particular in  $(z)$,
we infer that $z|R$. Thus the polynomial $P$ defining our regular del Pezzo surface
$V\subset\PP(1,1,1,2)$ takes the explicit form
$$
P(x,y,z,t)= t^2 - \alpha z^2 t -\beta(x^2-\beta y^2)^2 + z^2Q(x,y,z)
$$
for some homogeneous quadratic polynomial $Q(x,y,z)$.
This finally leads to the desired contradiction: 
Set $\PP=\PP(1,1,1,2)$, and consider the non-rational closed point $v=(\beta^{1/2}:1:0:0) $.
Then $v\in V\cap D_+(y)$, and the local ring is given by $\O_{V,v}=\O_{\PP,v}/(Py^{-4})$.
But we have 
$$
\frac{t}{y^2},\frac{z}{y},  \left(\frac{x}{y}\right)^2-\beta  \in\maxid_{\PP,v},
$$
which implies $Py^{-4}\in \maxid_{\PP,v}^2$. Consequently, the local ring  $\O_{V,v}$ is not regular, contradiction.
\qed

\medskip
Similar arguments apply to the twin case $h^0(\O_D)=h^1(\O_D)=1$ in characteristic three.
Then the curve $D$ is given by a homogeneous equation $
z^3-\varepsilon(x^3-\varepsilon y^3)=0$, for some non-cube $\varepsilon \in F^\times$.
However, we shall take another approach in   Section \ref{Non-existence p=3}.

\section{Non-existence  with irregularity in characteristic two}
\mylabel{Non-existence irregularity p=2}

Let $F$ be a ground field of characteristic $p=2$. Throughout, we use the
notation from Section \ref{Regular del Pezzo}:
Let $V$ be a regular del Pezzo  surface that is geometrically integral
but geometrically non-normal.
Let $N\subset \Sing(V/F)$ be the divisorial part of the locus of non-smoothness,
and $D=N_\red$ the underlying reduced scheme.

In this section we will rule out the other  sub-case from 
Proposition \ref{ramification non-reduced p=2} occurring in Proposition \ref{pdeg ramification non-reduced}.
Throughout, we assume that  $h^0(\O_D)=h^1(\O_D)=2$, and that
over the algebraic closure $\Omega=\bar{F}$, the normalization of $V_\Omega$
is the contracted Hirzebruch surface $\PP(1,1,2)$.
We then have $\nu^{-1}(D_\Omega)=R=2H$, such that $D^2=2$.
It appears to be the most challenging case, and the main result of this section is:

\begin{theorem}
\mylabel{nonexistence 22}
Regular del Pezzo surfaces $V$ as above do  not exist.
\end{theorem}

Throughout,  we assume that such a surface $V$ exists.
To reach the desired contradiction we may assume that the ground field $F$ is separably closed.
The main idea is to study   2-dimensional linear systems $ H^0(V,\shL)$
coming from   invertible sheaves $\shL=\shN(D)$ with $\shN$ numerically trivial. It turns out that these define genus-one
fibration   on blow-ups, and the Canonical Bundle Formula
will yield the desired  contradiction. We start by computing some relevant cohomology groups:

\begin{proposition}
\mylabel{cohomology C 22}
The conductor curve $C\subset Y$ has $h^0(\O_C)=h^1(\O_C)=1$.
\end{proposition}

\proof
The conductor curve $C$ is geometrically integral by Proposition \ref{reduced C}, thus $h^0(\O_C)=1$. 
It also follows that we have an inclusion $C\subset D_K$ as subschemes inside $Y=V_K$, and get
a short exact sequence
$$
0\lra\shI\lra\O_{D_K}\lra\O_C\lra 0
$$
for some sheaf of ideals $\shI\subset\O_{D_K}$ with $h^0(\shI)=1$. The long exact cohomology sequence
yields   surjections 
\begin{equation}
\label{surjection}
H^0(D_K,\O_{D_K})\lra H^0(C,\O_C)\quadand H^1(D_K,\O_{D_K})\lra H^1(C,\O_C).
\end{equation}
The former is  not injective.
Let $f\in H^0(D_K,\O_{D_K})$ be an element in the kernel.
By assumption,  $H^1(\O_D)$ is a one-dimensional vector space   over the field $H^0(\O_D)$,
so $H^1(D_K,\O_{D_K})$ is a free module of rank one over the ring $H^0(D_K,\O_{D_K})$.
For each cohomology class $\alpha$, the element $f\alpha$ becomes
zero in $H^1(C,\O_C)$. In turn, the surjection $H^1(D_K,\O_{D_K})\ra H^1(C,\O_C)$
is also not injective. This already gives $h^1(\O_C)\leq 1$. 

Seeking a contradiction, we assume that $h^1(\O_C)=0$. 
Since $D^2_K=2$, the invertible sheaf $\shL=\O_Y(D_K)$ has $\deg(\shL_C)=1$,
and it follows that $\shL_C$ is very ample with $h^0(\shL_C)=2$, thus defining an isomorphism $C\ra\PP^1$.
Consequently $h^1(\shI)=2$.
Let $\eta\in D_K$ be the generic point. Since $D_K^2=2$, the length of the local Artin ring
$\O_{D_K,\eta}$ is two, and it follows that $\shI$ annihilates itself.
Hence $\shI$ becomes a torsion-free $\O_C$-module, thus $\shI=\O_{\PP^1}(n)$ for some integer $n$.
Using $h^0(\shI)=1$ we see $n=0$, and this gives $h^1(\shI)=0$, contradiction.
The only remaining possibility is $h^1(\O_C)=1$.
\qed

\medskip
Using Proposition \ref{neron-severi sequence}, we infer the cohomological invariants  for $Y=V\otimes_FK$
and thus also for $V$:

\begin{proposition}
\mylabel{cohomology V 22}
The regular del Pezzo surface $V$ has cohomological invariants $h^0(\O_V)=h^1(\O_V)=1$ and $h^2(\O_V)=0$.
\end{proposition}
 
From this we deduce the structure of the Picard scheme:

\begin{proposition}
\mylabel{picard scheme}
The connected component $\Pic^0_{V/F}$ is a twisted form of $\GG_{a,F}$,
and it contains infinitely many rational points. In turn, there
are infinitely many isomorphism classes of numerically trivial invertible sheaves 
$\shN$ on $V$.
\end{proposition}

\proof
In light of $h^2(\O_V)=0$, the  Picard scheme $P=\Pic^0_{V/F}$ is smooth (\cite{Mumford 1966}, Lecture 27).
It is a  one-dimension scheme, because its tangent space $H^1(V,\O_V)$ is a one-dimensional vector space.
The  base-change to $K$ lies in the kernel for the pull-back map $\Pic_{Y/K}\ra\Pic_{X/K}=\ZZ$.
Using \cite{SGA 6}, Expos\'e XII, Theorem 1.1 we infer that $P$ is quasiaffine.
By the classification of one-dimensional commutative affine group schemes in \cite{Waterhouse; Weisfeiler 1980},
it must be a twisted form of $\GG_a$ or $\GG_m$.  

Now regard $P$ as a smooth connected algebraic curve, write $P\subset\bar{P}$ for the regular compactification,
and choose   some very ample invertible sheaf $\shL$ on $\bar{P}$.
Since the projective curve $\bar{P}$ is geometrically reduced, we may apply
Bertini Theorems (for example \cite{EGA V}, Proposition 4.3 or \cite{Jouanolou 1983}, Chapter I, Theorem 6.3)
and deduce that there are infinitely many reduced Cartier divisors $A\subset\bar{P}$
disjoint from the points at infinity. Since our ground field $F$ is separably closed,
these Cartier divisors are sums of rational points. They correspond to invertible sheaves on $V$,
because $\Br(F)=0$. In other words, there are infinitely many isomorphism classes of
numerically trivial sheaves on $V$.

Seeking a contradiction, we assume that $P$ is a twisted form of $\GG_m$. Choose some prime number $\ell\neq p$.
The kernel $P[\ell]$ is a twisted form of the finite \'etale group scheme $\mu_\ell$.
Since $F$ is separably closed, we conclude that $P[\ell]=\mu_\ell\simeq\ZZ/\ell\ZZ$.
In turn, we find an invertible sheaf $\shL$ of order $\ell$.
Choose a trivialization $\shL^{\otimes -\ell}\ra\O_V$, such that the locally free sheaf
$\shA=\O_V\oplus\shL^{\otimes -1}\oplus\ldots\oplus\shL^{\otimes 1-\ell}$ acquires the structure of
a finite \'etale $\O_V$-algebra. This shows that the algebraic fundamental group $\pi_1(V)$ is
non-trivial. However, the morphism $\nu:X\ra V$ is a universal homeomorphism, so the induced map $\pi_1(X)\ra\pi_1(V)$
is bijective by \cite{SGA 1}, Expos\'e IX, Theorem 4.10. However, $X=\PP(1,1,2)$ is simply-connected, contradiction.
\qed

\medskip
Note that the  structure of twisted forms of the additive group scheme 
$\GG_a$ were determined in \cite{Russell 1970}.
We now consider invertible sheaves $\shL=\shN(D)$ on $V$,
where $\shN$ is numerically trivial. 

\begin{proposition}
\mylabel{linear system}
The  linear system $H^0(V,\shL)$ is two-dimensional.
\end{proposition}

\proof
We compute this with the induced invertible sheaf $\shL_Y$ on the non-normal del Pezzo surface $Y=V\otimes_FK$,
by using the exact sequence
\begin{equation}
\label{linear system 22}
0\lra H^0(Y,\shL_Y)\lra H^0(X,\shL_X)\oplus H^0(C,\shL_C)\lra H^0(R,\shL_R).
\end{equation}
explained in \eqref{section sequence}.
The contracted Hirzebruch surface $X=\PP(1,1,2)$ has $H^1(X,\O_X)=0$, and
it follows that $\shL_X=\O_X(R)=\O_X(2H)$, regardless of the numerically trivial sheaf $\shN$.
Pulling back along the resolution of singularities $S\ra X$ and working on the
Hirzebruch surface $S\ra\PP^1$, one easily sees $h^0(\shL_X)=4$.

Since $\nu^{-1}(D_K)=R=2H$, we have $(D_K\cdot C)=1$, so $\deg(\shL_C)=1$. 
Serre Duality gives $h^1(\shL_C)=0$ and thus $h^0(\shL_C)=1$ by Riemann--Roch. In turn, the base locus for the invertible
sheaf $\shL_C$ comprises a unique rational base point $b\in C$, contained in the regular locus of the curve $C$.

Finally, the ramification divisor is the split ribbon $R=\PP^1\oplus\O_{\PP^1}(-1)$. In turn, 
we get  
$\shL_R=\O_{\PP^1}(1)\oplus\O_{\PP^1}$,
and therefore $h^0(\shL_R)=3$.
The short exact sequence $0\ra\O_X(-R)\ra\O_X\ra \O_R\ra 0$ yields an exact sequence
$$
H^0(X,\O_X)\lra H^0(R,\O_R)\lra H^1(X,\O_X(-R)).
$$
We have  $h^2(\O_X(-R))=h^0(\O_X(-R))=0$ by Serre Duality. Riemann--Roch gives
$$
\chi(\O_X(-R)) = \frac{R^2 + (R\cdot K_X)}{2} + \chi(\O_X) = \frac{R^2 -2R^2}{2} + 1= 0,
$$
and thus $h^1(\O_X(-R))=0$. In turn, the restriction map $H^0(X,\shL_X)\ra H^0(R,\shL_R)$ is surjective.
Using the exact sequence \eqref{linear system 22}, we infer $h^0(\shL)=2$.
\qed

\medskip
Next, we examine the 2-dimensional linear system on $Y=V_K$ induced by $\shL_Y$:

\begin{proposition}
\mylabel{strange divisor}
There is a non-zero global section $s'\in H^0(Y,\shL_Y)$ so that the resulting
effective Cartier divisor $D'\subset Y$ has $\Supp(D')=\Supp(D_K)$. Such a section
is unique up to invertible scalar. The  cohomological invariants for the scheme $D'$   are given by
$$
h^0(\O_{D'})=h^1(\O_{D'})=
\begin{cases}
2 & \text{if $\shL\simeq \omega_V^\vee$;}\\
1 & \text{else.}
\end{cases}
$$
\end{proposition}

\proof
According to Proposition \ref{sections vanishing}, the global sections of $\shL_Y$ vanishing
along $C$ correspond to the global sections of $\shL_X(-R)$.
In our situation,   the normalization $X$  is the contracted Hirzebruch surface $\PP(1,1,2)$
and $\shL_X(-R)=\O_X$. Existence and uniqueness follow.

Since the conductor curve $C$ is reduced, we have an inclusion $C\subset D'$ inside
the scheme $Y$. In turn, there is  a short exact sequence
\begin{equation}
\label{first sequence}
0\lra\shI\lra\O_{D'}\lra \O_C\lra 0
\end{equation}
for some sheaf of ideals $\shI\subset\O_{D'}$.
By definition of the conductor square,   $\nu^{-1}(C)=R$ holds as subschemes on $X$.
By assumption on our regular del Pezzo surface $V$ we also have   $\nu^{-1}(D')=R$. In turn,
the morphism $\nu:R\ra Y$ factors over $D'$, and the image of the homomorphism
$\O_{D'}\ra\O_R$ is the subsheaf $\O_C\subset\O_R$. Furthermore, we have a commutative
diagram
$$
\begin{CD}
0	@>>>	\O_Y(-D')	@>>>	\O_Y	@>>>	\O_{D'}	@>>> 	0\\
@.		@VVV			@VVV		@VVV\\
0	@>>>	\O_X(-R)	@>>>	\O_X	@>>>	\O_R	@>>> 	0.
\end{CD}
$$
The cokernels for the two vertical maps on the right are the same, and the Snake Lemma gives an
identification $\shI=\O_X(-R)/\O_Y(-D')$. In turn, we have a second short exact sequence
\begin{equation}
\label{second sequence}
0\lra \O_Y(-D')\lra \O_X(-R)\lra \shI\lra 0.
\end{equation}
Riemann--Roch gives
$$
\chi(\O_X(-R)) = \frac{(R^2) + (R\cdot K_X)}{2} + \chi(\O_X) = \frac{(R^2)-2(R^2)}{2} + 1 = 0
$$
and
$$
\chi(\O_Y(-D')) = \frac{(R^2) + (R\cdot \nu^{-1}(K_Y))}{2} + \chi(\O_Y)= \frac{(R^2) - (R^2)}{2} + \chi(\O_Y)  = 0.
$$
Consequently $\chi(\shI)=0$ and $\chi(\O_{D'})= \chi(\O_C) + \chi(\shI) = 0$.
Thus it suffices to verify the statement on $h^0(\O_{D'})$.

Obviously $h^0(\O_X(-R))=0$. Using $\omega_X(R)=\O_X(-R)$ also $h^2(\O_X(-R))=0$,
and thus   $h^1(\O_X(-R))=0$.
Likewise, we have $h^0(\O_Y(-D'))=0$. If $\shL=\omega_Y^\vee$, then Serre duality gives 
$h^2(\O_Y(-D'))=1$ and thus $h^1(\O_Y(-D'))=1$. 
On the other hand, if  $\shL\neq\omega_Y^\vee$, the invertible sheaf $\omega_Y(D')$ is numerically trivial
but non-trivial, and thus $h^i(\O_Y(-D'))=0$ for all $i\geq 0$.
From  \eqref{second sequence} we get an exact sequence
$$
0\lra H^0(\shI)\lra H^1(Y,\O_Y(-D'))\lra H^1(X,\O_X(-R)).
$$
We just saw that the term on the right vanishes, and thus $h^0(\shI)=1$
if $\shL=\omega_Y^\vee$, and $h^0(\shI)=0$ if $\shL\neq\omega_Y^\vee$.
To finish the argument, we use the exact sequence
$$
0\lra H^0(\shI)\lra H^0(D',\O_{D'})\lra H^0(C,\O_C).
$$
stemming from \eqref{first sequence}.
We have $h^0(\O_C)=1$, because  $C$ is geometrically reduced.
It follows that the map on the right is surjective, and the assertion follows.
\qed

\medskip
Note that we do not assert that  the above Cartier divisor $D'$ on $Y=V_K$
arises from an effective Cartier divisor on $V$.
In fact, this is impossible if  $\shN\not\simeq\O_V$,
because on the normal scheme $V$   effective Cartier divisors
are determined by their support and multiplicities.  

By assumption on our regular del Pezzo surface $V$, the divisor $D=N_\red$,
where $N\subset\Sing(V/F)$ is the divisorial part of the locus of non-smoothness, has
$h^0(\O_D)= h^1(\O_D)=2$. Let us record the following consequence:

\begin{corollary}
\mylabel{dualizing 22}
The dualizing sheaf for the  del Pezzo surface is $\omega_V=\O_V(-D)$.
\end{corollary}

\begin{proposition}
\mylabel{invariants D'}
Suppose  that  $\shN\not\simeq\O_V$. For each non-zero      $s'\in H^0(V,\shL)$
the resulting effective Cartier divisor $D'\subset V$ has cohomological invariants
$h^0(\O_{D'})=h^1(\O_{D'})=1$.
\end{proposition}

\proof
Riemann--Roch gives
$$
\chi(\O_{D'})=\chi(\O_V) - \chi(\O_V(-D')) = \chi(\O_V)-\chi(\O_V(-D)) = \chi(\O_D)= 0,
$$
so it suffices to check $h^0(\O_{D'})=1$.
The curve $D'\subset V$ is integral, because its class generates $\NS(V)$.
In turn, the preimage $\nu^{-1}(D'_K)\subset X$ is irreducible.
Since this is linearly equivalent to $R=2H$, it must be isomorphic to the split ribbon
$\PP^1\oplus\O_{\PP^1}(-1)$. Using that $\nu^{-1}(D'_K)\ra D'_K$ is schematically dominant,
we infer that $h^0(\O_{D'_K})=1$.
\qed

\medskip
For each vector space basis $s_1,s_2\in H^0(V,\shL)$ the resulting two effective
Cartier divisor $D_1,D_2\subset V$ define the base-locus $\Bs(\shL)=D_1\cap D_2$ as a closed subscheme.

\begin{proposition}
\mylabel{base locus}
The base-locus $\Bs(\shL)$ consists of a unique closed point $b\in V$
contained in $D\subset V$, and its residue field has  $[\kappa(b):F]=2$.
\end{proposition}

\proof
First note that  formation of base-loci commute with extension of ground fields.
The N\'eron--Severi group $\NS(Y)$ is cyclic, with $\shL$ as generator.
It follows that every non-zero global section of $\shL$ defines an integral
Cartier divisor. From this we infer that the base-locus $\Bs(\shL)=D_1\cap D_2$
is finite. It must have degree two, because $D^2=2$.

Using the restriction map $H^0(\shL_Y)\ra H^0(C,\shL_C)$ and the fact
that $\shL_C$ has degree one on the integral curve $C$ of genus one,
we see that $\Bs(\shL)$ consists of a unique closed point $b\in V$ contained in $D$.
Since our del Pezzo surface $V$ is regular and 
$D\subset N$ is the divisorial part of the locus of non-smoothness,
we see that the point $b\in V$ is not rational.
Since the base-locus has length two as closed subscheme, we infer that 
the field extension $F\subset\kappa(b)$ has degree two.
\qed

\medskip
The two-dimensional linear system $H^0(V,\shL)$   defines a rational map  
$$
X\dashrightarrow  \PP^1=\PP H^0(V,\shL)
$$
defined on the complement of the base-point $b\in V$.
Let $\tilde{V}\ra V$ be the blowing-up with center $b\in V$. The exceptional divisor $E\subset\tilde{V}$
is a projective line over the residue field $\kappa(b)$, and we get a fibration
$f:\tilde{V}\ra\PP^1$.
Note that all this depends on the numerically trivial sheaf $\shN$ in $\shL=\shN(D)=\shN\otimes\omega_V^\vee$, but we do 
not indicate this in notation.

\begin{proposition}
\mylabel{fibers integral}
If $\shN\neq \O_V$, then $f:\tilde{V}\ra\PP^1$ is a genus-one fibration with the following property:
for every point $a\in \PP^1$, the schematic fiber $\tilde{V}_a=f^{-1}(a)$ is integral.
\end{proposition}

\proof
We first establish that   $f:\tilde{V}\ra\PP^1$, which is proper, surjective and flat is indeed a genus-one fibration.
According to Proposition \ref{invariants D'}, for each rational point $a\in \PP^1$ the fiber $\tilde{V}_a= f^{-1}(a)$ 
has $h^0(\O_{\tilde{V}_a})=h^1(\O_{\tilde{V}_a})=1$. For the generic fiber this implies $\chi(\O_{\tilde{V}_\eta})=0$.
Semicontinuity gives $h^0(\O_{\tilde{V}_\eta})=1$, and thus also $h^1(\O_{\tilde{V}_\eta})=1$.  
Zariski's Main Theorem ensures that the canonical injection $\O_{\PP^1}\subset f_*(\O_{\tilde{V}})$ is an equality.
Summing up, the morphism $f:\tilde{V}\ra\PP^1$ is a genus-one fibration.

Now fix a closed point $a\in \PP^1$. It remains to verify that the scheme
$\tilde{V}_a=f^{-1}(a)$ is integral.
This is clear if $a\in\PP^1$ is a rational point, because the image
of    the fiber in $V$ generates $\NS(V)$.
For arbitrary $a\in\PP^1$, the fiber $\tilde{V}_a$ is at least irreducible,
because the group $\NS(\tilde{V})$ has rank $\rho=2$.
Write $\tilde{V}_a=m\Theta$ for some integral curve $\Theta\subset\tilde{V}$ and some
multiplicity $m\geq 1$. 
The exceptional divisor $E\subset\tilde{V}$ for the blowing-up $\tilde{V}\ra V$ is a projective line over the
quadratic field extension $F\subset \kappa(b)$ and the projection $E\ra\PP^1$ has degree two, we must have $m\leq 2$.

Seeking a contradiction, we now assume $m=2$. In particular, the closed point $a\in\PP^1$
is not rational. 
The intersection $E\cap \tilde{V}_a$ is a finite scheme, which has degree two over the
residue field $\kappa(a)$.
In turn, $E\cap \Theta$ is isomorphic to $\Spec\kappa(a)$. We infer that the field extension
$\kappa(a)\subset H^0(\O_\Theta)$ is bijective, and that 
the curve $\Theta$ is regular at the intersection point $E\cap \Theta$.
It follows that the composite morphism $\Theta\subset\tilde{V}\ra V$ is 
a closed embedding. Its image $D''\subset V$ passes through the center
$b\in V$, which is also the base-point for the invertible sheaf $\shL=\shN(D)$,
and has $[\kappa(b):F]=2$. Choose a non-zero global section of $H^0(V,\shL)$.
Recall that $s_1,s_2\in H^0(V,\shL)$ is a vector space basis, with corresponding 
effective Cartier divisors $D_1,D_2\subset V$. Since $(D_1\cdot D_2)=2$, the intersection $D_1\cap D_2$
is transversal, and in particular the curves $D_1,D_2$ are regular at the intersection point $b\in D_1\cap D_2$.
Without restriction, we may assume that the intersection $D_1\cap D''$ is transversal as well,
and thus $(D_1\cdot D'')=2$. It follows that $\shL=\shN(D)$ is numerically equivalent to $\O_V(D'')$.
Furthermore, the pullback to the conductor curve $C\subset Y=V_K$ become isomorphic.
Since both pullback maps  $\Pic(V)\ra\Pic(Y)$ and   restriction map $\Pic(V)\ra\Pic(C)$ are injective,
we conclude that $\O_V(D'')\simeq\shL$. It follows that $\Theta$ is a fiber 
over some rational point in $\PP^1$, contradiction.
\qed

\medskip
\emph{Proof for Theorem \ref{nonexistence 22}.}
According to Corollary \ref{dualizing 22}, we have $\omega_V=\O_V(-D)$.
Choose some numerically trivial invertible sheaf $\shN\not\simeq\O_V$,
and consider the genus-one fibration $f:\tilde{V}\ra\PP^1$ constructed above from  the invertible sheaf $\shL=\shN(D)=\O_V(D')$.
According to Proposition \ref{fibers integral}, all fibers over closed points of $\PP^1$ are integral. 
So we have $K_{\tilde V}=\sum m_i F_i$, where $m_i$ are integers and $F_i$ are fibers. 
Since the base of the fibration is $\PP^1$, we deduce that $K_{\tilde V}=mF$, for the fiber $F$
over  the rational point $\infty\in\PP^1$ 
and some integer $m< 0$. With $E^2=-2$ and $(F \cdot E)=2$, we deduce that $m=-1$.
Using that $\omega_V=f_*(\omega_{\tilde{V}})$, we infer that $\omega_V=\O_V(-D')$.
On the other hand, we have  $\omega_V=\O_V(-D)$  by Corollary \ref{dualizing 22}.
Thus the two divisors $D,D'\subset V$ are linearly equivalent, hence $\shN\simeq\O_V$, contradiction.
\qed

\section{Peculiar curves of higher genus}
\mylabel{Peculiar higher genus}

In this section we shall study the geometry  of certain peculiar algebraic curves  of higher genus that show up  in our study of
del Pezzo surfaces in characteristic three that are regular and geometrically non-normal.
The geometry of such curves, however, merits a study in all characteristics.
Our methods break down in characteristic two, so we exclude this from consideration.

Let $F$ be an imperfect ground field of characteristic $p\geq 3$.
For the sake of exposition, we assume that $F$ is separably closed and has $\pdeg(F)=1$.
In turn, every finite field extension is of the form $F^{1/p^n}$, which has degree $p^n$.
We write $ \tF=F^{1/p}$ for the unique field extension of degree $p$.
Throughout, $D$ denotes a proper integral curve with $h^0(\O_D)=1$ satisfying the following conditions:

\begin{enumerate}
\item The genus of the curve is $h^1(\O_D)=p+1$.
\item The local rings $\O_{D,a}$, $a\in D$ are Gorenstein and geometrically unibranch.
\item There is no rational point on $D$.
\item The normalization $\tD$ has genus $h^1(\O_\tD)=0$. 
\end{enumerate}

\noindent
In turn, we get the numerical invariants
$$
\chi(\O_D)=-p,\quad \deg(\omega_D)=2p,\quad h^0(\omega_D)=p+1\quadand h^1(\omega_D)=1.
$$
Note that for each closed point $x\in D$, the residue field $\kappa(x)$ is purely inseparable
finite extension of $F$, so its degree is a $p$-power.
As usual, we form  the conductor square
\begin{equation}
\label{conductor higher genus}
\begin{CD}
A	@>>>	\tD\\
@VVV		@VV\nu V\\
B	@>>>	D
\end{CD}
\end{equation}
for the normalization morphism $\nu:\tD\ra D$.

\begin{proposition}
\mylabel{normalization line}
The normalization $\tD$ is isomorphic to the projective line $\PP^1_\tF$,
and the conductor loci satisfy $h^0(\O_B)=2p$ and $h^0(\O_A)=4p$. More precisely, we have
$$
\O_A\simeq\tF[\epsilon]\times\tF[\epsilon]\;\;\text{and}\;\; \O_B\simeq\tF\times\tF,
\qquad\text{or}\qquad 
\O_A\simeq\tF[\eta]\;\;\text{and}\;\;  \O_B\simeq\tF[\epsilon].
$$
Here $\epsilon$ and $\eta$ denote indeterminates subject to the relation $\epsilon^2=0$ and $\eta^4=0$.
\end{proposition}
 
\proof
Seeking a contradiction, suppose that $h^0(\O_\tD)=1$. According to \cite{Schroeer 2010}, Theorem 2.3
the scheme $\tD$ is geometrically reduced, and thus the same holds for $D$.
According to the well-known Bertini Theorems (for example \cite{EGA V}, Proposition 4.3 or \cite{Jouanolou 1983}, Chapter I, Theorem 6.3), 
we find an effective Cartier divisor $C\subset D$ that is also geometrically reduced.
Since $F$ is separably closed,  such a subscheme  consists of rational points, contradiction.
So we have $h^0(\O_\tD)=p^\nu$ for some exponent $\nu\geq 1$.
The invertible sheaf $\omega_D|\tD$ has degree $2p$, which ensures $\nu=1$, because $p$ is odd.
Since $\pdeg(F)=1$, the two field extensions $\tF$ and $H^0(\tD,\O_\tD)$ coincide.

Now regard $\tF$ as the ground field for the curve $\tD$, which is geometrically reduced.
Applying Bertini Theorems as above, we find a rational point $x\in \tD$, and set  $\shL=\O_\tD(x)$.
The long exact sequence coming from   $0\ra\O_\tD\ra \shL\ra \kappa(x)\ra 0$
shows that $\shL$ is globally generated  with $\dim_{\tF}H^0(\tD,\shL)=2$, and it follows that the ensuing  
$\tD\ra\PP^1_\tF$ is an isomorphism.

From the conductor square \eqref{conductor higher genus}, we get an exact sequence
$$
0\ra H^0(\O_D)\ra H^0(\O_\tD)\oplus H^0(\O_B)\ra H^0(\O_A)\ra H^1(\O_D)\ra 0.
$$
Write  $h^0(\O_A)=pa$ for some integer $a\geq 1$. Since the curve $D$ is Gorenstein, we must have $h^0(\O_B)=pa/2$.
In turn, we get $ 1 - (p+pa/2) + pa - (p+1)=0$. Thus $a=4$, and the values for $h^0(\O_A)$ and $h^0(\O_B)$ follow.

Now decompose $A=\sum n_i a_i$ into prime divisors on the regular curve $\tD$, and write $[\kappa(a_i):\tF]=p^{\nu_i}$.
Then $4=\sum n_ip^{\nu_i}$. Each summand is even, because the local rings on $D$ are Gorenstein and geometrically unibranch, 
and $p$ is odd. 
One easily sees that the only solutions are $4=2+2$ or $4=4$, and the structure for $\O_A$ follows.
The Artin ring $\O_B$ of length $h^0(\O_B)=2p$ has an induced decomposition, with corresponding lengths
$2p=p+p$ or $2p=2p$, since the local rings on $B$ are geometrically unibranch. 
Since $D$ contains no rational point, all residue fields of $\O_B$ are copies of $\tF$.
The structure of $\O_B$ follows.
\qed
 
\medskip
Choose an identification $\tD=\PP^1_\tF$.
We now construct some partial denormalization $\hat{D}$, with morphisms $\PP^1_\tF=\tD\ra\hat{D}\ra D$.
First suppose that $\O_{A}=\tF[\epsilon]\times\tF[\epsilon]$.
Then define $\hD$ as the denormalization of $\tD=\PP^1_\tF$
with respect to $\O_B/(0\times\tF)$ inside the ring $\O_A/(0\times\tF[\epsilon])$.
Then $\hD$ is either a peculiar genus-one curve studied in Section~\ref{Peculiar curves}, or the rational cuspidal curve 
$\Spec\tF[T^2,T^3]\cup \Spec\tF[1/T]$.

Now suppose that  $\O_A=\tF[\eta]$, with $\eta^4=0$, and write $\tF=F(\beta^{1/p})$.
Since $h^0(\O_D)=1$, the subalgebra  $\O_B=\tF[\epsilon]\subset\tF[\eta]=\O_A$ is given by
$$
\beta^{1/p}\mapsto \beta^{1/p} +  \lambda_1\eta  + \lambda_2\eta^2 + \lambda_3\eta^3 
\quadand
\epsilon = \mu_2\eta^2+\mu_3\eta^3,
$$
where neither $(\lambda_1,\lambda_2,\lambda_3)$ nor $(\mu_2,\mu_3)$ are   zero tuples.
Consider $\O_A/\maxid_A^2=\tF[\eta]/(\eta^2)$ and the resulting subring $\O_B/(\O_B\cap\maxid_A^2)$.
Both residue class rings have residue fields isomorphic to $\tF$, with 
  $h^0(\O_A/\maxid_A^2)=2p$ and $h^0(\O_B/(\O_B\cap\maxid_A^2))=p$.
More precisely, $\O_B/(\O_B\cap\maxid_A^2)$ coincides with image 
of $\tF=H^0(\tD,\O_\tD)$ inside $\O_A/\maxid_A^2$ if and only if $\lambda_1=0$.
Now define $\hD$ as the denormalization of $\tD=\PP^1_\tF$
with respect to $\O_B/(\O_B\cap\maxid_A^2)$ inside $\O_A/\maxid_A^2$.
Again $\hD$ is either the rational cuspidal curve over $\tF$ or
a peculiar genus-one curve studied in   Section \ref{Peculiar curves}.
The latter holds if and only if $\lambda_1\neq 0$. 

The following statement, valid in all the above cases, is immediate if $\hD$ is a rational cuspidal curve over $\tF$,
and follows directly from Proposition \ref{peculiar gorenstein} if
$\hD$ is a peculiar genus-one curve:

\begin{proposition}
\mylabel{bound genus C}
The reduction $C=(D\otimes_F\tF)_\red$  of the base-change   has 
$$
\dim_\tF H^0(C,\O_C)=1\quadand  \dim_\tF H^1(C,\O_C)\leq 1.
$$
\end{proposition}
 
In combination with Maddock's bound \eqref{maddock bound},
this will play a crucial role in our analysis of   del Pezzo surfaces $V$ in characteristic three
that are regular but geometrically non-normal. To proceed, we study the linear
system attached to the dualizing sheaf:

\begin{proposition}
\mylabel{dualizing globally generated}
The invertible sheaf $\omega_D$ is globally generated.
\end{proposition}

\proof
Choose a non-zero global section $s\in H^0(D,\omega_D)$.
The zero-locus is an effective Cartier divisor  $C\subset D$ of   $\length(\O_C)=2p$.
As a set, it consists of either two points $c_1,c_2\in D$  or of a single point $c\in D$, because $F$
is separably closed and $C(F)=\varnothing$. In all cases, the respective residue fields $\tF_i$ are isomorphic to $\tF$,
and $\O_C=\tF\times\tF$ or $\O_C=\tF[\epsilon]$,
where  $\epsilon$ denotes an indeterminate subject to the relation $\epsilon^2=0$.
The case that $\O_C$ is a field does not occur, by    our assumption $p\neq 2$.

First consider the case that $C$   is reduced, such that $C=\{c_1,c_2\}$.
The short exact sequence $0\ra\O_D\ra\omega_D\ra \O_C\ra 0$ induces an exact sequence
$$
0\lra H^0(D,\O_D)\lra H^0(D,\omega_D)\lra H^0(C,\O_C)=\tF\times\tF.
$$
The linear  map on the right has rank $p$, so we may assume without restriction that $c_2\in D$ is not
a base point. Seeking a contradiction, we assume that $c_1\in D$ is a base point.
Set  $\shL=\O_D(c_2)=\omega_D(-c_1)$. Then the inclusion $H^0(D,\shL)\subset H^0(D,\omega_D)$ is
an equality, in particular $h^0(\shL)= p+1$, and furthermore $\shL$ is globally generated.
The short exact sequence $0\ra \shL(-c_1)\ra\shL \ra \kappa(c_1)\ra 0$ yields an exact sequence
$$
0\lra H^0(D,\shL(-c_1))\lra H^0(D,\shL)\lra \tF.
$$
The map on the right is not injective, and we conclude that the numerically trivial sheaf $\shL(-c_1)=\O_D(c_2-c_1)$
admits a non-zero section. In turn, $\O_D(c_1)\simeq\O_D(c_2)$ is globally generated, contradiction.

Now suppose that $C$ is non-reduced, such that $\O_C=\tF[\epsilon]$. Seeking a contradiction, we assume that $c\in D$
is a base point. Then the image of the restriction map $H^0(D,\omega_D)\ra H^0(C,\O_C)$ factors over
the maximal ideal $\tF\epsilon$ of the local Artin ring $\O_C$, and we get a short exact sequence
$$
0\lra H^0(D,\O_D)\lra H^0(D,\omega_D)\lra \tF\epsilon\lra 0.
$$
Choose two sections $s',s''$ of $\omega_D$ vanishing at $c\in D$, whose images in $\tF\epsilon$ are linearly independent over $F$.
Then the resulting closed subschemes $D',D''\subset D$ do not coincide, but we have $D'\cap C=D''\cap C$.
It follows that the local ring $\O_{D,c}$ is singular, with embedding dimension two.

We now distinguish two cases, according to the singularities on $D$.
First suppose that the ramification divisor $A\subset\tD$ is of the form $\O_A=\tF[\epsilon]\times\tF[\epsilon]$,
such that $\O_B$ is isomorphic to $\tF\times\tF$. We may assume that the first projection corresponds to the
inclusion $c\in B$. We now  regard the global sections $s\in H^0(D,\omega_D)$
as pairs $(s|\tD,s|B)$ whose entries coincide in $\O_A$. Since $c\in B$ is a base point, this means that
$s|B\in 0\times\tF$ and $s|\tD$ is a global section of $\O_\tD\subset \omega_D|\tD=\O_\tD(2)$.
In other words, we have an identification of $H^0(D,\omega_D)$ with the intersection
of two copies of $\tF$ as subalgebras inside $0\times\tF[\epsilon]$.
Such an intersection is a subalgebra inside $\tF$, which has dimension $d=p$ or $d=1$,
contradicting $h^0(\omega_D)=p+1$.
 
Now suppose that $\O_B=\tF[\epsilon]$ inside $\O_A=\tF[\eta]$, with $\eta^4=0$.
The inclusion is given by $\epsilon=\alpha\eta^2+\beta\eta^3$ for some scalars $\alpha,\beta\in F$, not
both zero. 
The residue class ring $\O_{A}/\epsilon\O_A$ is thus $\tF[\eta]/(\eta^n)$ with $n=2$ or $n=3$.
In both cases,  we may regard the global section $s\in H^0(D,\omega_D)$
as pairs $(s|\tD,s|B)$ where the first entry is a global section of $\O_\tD\subset\omega_D|\tD$ and 
the second entry lies is $\tF\epsilon$. As in the preceding paragraph, this contradicts $h^0(\omega_D)=p+1$.
\qed

\begin{proposition}
\mylabel{reduced zero-scheme}
There is a non-zero section of $\omega_D$ whose zero-scheme is reduced.
\end{proposition}

\proof
Suppose this would not hold. Then for each non-zero global section $s\in H^0(D,\omega_D)$,
the zero scheme $Z=Z(s)\subset D$ is local, with $\O_Z=\tF[\epsilon]$.
Choose two global sections $s,s'$ that generate $\omega_D$,
and consider the resulting finite flat  morphism $f:D\ra \PP^1 $ of degree $\deg(f)=2p$.
The morphism is not purely inseparable, by our standing assumption $p\neq 2$.
Hence the      fiber $f^{-1}(y)$ is geometrically disconnected for the generic point $y=\eta$.
By \cite{EGA IVc}, Proposition 9.7.8 this holds for almost all point $y\in \PP^1$.
Since $F$ is separably closed, whence infinite, there must be a rational point
$y\in \PP^1$ with $f^{-1}(y)$ geometrically disconnected. On the other hand,
this fiber is the spectrum of $\O_Z=\tF[\epsilon]$, which is geometrically connected, contradiction.
\qed

\section{Non-existence in characteristic three}
\mylabel{Non-existence p=3}
 
Let $F$ be a ground field of characteristic $p=3$, and
$V$ be a regular del Pezzo  surface that is geometrically integral
but geometrically non-normal.  
In Section \ref{Non-reduced ramification p=3}, we already narrowed down the possibilities.
Here our  main result is:

\begin{theorem}
\mylabel{non-existence p=3}
If   a del Pezzo surface $V$ as above exists, the ground field $F$
necessarily has $\pdeg(F)\geq 2$.
\end{theorem}

Seeking a contradiction, we assume throughout that such a del Pezzo surface $V$
exists over a ground field $F$ with $\pdeg(F)\leq 1$.
Without restriction, we may assume that the structure morphism $V\ra\Spec(F)$ is adapted,
and that $F$ is separably closed.
As usual, we use the notation from Section \ref{Regular del Pezzo}, such that  $Y=V\otimes_FK$ 
is a non-normal del Pezzo surface whose normalization $X$ is geometrically normal, with conductor
square
$$
\begin{CD}
R	@>>>	X\\
@VVV		@VVV\\
C	@>>>	Y.
\end{CD}
$$
Recall that the normalization is either the projective plane $X=\PP^2$ or the weighted
projective plane $X=\PP(1,1,3)$. We write $N\subset\Sing(V/F)$ for the divisorial part of the locus of non-smoothness
and $D=N_\red$  for its reduction.

\begin{proposition}
\mylabel{cohomology vanishes}
The conductor curve $C$ is isomorphic to $\PP^1$. Furthermore, we have  $h^1(\shL)=0$ for each $\shL\in\Pic(V)$, 
and $h^0(\O_A)=1$ for each effective Cartier divisor $A\subset V$,
and the restriction map $\Pic(V)\ra\Pic(X)$ is injective.  
\end{proposition}

\proof
According to Proposition \ref{no pair of lines} and Proposition \ref{line}, the ramification divisor $R$ is the split ribbon $\PP^1\oplus\O_{\PP^1}(-1)$.
Furthermore, $C$ is geometrically integral by Proposition \ref{reduced C}, hence  $C=(D_K)_\red$. 
To conclude that $C=\PP^1$, it suffices to check that $h^1(\O_C)=0$.
By Proposition \ref{neron-severi sequence}, we have  $h^1(\O_V)=h^1(\O_C)$.
In light of Maddock's bound \eqref{maddock bound}, it  thus suffices to check $h^1(\O_C)\leq 1$.
Since the canonical map
$$
H^1(D_K,\O_{D_K})\lra H^1((D_K)_\red,\O_{(D_K)_\red})=H^1(C,\O_C)
$$
is surjective, it is enough to verify that one of the two groups on the left is at most one-dimensional.
Indeed, we have $h^1(\O_D)=1$ if the normalization is a weighted projective plane $X=\PP(1,1,3)$,
according to the results in Section \ref{Non-reduced ramification p=3}.

Now suppose that $X=\PP^2$, and consider the integral curve $D\subset V$.
We first check that this is a peculiar curve of higher genus studied in Section \ref{Peculiar higher genus}.
According to the result in Section \ref{Non-reduced ramification p=3}, we have $h^0(\O_D)=1$ and $h^1(\O_D)=4=p+1$. 
This curve is an effective Cartier divisor in a regular surface, so 
all local rings $\O_{D,a}$  are Gorenstein. They are also geometrically unibranch,
since $R \ra D$ is a universal homeomorphism. The curve $D$ contains no rational points, because it lies
in the locus of non-smoothness $\Sing(V/F)$, by Corollary \ref{closed points geom non-reduced}. 
It remains to check that the normalization $\tD$ has
$h^1(\O_\tD)=0$.
The conductor square 
$$
\begin{CD}
A	@>>>	\tD\\
@VVV		@VVV\\
B	@>>>	D
\end{CD}
$$
for the normalization yields an exact sequence
$$
0\ra H^0(\O_D)\ra H^0(\O_\tD)\oplus H^0(\O_B)\ra H^0(\O_A)\ra H^1(\O_D)\ra H^1(\O_\tD)\ra 0.
$$
The field extension $\tF=H^0(\tD,\O_\tD)$ is purely inseparable, of  degree $3^\nu$ for some exponent $\nu\geq 0$.
We have $h^1(\O_\tD)=3^\nu g$ and $h^0(\O_A)=3^\nu a$ for some integers $g\geq0 $ and $a\geq 1$.
Moreover $h^0(\O_B)=3^\nu a/2$, according to Proposition \ref{necessary condition gorenstein}.
The above exact sequence gives $1-(3^\nu+3^\nu a/2)+3^\nu a - 4+ 3^\nu g=0$, in other words
\begin{equation}
\label{genus condition}
3^\nu(a/2+g-1)= 3.
\end{equation}
In particular $0\leq \nu\leq 1$.
If $\nu=0$, the normal curve $\tD$ is geometrically normal, according to \cite{Schroeer 2010}, Theorem 2.3.
Here we use that $\pdeg(F)\leq 1$.
It follows that $D$ is also geometrically normal, because the normalization map
$\tD\ra D$ is birational. This contradicts that the Cartier divisor $D\subset V$
is contained in the locus of non-smoothness $ \Sing(V/F)$, by Corollary  \ref{prime divisors geom non-reduced}.
We thus have $\nu=1$, and equation \eqref{genus condition} only has the following two solutions:
$$
g=0,  a=4\quadand g=1, a=2.
$$
In the former case, $D$ is a peculiar curve of higher genus, and 
Proposition \ref{bound genus C} tells us that  $h^1(\O_C)\leq 1$.
In the latter case,   $\tF=H^0(\tD,\O_{\tD})$ is a field extension
of the form $\tF=F(\beta^{1/3})$ for some non-cube $\beta\in F^\times$.
Note that this field extension is unique up to isomorphism, because now $\pdeg(F)=1$.
Moreover, the cohomology group $H^1(\tD,\O_\tD)$ is a one-dimensional vector space
over $\tF$, and we fix a cohomology class $\alpha\neq 0$.
Now consider the base-change $\tF\otimes_FK=K[t]/(t^3)$.
The resulting nilpotent element $t\in H^0(\tD_K,\O_{\tD_K})$ vanishes on the reduction $C=(\tD_K)_\red$,
and in turn the image of $\alpha$ in the cohomology group $H^1(C,\O_C)$
is annihilated by $t$. It follows that the image of the surjective restriction map 
$$
H^1(D_K,\O_{D_K})\lra H^1((D_K)_\red,\O_{(D_K)_\red})=H^1(C,\O_C)
$$
is at most one-dimensional, and we get $h^1(\O_C)\leq 1$ again.

Now let $\shL$ be an invertible sheaf on the del Pezzo surface $V$.
The short exact sequence $0\ra\shL_Y\ra\shL_X\oplus\shL_C\ra\shL_R\ra 0$ from
the conductor square yields an exact sequence
$$
H^0(\shL_X)\oplus H^0(\shL_C)\ra H^0(\shL_R)\ra H^1(\shL_Y)\ra H^1(\shL_X)\oplus H^1(\shL_C)\ra H^1(\shL_R).
$$
Now we use some general facts:
For each weighted projective space $\PP=\PP(q_0,\ldots,q_n)$ with $\gcd(q_0,\ldots,q_n)=1$,
all tautological sheaves $\O_{\PP}(n)$ have trivial cohomology in intermediate degrees $0<i<n$,
according to  \cite{Dolgachev 1981}, Theorem in Section 1.4. 
Moreover, it is well-known that $\O_\PP(1)$ generates
the  group of reflexive rank-one sheaves
$\APic(\PP)$, compare  \cite{Cox; Little; Schenck 2011}, Section 4.1.
It follows that $H^1(X,\shL_X)=0$, and also that the restriction map $H^0(X,\shL_X)\ra H^0(R,\shL_R)$
is surjective.

Recall that $R=\PP^1\oplus\O_{\PP^1}(-1)$ is the split ribbon, hence $R_\red\ra C$ is an isomorphism.
As a consequence, the pull-back  map $H^1(\shL_C)\ra H^1(\shL_R)$ is injective. 
In light of the above long exact sequence, the group $H^1(Y,\shL_Y)$ vanishes, thus  $h^1(\shL)=0$.
Applying this to $\shL=\O_X$ we infer that the Picard scheme $\Pic^0_{V/F}=0$,
such that $\Pic(V)=\NS(V)$. 
The injectivity of $\Pic(V)\ra\Pic(X)$ follows from Proposition \ref{neron-severi sequence}.
\qed

\medskip
For $n\geq 0$, Serre Duality gives  $h^2(\omega_V^{\otimes -n})=0$, because $\omega_V^{\otimes n+1}$ is antiample,
so Riemann--Roch yields 
\begin{equation}
\label{global sections}
h^0(\omega_V^{\otimes -n}) = \frac{(-nK_V)^2 - (-nK_V)\cdot K_V}{2} +\chi(\O_V) = \frac{n(n+1)}{2}K_V^2 + 1.
\end{equation}
To proceed, we consider the invertible sheaves 
$$
\shL=\begin{cases}
\omega_V^{\otimes -2} & \text{if $X=\PP^2$;}\\
\omega_V^{\otimes -1} & \text{if $X=\PP(1,1,3)$.}
\end{cases}
$$
which both have $h^0(\shL)=4$.

\begin{proposition}
The  ample invertible sheaf $\shL$ is globally generated.
\end{proposition}

\proof
Suppose first that $X=\PP^2$. Consider first the invertible sheaf $\omega_X^{\otimes -1}$.
Since it generates the Picard group, the effective Cartier divisor $A\subset V$ attached to any   non-zero   $s\in H^0(X,\omega_X^{\otimes -1})$
is integral. Moreover, the self-intersection number is  $K_V^2=1$,  and $h^0(\omega_V^{\otimes -1})=2$. Thus there is another non-zero global section $s'$
such that the corresponding effective Cartier divisor $A'$ has $A\cap A'=\Spec\kappa(v)$ for some rational point $a\in V$.
Moreover, the local ring $\O_{A,a}$ is regular.
This already shows that $\Bs(\shL)\subset\{a\}$.
The restriction map $H^0(V,\shL)\ra H^0(A,\shL_A)$ is surjective, because $H^1(X,\shL(-A))=0$.

It thus suffices to check that $a\in A$ is not a base point for the restriction $\shL_A$. This would hold
if $h^1(\shL_A(-a))=h^0(\omega_A\otimes\shL^{\otimes -1}(a))$ vanishes.
The adjunction formula shows that the degree of $\omega_A\otimes\shL^{\otimes -1}_A(a)$ is
$$
(K_V+A)\cdot A - (\shL\cdot A) + 1 = 0 - 2+1=-1.
$$
In turn, the invertible sheaf $\omega_A\otimes\shL^{\otimes -1}(a)$ on the integral scheme $A$ has no non-zero global sections.
This shows that $\shL$ is globally generated.
In the case $X=\PP(1,1,3)$, one argues as above with $\omega_V^{\otimes-1}=\shL$.
\qed

\medskip
Consider the finite morphism $V\ra \PP^3$ resulting from
the linear system $H^0(V,\shL)$. Its image is an integral surface $V'\subset \PP^2$,
and we get a finite surjection  $f:V\ra V'$. There are only few possibilities:

\begin{proposition}
\mylabel{finite morphism birational}
The finite morphism $f:V\ra V'$ is birational.
The image $V'\subset\PP^3$ is a quartic surface if $X=\PP^2$,
and  a cubic surface if $X=\PP(1,1,3)$.
\end{proposition}

\proof
This follows from $(\shL\cdot\shL)=\deg(f)\cdot\deg(V')$ and $\deg(V')>1$.
In case $X=\PP(1,1,3)$, the selfintersection is $(\shL\cdot\shL)=3$, and the assertion   is immediate.

Now suppose that $X=\PP^2$. Seeking a contradiction, we assume that $f:V\ra V'$ is a double covering
over some quadric surface $V'$ in $\PP^3=\Proj F[T_0,\ldots,T_3]$.
As $p\neq 2$, we may assume that $V'$ is defined by the homogeneous equation
$T_0^2+\ldots +T_n^2=0$ for some $0\leq n\leq 3$. 
We have $n\neq 0$, because $V'$ is reduced.
Moreover  $n\neq 1$ because $V'$ is irreducible.
In case $n=3$, we would have $V'\simeq\PP^1\times\PP^1$,
in contradiction to $\rho(V)=1$. Thus $n=2$, and   $V'=\PP(1,1,2)$ is isomorphic to the contracted Hirzebruch
surface with numerical invariant $e=2$.

To simplify notation,  set   $\PP=\PP(1,1,2)$ and write
 $\O_\PP(1)$ for the tautological sheaf 
on the weighted homogeneous spectrum $\PP=\Proj k[T_0,T_1,T_2]$.
This sheaf generates the group $\APic(\PP)$ of isomorphism classes of reflexive rank-one sheaves.
Note that $\O_\PP(2)$ corresponds to the restriction of the invertible sheaf $\O_{\PP^3}(1)$ to the quadric surface $\PP=V'$,
and that $\O_\PP(1)$ comes from  the fiber on the resolution of singularities,
which is a Hirzebruch surface. The latter shows that 
the rational selfintersection number of $\O_\PP(1)$ is $1/2$, hence
the induced map $f^*:\APic(\PP)\ra\Pic(V)$ under the double covering   is bijective.
In particular, $\omega_V^{\otimes -1}$ is the reflexive hull of  the rank-one sheaf $f^*\O_\PP(1)$.

Since $p\neq 2$, we may view  $V$ as the relative spectrum of $\shA=\O_\PP\oplus\O_\PP(-d)$ for some integer $d$.
The multiplication law for $\shA$ is given by some homomorphism  $\varphi:\O_\PP(-2d)\ra\O_\PP$, which can also be viewed as a global
section $\varphi\in\Gamma(\PP,\O_\PP(2d))$.
We have $d>0$ because $h^0(\O_V)=1$.
Moreover, the integer $d$ is odd because $V$ is regular whereas  $\PP$ is singular, such that the finite morphism $f:V\ra\PP$ is not flat.
We have  $\omega_\PP=\O_\PP(-4)$, and $H^2(V,\O_V)=H^2(\PP,\shA)$ is Serre dual to $\Hom(\shA,\omega_\PP)=H^0(\PP,\O_\PP(-4)\oplus\O_\PP(d-4))$,
hence  $d\leq  3$.
Furthermore, for each integer $t$ the Projection Formula yields $f_*(\omega_V^{\otimes-t})=\O_\PP(t)\oplus\O_\PP(t-d)$, thus
$$
2= h^0(\omega_V^{\otimes-1}) = h^0(\O_\PP(1)) + h^0(\O_\PP(1-d)).
$$
Using $h^0(\O_\PP(1))=2$ we conclude $1-d<0$. Thus  $d=3$ is the only remaining possibility.

Summing up, $\shA=\O_\PP\oplus\O_\PP(-3)$, and the multiplication law comes from some non-zero
section $\varphi\in\Gamma(\PP,\O_\PP(6))$. The latter defines an effective Cartier divisor $B\subset\PP$.
In light of the multiplication law on $\shA$, the schematic preimage takes
the form $f^{-1}(B)=2A$ for some effective Cartier divisor $A\subset V$.
Let $U\subset \PP$ be the complement of $B\cup\Sing(\PP)$.
The scheme $U$ is smooth and the  double covering $f:V\ra\PP$ is \'etale over $U$,
which implies $D\subset A$. We actually have $D=A$, because $\O_V(D)=\omega_V^{\otimes-3}=\O_V(A)$.
In turn, the induced morphism $f:D\ra B$ is birational.  
Since $D$ is singular,  the integral curve $B$ is singular as well.

We claim that $B$ does not pass through the singular point $y\in \PP$, where the local ring $\O_{\PP,y}$ is a rational double point
of type $A_1$, whose local fundamental group $\pi_1^\loc(\O_{\PP,y})$ is cyclic of order two.
Let $x\in f^{-1}(y)$ be the preimage, and suppose that $y\in B$.
Applying \cite{Ito; Schroeer 2015}, Proposition 2.3 for the inclusion of local rings
$\O_{\PP,y}\subset\O_{V,x}$, we conclude that $\pi_1^\loc(\O_{\PP,y})$ is trivial, contradiction.

Fix a closed point $b\in B$ where the local ring $\O_{B,b}$ is singular.
By the preceding paragraph we have  $\varphi_b\in\maxid_b^2\O_{\PP,b}(2d)$.
Now choose an identification  $\O_{\PP,b}(2d)=\O_{\PP,b}$ and a  regular system of parameters $x,y\in\O_{\PP,b}^\wedge$,
and write $\varphi_b=g(x,y)$ as a formal power series with neither constant nor linear terms.
For the point $a\in f^{-1}(b)$, we get  $\O_{V,a}^\wedge=\kappa(b)[[x,y,z]]/(z^2-g)$.
Since $z^2-g\in(x,y,z)^2$, the local ring $\O_{V,a}$ is singular, 
contradiction.
\qed

\medskip
\emph{Proof for Theorem \ref{non-existence p=3}.}
Consider first the case that $X=\PP^2$.
Then $\shL=\omega_V^{\otimes-2}$ defines a   finite birational  morphism $f:V\ra V'$
onto some quartic surface $V'\subset\PP^3$, with $\shL=f^*(\O_{V'}(1))$.
The Adjunction Formula gives $\omega_{V'}=\O_{V'}$, and consequently $\omega_V=\O_V(-A)$
where $A\subset V$ is the ramification divisor.
We obtain a contradiction by showing that $f:V\ra V'$ does not ramify along $A$.
Since $\omega_V^{\otimes-1}$ generates the Picard group and $V$ is regular, the curve 
$A$ must be integral.
The long exact sequence for $0\ra\omega_V\ra\O_V\ra\O_A\ra 0$, together with Proposition \ref{cohomology vanishes}
yields $h^0(\O_A)=1$, and the Adjunction formula gives $\omega_A=\O_A$,
whence $h^1(\O_A)=1$. The restriction $\shL|A$ has degree $(\shL\cdot A)=2$,
and Riemann--Roch gives $h^0(\shL|A)=2$. So the globally generated sheaf $\shL|A$ defines
a morphism $A\ra\PP^1$ of degree two. Consequently,
the morphism  $f:V\ra V'$ does not ramify along $A$, contradiction.

Now consider the case $X=\PP(1,1,3)$.
Then $\shL=\omega_V^{\otimes -1}$ defines a finite birational morphism $f:V\ra V'$
onto some cubic surface $V'\subset\PP^3$. The Adjunction Formula gives
$\omega_{V'}=\O_{V'}(-1)$, and we have $\omega_V=f^*(\omega_{V'})\otimes \O_V(-R)$,
where $R\subset V$ is the ramification divisor.
Here we have  $(\omega_V\cdot\omega_V)=3$ and $(\omega_{V'}\cdot\omega_{V'})=3$.
Together with $\Pic(V)=\ZZ$, one easily infers that the ramification divisor is empty,
such that $V=V'$ itself is a cubic surface.
For each anticanonical divisor $A\subset V$, we have $h^0(\O_A)=1$ and $\omega_A=\O_A$,
whence also $h^1(\O_A)=1$. Riemann--Roch gives $h^0(\shL|A)=3$, showing
that $A$ becomes a quadric curve inside some  plane $\PP^2\subset\PP^3$ under the embedding $V\subset\PP^3$.

To proceed, we consider the special case that $A=D$ is the reduced locus of non-smoothness,
which is geometrically non-reduced. The cubic curve $D\subset\PP^2$ is given by 
some irreducible homogeneous polynomial $\bar{P}\in F[T_0,T_1,T_2]$ of degree three.
Since $\Sing(D/F)=D$, all partial derivatives $\bar{P}_i=\partial \bar{P}/\partial T_i$
are multiplies of $\bar{P}$.  For degree reasons we must have $\bar{P}_i=0$,
and thus our polynomial takes the form $\bar{P}=\sum_{i=0}^2\lambda_iT_i^3$ for some
scalars $\lambda_i\in F$. Thus $D\subset\PP^2$ is a $p$-Fermat hypersurface,
which were studied in \cite{Schroeer 2010}, \S 3. 

In turn,   the cubic surface $V\subset\PP^3$ is given by a homogeneous polynomial of the form
$$
P(T_0,T_1,T_2,T_3)= \sum_{i=0}^3\lambda_iT_i^3  + T_3^2L  + T_3Q, 
$$
for some linear term $L=L(T_0,T_1,T_2)$ and some quadratic term $Q=Q(T_0,T_1,T_2)$.
Here the curve $D\subset V$ is given by $T_3=0$.
Since $D\subset\Sing(V/F)$, all partial derivatives $\partial P/\partial T_i$ are multiples
of $T_3$. Using $\partial P/\partial T_3 = 2T_3 L + Q$, we infer $T_3|Q$ and thus may assume $Q=0$.

The curve $D$ is integral with $h^0(\O_D)=1$, and our ground field has $\pdeg(F)\leq 1$.
According to \cite{Schroeer 2010}, Theorem 2.3 there must be some closed point $a\in D$ where
the local ring $\O_{D,a}$ is singular. This point has homogeneous coordinates 
$(\alpha_0:\alpha_1:\alpha_2:0)$
with $\alpha_i\in F^\alg$. Without restriction, we may assume that $\alpha_0=1$.
Then the homogenization $\bar{P}T_0^{-3}$ lies in   the square of 
the maximal ideal $\O_{\PP^2,a}$.   In light of the form $P=\sum\lambda_iT_i^3+T_3^2L$, it follows that $PT_0^{-3}$ lies in the 
square of the maximal ideal of $\O_{\PP^3,a}$, thus the local ring $\O_{V,a}$ is singular,
contradiction.
\qed

\section{Del Pezzo surfaces and fibrations}
\mylabel{Del Pezzo fibrations}

We now can state and proof the main result of our paper:

\begin{theorem}
\mylabel{main pdeg}
Let $F$ be a ground field of characteristic $p>0$, with $p$-degree $\pdeg(F)\leq 1$.
Then every regular del Pezzo surface $V$ over $F$ with Picard number $\rho(V)=1$ is geometrically normal.
\end{theorem}

\proof 
Seeking for a contradiction, assume there exists a del Pezzo surface $V$  with Picard number $\rho(V)=1$ that  is geometrically non-normal. 
Since $\pdeg(F)\leq 1$, we can apply \cite{Schroeer 2010}, Theorem 2.3 and conclude that $V$ is geometrically integral. 
According to Theorem \ref{higher p}, we must have $p\leq 3$.
Choose some finite separable field extension $F\subset F'$ and some finite purely inseparable
extension $F\subset K'$ so that $V'=V\otimes_FF'$ and $F'\subset F'\otimes_FK'$ are adapted,
as explained in Section \ref{Geometrically non-normal}.

The geometry of  $V'$ was described in the tables of  Sections 
\ref{Smooth ramification}, \ref{Non-reduced ramification p=2} and \ref{Non-reduced ramification p=3}.
In all but one case, both surfaces $V$ and $V'$ have Picard number $\rho=1$.
The case with $\rho(V')>1$ is treated in Theorem \ref{main smooth ramification},
and then $\rho(V')=2$.
In this case, we also have $\rho(V)=2$, as  remarked in the end of Section \ref{Smooth ramification},
contradiction.

Thus we may assume from the start that $F$ is separably closed, and that we have a finite purely inseparable
field extension $F\subset K$ such that $V$ is adapted.
The possibilities for $V$ with arbitrary $p$-degree were first narrowed down in 
  Propositions \ref{D smooth R}, \ref{pdeg ramification non-reduced} and \ref{char 3 irregularity}.
Finally, the  non-existence  of the remaining possibilities for $\pdeg(F)=1$
follow from  Theorems \ref{main smooth ramification}, \ref{non-existence without irregularity}, \ref{nonexistence 22} and \ref{non-existence p=3}.
\qed

\medskip
The previous result is   sharp: In the rest of this section we will describe two new examples $V$ and $W$ of regular del Pezzo surfaces defined over any imperfect field 
$F$ in characteristic $p=2$ with   the following   properties:
\begin{enumerate}
\item  $W$ geometrically non-regular, with Picard number  $\rho(W)=1$;
\item  $V$ geometrically non-normal, with Picard number  $\rho(V)=2$.
\end{enumerate}

\medskip
Recall that $z^p-xy=0$ is the equation, in normal form, for 
the rational double point of type $A_{p-1}$.
In what follows  we shall use  twisted forms of the corresponding local ring 
that are regular. It is very easy to determine when this happens:

\begin{lemma}
\mylabel{regular twisted form}
Let $R$ be a   local ring that is essentially of finite type over $F$.
Suppose there is a finite purely inseparable extension $F\subset K$ and an isomorphism
$$
\hat{R}\otimes_FK\simeq K[[x,y,z]]/(z^p-xy).
$$
Then the degree $d=[\kappa:F]$ of the residue field $\kappa=R/\maxid_R$
is either $d=1$ or $d=p$.  The latter holds if and only if the local  $R$ is regular.
\end{lemma}

\proof
Let $J\subset R$ be the jacobian ideal. Then $M=R/J$ has finite length.
Computing with $K[[x,y,z]]/(z^p-xy)$, one sees that its vector space dimension is $\dim_F(M)=p$.
Since there is a filtration on $M$ whose subquotients are copies of the unique
simple $R$-module $\kappa=R/\maxid_R$, we conclude that either $d=1$ or $d=p$.
In the latter case, we must have  $M\simeq \kappa$, so the jacobian ideal coincides with the maximal
ideal. Computing again with $K[[x,y,z]]/(z^p-xy)$, we see that $M$ has finite projective dimension.
By the homological characterization of regularity,  the local ring $R$ must be regular.
\qed

\medskip
If a local ring $R$ is regular and satisfies the assumption of the lemma,
we call it  a   \emph{regular twisted forms of the rational double point of type $A_{p-1}$}.
See \cite{Schroeer 2008} for more on this.
If the ground field $F$ is imperfect, such rings indeed do exist:
Choose a scalar $\beta\in F$ that is not a $p$-th power, and consider the
local ring $R$ coming from the residue class ring  $K[x,y,z]/(z^p-xy-\beta)$ with respect to the maximal ideal $\maxid=(x,y)$.
If $\beta$ becomes a $p$-th power  in the field extension $F\subset K$,
then the substitution $z=z'+\beta^{1/p}$ yields the desired isomorphism to the
rational double point of type $A_{p-1}$.

Suppose that  $R$ is a regular twisted form of a rational double point of type $A_{p-1}$,
with residue field $\kappa=R/\maxid_R$. Let 
$r:X\ra\Spec(R)$ be the blowing-up of the reduced closed point, 
and write $D\subset X$ for the resulting exceptional divisor.

\begin{proposition}
\mylabel{blowing-up non-normal}
As a scheme, the   exceptional divisor     $D$ is isomorphic to $\PP^1_\kappa$, and its  selfintersection number
is $D^2=-p$. The surface $X$ is regular but not geometrically normal. For $p=2$, the locus of non-smoothness $N=\Sing(X/F)$
 is given by  $N=2D$. Otherwise, it contains two embedded associated points, and its divisorial part is $N_\text{\rm div}=D$.
\end{proposition}

\proof
Regarding the blowing-up as a closed subscheme
$X\subset \PP^1\otimes_FR$ as in \cite{SGA 6}, Expos\'e VII, Proposition 1.8, we see that the exceptional divisor is given by $D=\PP^1_\kappa$,
and it follows that the scheme $X$ is regular.
To see that it is not geometrically normal, we may compute with $A=K[[x,y,z]]/(z^p-xy)$.
According to the proof of Lemma \ref{regular twisted form}, the center of the blowing-up is defined by
the jacobian ideal, which   is generated by $x,y\in A$. The $x$-chart of the blowing-up
is given by the variables $x,y/x,z$ modulo the relation 
$z^p-x^2(y/x)=0$, and the exceptional divisor has equation $x=0$.
The module $\Omega^1_{A/K}$  is generated by the differentials $dx,d(y/x),dz$ modulo
the relation $2x(y/x)dx + x^2d(y/x)=0$. Consequently, the locus of non-smoothness is defined by the ideal $(2x(y/x),x^2)\subset A$,
and the assertion on  $N=\Sing(X/F)$ follows.
\qed

\medskip
Now fix a scalar $\beta\in F$, and 
consider the hypersurfaces $W=W_\beta$ inside $\PP^3$ of degree $\deg(W)=p$ over the ground field $F$ given by the homogeneous equation
\begin{equation}
\label{hypersurface equation}
T_2^p+T_0T_1T_3^{p-2} - \beta T_3^p =0.
\end{equation}
These are projective Gorenstein surfaces  with $h^0(\O_W)=1$ and $h^1(\O_W)=0$.
Furthermore, the dualizing sheaf is $\omega_W=\O_W(p-4)$.
Note that  $W$ is a del Pezzo surface if and only if $p\leq 3$. In any case, $K_W^2=(p-4)^2p$. In characteristic two, this
becomes $K_W^2=8$.

Clearly, the substitution $T_2=T'_2+\beta^{1/p}T'_3$ shows that the subschemes $W_\beta,W_0\subset\PP^3$ become projectively
equivalent over the field extension $K=F(\beta^{1/p})$. In particular, each $W=W_\beta$ is a twisted form of $W_0$.

\begin{proposition}
\mylabel{picard group}
In characteristic $p=2$, the restriction map $\Pic(\PP^2)\ra\Pic(W)$ is bijective.
For $p\geq 3$, this holds  true up to torsion in $\Pic(W)$.
\end{proposition}

\proof
Since $H^1(W,\O_W)=0$, the Picard scheme is zero-dimensional, thus the canonical map $\Pic(W)\ra\NS(W)$
is bijective, and this group is finitely generated.

Consider the hyperplane section $H=W\cap V_+(T_3)\subset W$.
Clearly, the    affine scheme $U=W\smallsetminus H$ is isomorphic to the spectrum of 
$k[x,y,z]/(z^p-xy)$. This can be regarded as an affine toric variety,
and thus its Picard group vanishes (\cite{Kempf et al 1973}, Theorem 9 on page 28, compare also
\cite{Cox; Little; Schenck 2011}, Proposition 4.2.2).
Consequently, every Cartier divisor $D\in\Div(W)$
is linearly equivalent to some Cartier divisor supported by the effective Cartier divisor $H\subset W$.

Using Equation \eqref{hypersurface equation} and treating the case $p=2$ and $p\geq 3$ separately, one easily sees that the scheme $H$ is irreducible.
Let $\eta\in H$ be the generic point,
and set $n=\length(\O_{H,\eta})$. Then $(D\cdot H)= n \deg(\shL|H_\red)$, where $\shL=\O_W(D)$.
Since $W$ is projective, the intersection form on $\NS(W)$ modulo its torsion subgroup is non-degenerate
(\cite{SGA 6}, Expos\'e XIII, Corollary 7.4),
and we infer that every Cartier divisor $D\in \Div(W)$ is numerically equivalent to
some rational multiple of $H$. Since $H^2=p$, the element $H\in\Pic(W)$ is primitive,
and thus $H\in \NS(W)$ is a generator up to torsion.
In case $p=2$, the hyperplane $H$ is isomorphic to the quadric curve in $\PP^2$ given by $T_2^2+T_0T_1=0$,
which is integral, and it follows that $H\in\NS(W)$ indeed is a generator.
\qed

\medskip
Consider the closed point $x=(0:0:\beta^{1/p}:1)$.
If the scalar $\beta\in F$ is a $p$-th power, then   $x\in W$ is a rational point,
and the local ring $R=\O_{W,x}$ is a rational double point of type $A_{p-1}$.
Now suppose that $\beta\in F^\times $ is not a $p$-th power.
Then $x\in W$ is   non-rational, with residue field $K=F(\beta^{1/p})=\kappa(x)$. By Lemma \ref{regular twisted form}, the
local ring $R=\O_{W,x}$ is a regular twisted form of the rational double point
of type $A_{p-1}$.  

The two hyperplanes $H_0=V_+(T_0)$ and $H_1=V_+(T_1)$ yield two linearly equivalent Cartier divisors $C_i=W\cap H_i$
whose  intersection  $C_1\cap C_2$ consists of this  point $x\in W$, viewed as a reduced subscheme.
Let $r:V\ra W$ be the blowing-up with reduced center $x\in W$.
The exceptional divisor $D=r^{-1}(x)$ is a copy of $\PP^1_K$, and the strict transforms of the Cartier divisors $C_1,C_2$
yield a fibration $\varphi:V\ra \PP^1$, for which $D\subset X$ is horizontal, of relative degree $\deg(D/\PP^1)=p$.

\begin{proposition}
\mylabel{blowing-up}
Suppose   $p=2$. Then $V$ is a regular del Pezzo surface of degree $K_V^2=6$.
Its Picard number is $\rho(V)=2$,  and the Picard group is freely generated by the fiber  $F=\varphi^{-1}(\infty)$
and the exceptional divisor $D=r^{-1}(x)$. The Gram matrix is 
$(\begin{smallmatrix}0& 2\\2&-2\end{smallmatrix})$, and we have $K_V=-(D+2F)$. Furthermore,
 $\Sing(V/F)=D$, and the locus of non-smoothness is $N=2D$.
\end{proposition}

\proof
Write $K_{V/W}=nD$ for some integer $n\in\ZZ$. The Adjunction Formula gives
$$
-2(n+1) = (n+1)D^2= (K_{V/W}+D)\cdot D = \deg(K_D) = -2\chi(\O_{\PP^1_K}) = -4,
$$
thus $n=1$. Consequently $K_V^2=K_W^2+K_{V/W}^2=6$. Using $r^{-1}(H_0)=F+D$ and $K_W=-2H_0$,
we infer $K_V=-(D+2F)$.
By Proposition \ref{picard group}, the Picard group of $W$ is
generated by $H_0$, so $\Pic(V)$ is generated by the strict transform $F$ and the exceptional divisor $D$.
We compute $(K_V\cdot D) = D^2 = -2$ and 
$$
(K_V\cdot F) = K_V\cdot (r^{-1}(H_0) - D) = (K_W\cdot H_0) - D^2 = -4+2=-2.
$$
The cone of curve is generated by the two extremal rays coming from $D,F\subset V$.
By the Nakai Criterion (\cite{Hartshorne 1970} Chapter I, Theorem 5.1), the anticanonical divisor  $-K_V$ is ample.
\qed

\medskip
Let us examine the fibration $\varphi:V\ra \PP^1$ for $p=2$ in more detail.
The fiber over a rational point $(\lambda_0:\lambda_1)\in\PP^1$, say with $\lambda_0=1$,
can be regarded as the zero-scheme inside $\PP^2$ for the equation
$$
T_2^2+\lambda_1T_1^2 - \beta T_3^2 =0.
$$
According to \cite{Schroeer 2010}, Theorem 3.3  this is  regular 
provided that the field extension $F\subset F(\lambda_1^{1/2},\beta^{1/2})$
has degree four. In our situation, this means $\lambda_1^{1/2}\not\in K$.
On the other hand, if the condition does not hold, the curve $C=\varphi^{-1}(\lambda_0:\lambda_1)$ is at least integral, and its
normalization is given by the conductor square 
$$
\begin{CD}
\Spec (K) 	@>>> 	\PP^1_K\\
@VVV			@VVV\\
\Spec(F)	@>>>	C.
\end{CD}
$$
Roughly speaking, a $K$-rational point on $\PP^1_K$ is replaced by an $F$-rational point on $C$.
In any case, the fiber is   a twisted form of the ribbon $\PP^1\oplus\O_{\PP^1}(-1)$.
In particular,  all fibers are proper curves $C$ with $h^0(\O_C)=1$ and $h^1(\O_C)=0$.

Let us call a proper morphism between integral scheme $q:Z\ra B$ a \emph{genus-zero fibration} if
$q_*(\O_Z)=\O_B$ and the generic fiber $C=Z_\eta$ is a curve with $h^0(\O_C)=1$ and $h^1(\O_C)=0$.
If the generic fiber is smooth, we say that the fibration is a \emph{ruling}.
Otherwise, we call it a \emph{quasiruling}, in analogy to the
quasielliptic fibrations.

The regular del Pezzo surfaces $W$ and $ V$ constructed above actually occur as generic fiber of some \emph{del Pezzo fibrations}.
Fix some algebraically closed ground field $k$ of characteristic $p=2$.
Consider the divisor $Y\subset\PP^1\times\PP^3$ of bidegree $\deg(Y)=(1,2)$  given by the bihomogeneous equation
\begin{equation}
\label{bihomogeneous equation}
S_0(T_2^2+T_0T_1) -  S_1T_3^2 = 0,
\end{equation}
where the $S_0,S_1$ and $T_0,\ldots,T_3$ are the homogeneous coordinates for $\PP^1$ and   $\PP^3$,
respectively. Write $\O_{\PP^1\times\PP^3}(d_1,d_2)=\pr_1^*(\O_{\PP^1}(d_1))\otimes \pr_2^*(\O_{\PP^1}(d_2))$ to simplify notation.
The Adjunction Formula gives 
$$
\omega_Y= \O_Y(-1,-2)=\O_Y(-Y),
$$
so our $Y$ is a Fano threefold of degree $-K_Y^3=(A+2B)^4=A\cdot(2B)^3=8$, where $A,B\subset \PP^1\times\PP^3$
are divisors corresponding to $\O_{\PP^1\times\PP^3}(1,0)$ and $\O_{\PP^1\times\PP^3}(0,1)$

\begin{proposition}
\mylabel{scheme Y}
The Fano threefold $Y$ is normal. Its  singular scheme $\Sing(Y)$ is given by $V_+(S_0,T_3,T_2^2+T_0T_1)$,
which  is isomorphic to $\PP^1$.
Moreover, for each closed point $a\in \Sing(Y)$, the corresponding complete local ring $\O_{Y,a}^\wedge$
is isomorphic to $k[[x,y,z,w]]/(z^2-xy)$. Furthermore, the blowing-up $Y'\ra Y$
with reduced center $\Sing(Y)$ is a crepant resolution of singularities.
\end{proposition}

\proof
The assertion on $\Sing(Y)$ follows from computing the partial derivatives
in \eqref{bihomogeneous equation}. We see that $Y$ is regular in codimension one.
Being a hypersurface, it is also  Cohen--Macaulay, and Serre's Criterion ensures
that $Y$ is normal. 

Each singular point $a\in Y$ must have  homogeneous coordinates $(0:1)$ in the   homogeneous coordinates with respect to the first factor $\PP^1$.
Without restriction, we may assume that its homogeneous coordinates for the second factor $\PP^3$
are of the form $(1:\lambda^2:\lambda:0)$. Setting 
$$
z=T_3/T_0,\quad x=S_0/S_1,\quad y=(T_2/T_0)^2-(T_1/T_0)\quadand w=T_2/T_0
$$
yields the assertion on $\O_{Y,a}^\wedge$. This is basically the equation for a rational double point of
type $A_1$, up to the additional variable $w$, so the blowing-up of the reduced singular locus
gives a crepant resolution of singularities.
\qed

\medskip
The singularities $a\in Y$ are very simple special cases of the so-called \emph{compound du Val singularities}
introduced in \cite{Reid 1980}. In particular, these are rational and canonical singularities.

The first projection induces a morphism $f:Y\ra\PP^1$.  
All fibers $f^{-1}(\lambda_0:\lambda_1)$ are hypersurfaces in $\PP^3$, and it follows that
the map $\O_{\PP^1}\ra f_*(\O_Y)$ is bijective, and $R^1f_*(\O_Y)$ vanishes.
Moreover, the dualizing sheaf $\omega_Y=\O_Y(-1,-2)$ is $f$-antiample. Set $F=k(\PP^1)$.

\begin{proposition}
\mylabel{first Fano--Mori  fibration}
The generic fiber $Y_F$  is   the regular geometrically normal but non-smooth del Pezzo surface $W=W_\beta$ of degree $K_W^2=8$
and Picard number $\rho(W)=1$ constructed above, with scalar $\beta=S_1/S_0\in F$.
The closed fiber $f^{-1}(\infty)$  over the point $\infty=(0:1)$  is isomorphic to the split ribbon
$\PP^2\oplus\O_{\PP^2}(-1)$.
All other closed fibers are isomorphic to the normal del Pezzo surface $W=W_0$.
\end{proposition}

\proof
This follows immediately from the defining equation \eqref{bihomogeneous equation}, and we merely
need to clarify the structure of the closed fiber  $ f^{-1}(\infty)\subset Y$. By definition, this is 
the hypersurface $H\subset \PP^3$ given by $T_3^2=0$,
thus $H_\red=\PP^2$.
The structure sheaf sits in a short exact sequence
$0\ra\shL\ra \O_H\ra \O_{\PP^2}\ra 0$,
with $\shL=\O_{\PP^2}(-1)$. In other words, $H$ is  a ribbon on $\PP^2$ with invertible sheaf $\shL$.
According to \cite{Bayer; Eisenbud 1995}, Theorem 1.2,  these are classified by elements in
$\Ext^1(\Omega^1_{\PP^2/k},\shL) = H^1(\PP^2,\Theta_{\PP^2/k}(-1))$.
But this cohomology group vanishes, which easily follows from the long exact cohomology sequence for the  Euler sequence
$0\ra\O_{\PP^2}(-1)\ra \O_{\PP^2}^{\oplus 3} \ra\Theta_{\PP^2/k}(-1)\ra 0$.
Thus $H=f^{-1}(\infty)$ is isomorphic to the split ribbon $\PP^2\oplus\O_{\PP^2}(-1)$.
\qed

\medskip
Note that $W_0$ is a contracted Hirzebruch surface with invariant $e=2$.
Now consider the closure
$$
Z=Y\cap  V_+(T_0,T_1) = V_+(T_0,T_1, S_0T_2^2 - S_1T_3^2)\subset Y 
$$
of $\Sing(Y_F/F)$ inside the total space $Y$.
Clearly, $Z\subset Y$ is a smooth complete intersection
of two effective Cartier divisors. We may regard this as a divisor in $\PP^1\times\PP^1$, and the second projection gives
an isomorphism $\pr_2:Z\ra\PP^1$, whereas the first projection $\pr_1:Z\ra\PP^1$ can be seen
as the relative Frobenius morphism. Note  that $Z\subset\Reg(Y)$, which is also a consequence of 
Proposition \ref{scheme Y}. 

Now let $ X\ra Y$ be the blowing-up with reduced center $Z\subset Y$.
This gives a commutative diagram
$$
\begin{CD}
X'	@>>>	X\\
@VVV		@VVV\\
Y'	@>>>	Y,
\end{CD}
$$
where the horizontal morphisms are crepant resolutions  described in Proposition \ref{scheme Y}.
Now consider the induced fibrations $f:Y'\ra \PP^1$ and  $g:X'\ra\PP^1$.
Using Proposition \ref{blowing-up}, we immediately get:
 
\begin{theorem}
\mylabel{main examples}
The morphism $g:X'\ra\PP^1$ is a fibration from a smooth threefold, 
whose generic fiber is a regular and geometrically non-normal del Pezzo surface with Picard number two.
The morphism $f:Y'\ra\PP^1$ is a fibration from a smooth threefold 
whose generic fiber is a regular, geometrically normal and geometrically non-regular del Pezzo surface with Picard number one.
\end{theorem}

\section{Mori  fiber spaces}
\mylabel{Mori  fiber spaces}

We now examine our results on del Pezzo surfaces in the context of the minimal model program.
We refer to the monographs of Koll\'ar and Mori \cite{Kollar; Mori 1998} and Matsuki \cite{Matsuki 2002}
for   general expositions, and the book of Koll\'ar \cite{Kollar 2013} on the singularities occurring in the
minimal model program. In this section, $k$ denotes an algebraically ground field of
characteristic $p\geq 0$.

Let us call a morphism $f:Z\ra B$ between proper normal integral scheme a \emph{fibration}
if $f_*(\O_Z)=\O_B$.
A fibration is  called a \emph{Mori fibration} if the following four properties hold:
\begin{enumerate}
\item The generic fiber $V=Z_\eta$ has dimension $\dim(V)\geq 1$.
\item The   local rings $\O_{Z,a}$, $a\in Z$ are $\QQ$-factorial  klt  singularities.
\item The Picard numbers satisfy  $\rho(Z)=\rho(B)+1$.
\item The $\QQ$-Cartier  divisor $K_Z$ is relatively antiample for $f:Z\ra B$.
\end{enumerate}
In this situation, one also says that the total space $Z$ is a \emph{Mori fiber space}.
Mori fibrations $f:Z\ra B$ arise from the contractions of \emph{extremal rays of fiber type}
in the cone of curves $\overline{\operatorname{NE}}(Z)$, and play a central role in the minimal model
program.

A local noetherian ring  $R$, such as $R=\O_{Z,a}$, is called a \emph{klt singularity} if   
it is  normal, $\QQ$-Gorenstein, and 
for every proper birational modification  $r:X\ra  \Spec(R)$ with $X$ normal, the \emph{discrepancies} $\mu_i\in\QQ$ defined by 
$$
K_X= r^*(K_R) + \sum \mu_i E_i
$$
are bounded by $\mu_i> -1$. Here $E_i\subset X$ are the exceptional divisors,
and the equality holds in $\Pic(X)\otimes\QQ$.
If the discrepancies satisfy  the stronger condition $\mu_i\ge 0$, the local ring $R$ is a  \emph{canonical singularity}.
For $\mu_i>0$, the local ring is a \emph{terminal singularity}.
For further details on the singularities appearing in the context of the minimal model program, see \cite{Kollar 2013}, Section 2.

For every   fibration $f:Z\ra B$, the generic fiber $V=Z_\eta$ is a proper normal scheme over the function field
$F=\O_{B,\eta}=\kappa(\eta)=k(B)$, 
with $h^0(\O_V)=1$. It has dimension $n=\dim(Z)-\dim(B)$.
In Mori fiber spaces, it is $\QQ$-factorial, with  Picard number $\rho(V)=1$. The latter holds because
the scheme $Z$ is $\QQ$-factorial with $\rho(Z)-\rho(B)=1$. Moreover, the $\QQ$-Cartier divisor $K_V$ is antiample.
If the non-Gorenstein locus of the total space  $Z$ does not dominate the base $B$, then the
$V$  is an $n$-dimensional  Fano variety.
For $n=2$, the generic fiber $V$ is a del Pezzo surface. 
Let us record the following general fact:

\begin{lemma}
Suppose that $Z$ is a proper normal integral scheme,
and $Z\ra B$ be a fibration of relative dimension $\dim(Z)-\dim(B)=2$.
Then the generic fiber $V=X_\eta$ is normal and geometrically irreducible
over the function field $F=k(B)$.
If $Z$ has only canonical or terminal singularities, then
$V$ is Gorenstein or regular, respectively.
The scheme $V$ is geometrically reduced, provided that $\dim(B)=1$.
\end{lemma}

\proof
According to the arithmetical proof in \cite{Alexeev 1992}, each two-dimensional  canonical singularity is Gorenstein,
and each two-dimensional terminal singularity is regular.
The last statement follows from \cite{Schroeer 2010}, Theorem 2.3.
\qed

\medskip
Combining this with  Theorem \ref{main pdeg}, we get the following immediate consequence,
which answers the questions originating in Koll\'ar's study of extremal rays on threefolds
(\cite{Kollar 1991}, Remark 1.2):

\begin{theorem}
\mylabel{main fibration}
Let $Z$ be a threefold with only terminal singularities, and $f:Z\ra B$ be  a Mori fibration of relative dimension $n=2$.
Then the generic fiber $V=X_\eta$ is a  del Pezzo surface over the function field $F=k(B)$ that is geometrically normal.
\end{theorem}

\appendix
\section{Finite modifications and Gorenstein conditions}
\mylabel{Finite modifications}

In this appendix, we collect and slightly generalize some well-known facts pertaining to finite modifications, gluings of schemes,
and Gorenstein conditions that go back to Serre \cite{Serre 1988}, and were further developed 
by Reid \cite{Reid 1994}.

Suppose that $Y$ is a noetherian scheme, and let 
$\nu:X\ra Y$ be a \emph{finite modification}. This means that the morphism 
is finite and schematically dominant, and the  injective homomorphism $\O_Y\ra\nu_*(\O_X)$ is  
generically bijective.
The \emph{conductor ideal} $\mathfrak{C}\subset\O_Y$ is the annihilator   of the
coherent sheaf $\nu_*(\O_X)/\O_Y$. 
This is a quasicoherent ideal sheaf, and the corresponding closed subscheme $C\subset Y$ is called the \emph{conductor scheme}.
The  subsheaf $\mathfrak{C}\subset\nu_*(\O_X)$
is an ideal sheaf, in fact the largest ideal sheaf in $\O_Y$ that is also an ideal sheaf in $\nu_*(\O_X)$.
We call the resulting closed subscheme $R\subset X$ the \emph{ramification locus}.
The diagram
\begin{equation}
\label{conductor square}
\begin{CD}
R	@>>>	X\\
@VVV		@VV\nu V\\
C	@>>>	Y
\end{CD}
\end{equation}
is both cartesian and cocartesian.
Note that the finite morphisms $X\ra Y$ and the induced morphism $R\ra C$ are   schematically dominant.
We refer to the restriction $R\ra C$ as the \emph{gluing map},
and call the above diagram the \emph{conductor square}.
The sequence
\begin{equation}
\label{exact coherent}
0\lra\O_Y\lra\nu_*(\O_X)\oplus\O_C\lra\nu_*(\O_R)\lra 0
\end{equation}
of coherent sheaves is exact, where the arrow on the left is the diagonal map, and the arrow on the right
is the difference map  $(s,t)\mapsto s_R-t_R$. 
Using the Snake Lemma, one infers  that   the sequence 
\begin{equation}
\label{exact coherent short}
0\lra\O_Y\lra\nu_*(\O_X)\lra\nu_*(\O_R)/\O_C\lra 0
\end{equation}
is exact as well.
As explained in \cite{Schroeer; Siebert 2006}, Proposition 4.1 the exact sequence \eqref{exact coherent} of coherent sheaves induces an exact sequence
of multiplicative abelian sheaves
\begin{equation}
\label{exact multiplicative}
1\lra\O_Y^\times\lra\nu_*(\O_X)^\times\oplus\O_C^\times\lra\nu_*(\O_R)^\times\lra 1.
\end{equation}
where the map on the right is $(s,t)\mapsto s_R/t_R$.

Conversely, suppose we start with  a  noetherian scheme $X$,
a closed subscheme $R\subset X$ that contains no generic point $\eta\in X$,
and a   schematically dominant finite morphism $\nu:R\ra C$.
Then there exists a morphism $\nu:X\ra Y$ onto a noetherian \emph{algebraic space} $Y$,
making the diagram \eqref{conductor square}  cocartesian
(\cite{Artin 1970}, Theorem 6.1). The sequence \eqref{exact coherent} is exact,
and the diagram is also cartesian.
After replacing $C\subset Y$ and $R\subset X$, we may assume that these subscheme are the conductor scheme
and the ramification locus.
One also says that $Y$ is obtained by $X$ by \emph{gluing} with respect to the gluing map $R\ra C$.
For a general discussion of this process, we refer to Ferrand's paper \cite{Ferrand 2003}.

In this paper, we are mainly interested in the case where $X\ra Y$, or equivalently the gluing map 
$R\ra C$, are universal homeomorphisms (Proposition \ref{universal homeomorphism}). In this case, the algebraic space $Y$ 
is automatically a scheme, as follows from \cite{Olsson 2016}, Theorem 6.2.2. To simplify the exposition, we have decided to 
restrict our discussion to schemes. 
In what follows, we often write $\O_X$ and $\O_R$ instead of the more precise $\nu_*(\O_X)$ and $\nu_*(\O_R)$,
which should not cause any confusion.

If the the ramification locus $R$ has no embedded components, it is customary to  write $\shM_R$ for the quasicoherent
sheaf of meromorphic functions on $R$. In other words,  $\shM_R=i_*(\O_{R^{(0)}})$
where $i:R^{(0)}=\amalg_{\zeta\in R}\Spec(\O_{R,\zeta})\ra R$
is the inclusion of the scheme of generic points.
The gluing map induces an inclusion $\O_C\subset\shM_C\cap\O_R$, where the intersection takes
place in $\shM_R$. The following is a   generalization of
Reid's observation \cite{Reid 1994}, Proposition 2.2:

\begin{proposition}
\mylabel{condition S2}
Suppose that both $X$ and $Y$ satisfies Serre's Condition $(S_2)$. Then the following three conditions hold:
\begin{enumerate}
\item The ramification locus $R$ has no embedded components.
\item For each generic point $\zeta\in R$, the local ring $\O_{X,\zeta}$ has dimension one.
\item For each  point $y\in C$ with $\dim(\O_{Y,y})\geq 2$,  the inclusion
$\O_{C,y}\subset(\shM_C\cap\O_R)_y$ is an equality.
\end{enumerate}
\end{proposition}

\proof
For the first assertion, suppose there is an   associated point  $x\in\Ass(\O_R)$ that is embedded. Let $y\in Y$ be its image.
Clearly, $\dim(\O_{X,x})=\dim(\O_{Y,y})\geq 2$.
Replacing $Y$ by the spectrum of the   ring $\O_{Y,y}$,
we may assume that $y\in Y$ is a closed point.
Applying local cohomology to the short exact sequence \eqref{exact coherent} yields an
exact sequence of groups
$$
H^0_y(\O_Y)\lra H^0_y(\O_X)\oplus  H^0_y(\O_C)\lra H^0_y(\O_R)\lra H^1_y(\O_Y).
$$
The outer terms and $H^0_y(\O_X)$ vanish, because both $X$ and $Y$ satisfy $(S_2)$.
It follows that $H^0_y(\O_C)\ra H^0_y(\O_R)$ is bijective.
Now choose some non-zero $\bar{s}\in\O_{R,y}$ with support $\Supp(\bar{s})=\{x\}$. Then $\bar{s}\in  H^0_y(\O_R)=H^0_y(\O_C)$,
in particular $\bar{s}\in\O_{C,y}$.
Let $s\in\O_{X,y}$ be a representative.
For every $f\in \O_{X,y}$, the short exact sequence \eqref{exact coherent short} implies $sf\in\O_{Y,y}$, thus 
$s\in\mathfrak{C}$ and finally $\bar{s}=0$, contradiction.

For the second assertion, suppose that $\zeta\in R$ is a generic point, with image $y\in C$. 
Again we may assume that $y\in Y$ is closed.
Since $\nu:X\ra Y$ is birational,
we have $\dim(\O_{X,x})\geq 1$. In the above exact sequence, $H^0_y(\O_X)=0$ and
the inclusion $H^0_y(\O_C)=\O_{C,y}\subset \O_{R,y}=H^0_y(\O_R)$ is not an equality.
Therefore $H^1_y(\O_Y)\neq 0$, which implies that $\dim(\O_{Y,y})=1$.

Finally, we check condition (iii). 
Again  we may assume that $Y$ is local, with closed point  $y\in Y$.
Suppose that  $s\in \O_{R,y}$ 
lies in $\O_{C}$ at all generic points  $\zeta\in R$. Choose an open dense subset
$V\subset Y$ so that $s\in H^0(V,\O_C)$, and write $A=Y\smallsetminus V$ for the complementary closed
subset. 
Furthermore,  we have a commutative diagram  
$$
\begin{CD}
0	@>>>	H^0(C,\O_C)	@>>> 	H^0(V,\O_C)	@>>>	H^1_A(\O_C)	@>>>	0\\
@.		@VVV			@VVV			@VVV\\
0	@>>>	H^0(R,\O_R)	@>>> 	H^0(V,\O_R)	@>>>	H^1_A(\O_R)	@>>>	0.\\
\end{CD}
$$
The rows are exact, because the schemes $C$ and $R$ are affine. Since   $H^1_A(\O_Y)=H^1_A(\O_X)=0$,
the long exact sequence coming from  \eqref{exact coherent} ensures that  the vertical map on the right is injective. Since $C$ is  affine,
the vertical maps on the left has $H^0(C,\shF)$ as cokernel, with $\shF=\O_R/\O_C$. The cokernel
of the vertical map in the middle at least is contained in $H^0(V,\shF)$.
In turn, the Snake Lemma implies that
the restriction map  $H^0(C,\shF)\ra H^0(V,\shF)$ is injective. Thus the class of $s$ in $\shF$ vanishes,
in other words $s\in\O_{C,y}$.
\qed

\medskip
Next, we discuss Gorenstein conditions.
Let $B$ be a complete local noetherian ring.
Then the contravariant functor 
$$
M\longmapsto \Hom_B(H^d_\maxid(M), E)
$$
on the category of $B$-modules is representable. Here 
$\maxid=\maxid_B$ is the maximal ideal, $d=\dim(B)$ is the Krull dimension, 
and $E$ is the injective hull of the residue field $\kappa(B)$.
According to  Aoyama \cite{Aoyama 1983}, any $\omega_B$ representing this functor
is called a \emph{canonical module}. Note that   $\omega_B$ is unique up to isomorphism,
and that we do not demand that $B$ is Cohen--Macaulay.
The complete local ring  $B$ is called \emph{quasi-Gorenstein} if
the canonical module $\omega_B$ is invertible. If $B$ is additionally
Cohen--Macaulay, one says that $B$ is \emph{Gorenstein}.

For local noetherian rings $B$ that are not necessarily complete,
a module $\omega_B$ is called canonical if   $\omega_{\hat{B}}\simeq \omega_B\otimes_B\hat{B}$.
We then say that $B$ \emph{admits a canonical module}. Note that this is not always the case.
However, this condition obviously holds if $ \hat{B}$ is quasi-Gorenstein.
Extending the notions from the complete local case to the local case, 
we say that $B$ is Gorenstein or quasi-Gorenstein if 
the respective property holds for the completion $\hat{B}$.

We say that our noetherian  scheme $Y$ is Gorenstein or quasi-Gorenstein at a point $y\in Y$,
if the respective property holds for the local ring $B=\O_{Y,y}$, or equivalently
the complete local ring $\hat{B}=\O_{Y,y}^\wedge$.
We say that $Y$ is \emph{Gorenstein or quasi-Gorenstein} if 
the respective property holds for each point $y\in Y$.
If this holds for   all points with $\dim(\O_{Y,y})\leq n$,
we say that $Y$ is \emph{Gorenstein or quasi-Gorenstein in codimension $n$},
or that $Y$ satisfies condition $(G_n)$ or $(qG_n)$, respectively.

If the point $y\in Y$ of codimension one is contained  in the conductor scheme, 
$\O_{C,y}\subset\O_{R,y}$ is a finite extension of Artin rings,
and we denote by
$$
\length(\O_{C,y})\leq \length(\O_{R,y})
$$
their lengths  as modules over $\O_{C,y}$.
The following two results, under various additional assumptions, are due
to Samuel (\cite{Samuel 1951}, Theorem 5), Gorenstein (\cite{Gorenstein 1952}, Theorem 6), 
Rosenlicht (\cite{Rosenlicht 1952}, Theorem 14), 
Roquette \cite{Roquette 1961},  Serre 
(\cite{Serre 1988}, Chapter IV, \S 3, Section 11),
Kunz \cite{Kunz 1970}
and Reid (\cite{Reid 1994}, Theorem 3.2). 
In our general form,  we merely demand suitable Gorenstein assumptions 
and make no reference to ground rings:

\begin{proposition}
\mylabel{necessary condition gorenstein}
Suppose that  $Y$ is  Gorenstein in codimension one.
Then the length formula
$$
\length(\O_{R,y})=2\length(\O_{C,y})
$$
holds for each point $y\in Y$ of codimension one contained in $C$.
\end{proposition}

\proof
It suffices to treat the case that $Y=\Spec(B)$ is the spectrum of a  complete one-dimensional local noetherian ring.
Set $l=\length(B/\mathfrak{C})$ and  choose a sequence of ideals 
\begin{equation}
\label{first sequence g}
\mathfrak{C}  = \mathfrak{b}_0\subset  \mathfrak{b}_1\subset \ldots\subset \mathfrak{b}_l=B
\end{equation}
with simple sub-quotients. Note that this corresponds to
a Jordan--H\"older sequence for $B/\mathfrak{C}$, which is a module of finite length.
Clearly,  each $\mathfrak{b}_i$ is a maximal Cohen--Macaulay
$B$-module. According to \cite{Eisenbud 1995}, Theorem 21.21,
the  contravariant functor $M\mapsto \Hom_B(M,\omega_B)$ is an exact antiequivalence 
from the category of maximal Cohen--Macaulay modules to itself. In
fact, the biduality maps 
$$
M\lra \Hom_B(\Hom_B(M,\omega_B),\omega_B),\quad m\longmapsto (f\mapsto f(m))
$$
are bijective. Since    $B$ is Gorenstein, the module   $\omega_B=B$ is canonical,
and we conclude that \eqref{first sequence g} induces another sequence
\begin{equation}
\label{second sequence g}
\Hom_B(\mathfrak{b}_l,B)\subset\Hom_B(\mathfrak{b}_{l-1}, B)\subset\ldots\subset \Hom_B(\mathfrak{b}_0, B) 
\end{equation}
of the same length, with simply sub-quotients. 
The term on the left is $\Hom_B(B,B)= B$. Thus we may splice the
two sequences \eqref{first sequence g} and \eqref{second sequence g}
and obtain a sequence of length $2l$ with simple sub-quotients,
starting with $\mathfrak{C}=\mathfrak{b}_0$ and ending with $\Hom_B(\mathfrak{C},B)=\Hom_B(\mathfrak{b}_0, B)$.
It remains to identify the latter with $A$, in such a way that the resulting inclusion $\mathfrak{C}\subset A$
coincides with the canonical inclusion. To achieve this, consider the commutative diagram
$$
\begin{CD}
B		@>\can>>		A\\
@VVV				@VVV\\
\Hom_B(B,B)	@>>\res>	\Hom_B(\mathfrak{C},B),
\end{CD}
$$
where the upper map is the canonical inclusion, the lower map 
is given by restriction, 
and the vertical maps are given by $a\mapsto (x\mapsto ax)$.
We need to verify that the map $A\ra \Hom_B(\mathfrak{C},B)$ is bijective.
Applying the antiequivalence of categories again and using biduality, we have to check that
$$
\mathfrak{C}=\Hom_B(\Hom_B(\mathfrak{C},B),B) \lra \Hom_B(A,B),\quad
b\longmapsto (y\mapsto by)
$$
is bijective. Clearly,  the   map $f\mapsto f(1)$ 
is a left inverse. Thus $\mathfrak{C}\ra\Hom_B(A,B)$ is injective and admits a complement $M\subset\Hom_B(A,B)$.
Since $B$ is Cohen--Macaulay, thus torsion-free, the Hom module and thus also $M$ are torsion-free.
But the canonical inclusions $\mathfrak{C}\subset B\subset A$ become equalities after inverting
any regular $s\in\maxid_B$. It follows that $M_s=0$, and thus $M=0$.
Summing up, we have shown that $\length(A/\mathfrak{C})=2l$.  
\qed

\medskip
The converse statement takes the following form: 

\begin{proposition}
\mylabel{sufficient condition gorenstein}
Let $y\in Y$ be a point of codimension one. Suppose  that the following three conditions hold:
\begin{enumerate}
\item The semilocal ring $\O_{X,y}$ is Gorenstein
\item The conductor ideal $\mathfrak{C}_y\subset\O_{X,y}$ is invertible.
\item The  modules $\O_{C,y}$ and
$(\O_R/\O_C)_y$ are free of the same rank over some local Artin subring $W\subset\O_{C,y}$ with $W$ Gorenstein.
\end{enumerate} 
Then the local ring $\O_{Y,y}$ is Gorenstein.
\end{proposition}

\proof
To simplify notation,   set $A= \O_{X,y}$ and $B=\O_{Y,y}$. It suffices to treat the
case that the rings in question are complete, whence $A$ is a product of complete
local rings.  
By assumption, the  conductor ideal has the form  $fA=\mathfrak{C}_y$ for some regular element $f\in A$,
and we write   $\bar{A}, \bar{B}$ for the resulting residue class rings.
Now the exact sequence \eqref{exact coherent short} takes the form
\begin{equation}
\label{cartesian sequence}
0\lra B\lra A\lra \bar{A}/\bar{B}\lra 0.
\end{equation}
Since the semilocal ring $A$ and the   subring $W\subset\bar{B}$ are  Gorenstein, 
the modules $\omega_A=A$ and $\omega_W=W$
are dualizing. In turn,   both modules 
$$
\Ext^1_A(\bar{A},A)=\Ext^1_A(\bar{A},\omega_A)\quadand 
\Hom_W(\bar{A},W)=\Hom_W(\bar{A},\omega_W)
$$
are  dualizing for the semilocal ring $\bar{A}$, and it follows that these $\bar{A}$-modules are isomorphic. Likewise,
$\omega_{\bar{B}}=\Hom_W(\bar{B},W)$ is dualizing for $\bar{B}$.
Using the short exact sequence $0\ra fA\ra A\ra \bar{A}\ra 0$,
we get an exact sequence
$$
0\lra \Hom(A,\omega_A)\lra  \Hom(fA,\omega_A) \stackrel{\partial}{\lra} \Ext^1(\bar{A},\omega_A) \lra 0.
$$
According to \cite{Reid 1994}, the theorem in 2.6,  the kernel of the composite map
\begin{equation}
\label{reid's sequence}
\Hom(fA,\omega_A)\stackrel{\partial}{\lra} \Ext^1(\bar{A},\omega_A) 
\simeq
\Hom_W(\bar{A},\omega_W)\stackrel{\res}{\lra} \Hom_W(\bar{B},W)
\end{equation}
is a dualizing module $\omega_B$ for $B$. 
The arrow on the left is the connecting map, and the map on the right the restriction map, which is
often referred to as the \emph{trace map}.
Up to isomorphism, the  kernel does not depend on the
chosen isomorphism in the middle, because the connecting map  is isomorphic
to the residue class map  $A\ra\bar{A}$, and the induced homomorphism  of multiplicative groups $A^\times\ra\bar{A}^\times$ 
for the semilocal rings is surjective. The idea now is to  choose a particular
isomorphism in the middle of \eqref{reid's sequence}   that is adapted to our problem.

By assumption, the finitely generated $W$-modules $\bar{B}$ and $\bar{A}/\bar{B}$ are free.
Hence the inclusion $\bar{B}\subset\bar{A}$ admits a complement,
and we may write $\bar{A}=\bar{B}\oplus\bar{A}/\bar{B}$ as $W$-modules.
Choosing a $W$-basis in these summands, 
we obtain   a unimodular symplectic  form  $\Phi:\bar{A}\times\bar{A}\ra W$ so that $\bar{A}/\bar{B}\subset\bar{A}$ becomes
a Lagrangian, with orthogonal complement $\bar{B}\subset\bar{A}$. This   crucial step  
hinges on the assumption that $\rank_W(\bar{A})=\rank_W(\bar{A}/\bar{B})$.
In turn, the diagram
\begin{equation}
\label{adjoint maps}
\begin{CD}
\bar{A} 	@>\Phi>>	\Hom_W(\bar{A},W)\\
@V\pr VV				@VV\res V\\
\bar{A}/\bar{B}	@>>\Phi>		\Hom_W(\bar{B},W)
\end{CD}
\end{equation}
becomes commutative, where the horizontal  maps are given by  $a\mapsto (a'\mapsto\Phi(a,a'))$,
the map on the left is the canonical projection, and the map on the right is given by restriction.

We now can make the sequence \eqref{reid's sequence} explicit:
Using the homomorphism $f\mapsto 1$ as an $A$-basis $e\in\Hom(fA,\omega_A)$ 
and its image under the connecting map as an $\bar{A}$-basis $\bar{e}\in \Ext^1(\bar{A},\omega_A)$, where   $\omega_A=A$,
we may regard $\omega_B$ as   the kernel for the composite map
$$
Ae \stackrel{}{\lra} \bar{A}\bar{e}
=
\bar{A} \stackrel{\Phi}{\lra} \Hom_W(\bar{A},W)\stackrel{\res}{\lra} \Hom_W(\bar{B},W).
$$
By diagram \eqref{adjoint maps},  this composite mapping is isomorphic to 
the canonical projection $A\ra\bar{A}/\bar{B}$. The exact sequence
\eqref{cartesian sequence} now gives the desired isomorphism $\omega_B\simeq B$.
\qed

\medskip
Note that if the local Artin  ring $\O_{C,y}$ contains a field, it is natural to choose  for $W\subset\O_{C,y}$
a coefficient field, and   Condition (iii) becomes equivalent
to the length condition $\length(\O_{R,y})=2\length(\O_{C,y})$.

\begin{proposition}
\mylabel{quasi-gorenstein}
Suppose that $X$ and $Y$ satisfies Condition $(S_2)$, and that
the three assumptions of Proposition \ref{sufficient condition gorenstein} hold for each point $y\in Y$ of codimension one.
Then $Y$ is quasi-Gorenstein   if for each closed point $x\in X$, the classes of $\omega_A$ and $R$
in the local Class group $\Cl(A)$ for the local ring $A=\O_{X,x}$ are inverse to each other.
\end{proposition}

\proof
According to \cite{Aoyama 1983}, Corollary 2.4 it suffices to treat the case
that  $Y=\Spec(B)$ is a complete local scheme, with closed point $y\in Y$.
The canonical module $\omega_B$ and $B$  satisfies Serre's Condition $(S_2)$ 
(see \cite{Aoyama 1983} 1.10).
Using the arguments for \cite{Hartshorne 1994}, Theorem 1.12,
it suffices to check that the coherent sheaves $\shF=\widetilde{\omega_B}$
and $\O_Y$ are isomorphic over some open subset $V\subset Y$ that contains all points of codimension one.
Furthermore, by the assumption on the local class group
we may choose as canonical module $\omega_A=\mathfrak{C}$, such that we get
an identification $A=\Hom(\mathfrak{C},\omega_A)$.
As in the proof of Proposition \ref{sufficient condition gorenstein}, we consider the composite map 
\begin{equation}
\label{adjunction}
A=\Hom_A(\mathfrak{C},\omega_A)\lra\Ext^1_A(\bar{A},\omega_A)\simeq 
\Hom_{\bar{B}}(\bar{A},\omega_{\bar{B}})\lra \omega_{\bar{B}}.
\end{equation}
The map on the left is the connecting map  from   $0\ra\mathfrak{C}\ra A\ra \bar{A}\ra 0$,
whereas the map on the right is the trace map.
Note that the term on the right is supported by the conductor $C\subset Y$.
According to \cite{Reid 1994}, 2.6, the theorem,  the kernel $I$ of the above composition  is a canonical module $\omega_B$.

By assumption, both schemes $X$ and $Y$ satisfies Serre's Condition $(S_2)$. According
to Proposition \ref{condition S2}, the ramification locus $R\subset X$ is purely one-codimensional,
and contains no embedded components.  
Let $\zeta_1,\ldots,\zeta_r\in C$ be the generic points.
Then $\dim(\O_{X,\zeta_i})=1$, and we saw in the proof for Proposition \ref{sufficient condition gorenstein}
that the $\bar{B}$-module $\Ext^1(\bar{A},\omega_A)$ and
$\bar{A}/\bar{B}$ are isomorphic at these points $\zeta_i\in Y$.
It follows that there is an open neighborhood $V\subset Y$ containing all points of codimension one
so that   the quasicoherent sheaf attached to $\omega_{\bar{B}}$ and the sheaf $\O_R/\O_C$     have
isomorphic restriction to $V$. Let   $\shI=\widetilde{I}$ the   quasicoherent sheaf on $Y$
corresponding to the $B$-module $I=\omega_B$.
By our construction, $\shI|V\simeq\O_V$.
Thus we see that $\widetilde{\omega_B}|V\simeq\O_V$.
\qed


\end{document}